\newif\iffinal
\finalfalse	
\finaltrue	

\documentclass[letterpaper,12pt,reqno]{amsart}
\RequirePackage[utf8]{inputenc}
\usepackage[portrait,margin=2.95cm]{geometry}
\iffinal\else\usepackage[notref,notcite]{showkeys}\fi
\usepackage{mathrsfs,soul,mathtools} 
\usepackage{extarrows} 
\usepackage{hyperref}
\usepackage[foot]{amsaddr}
\usepackage{amssymb,amsthm,amsbsy,latexsym,amsfonts}
\usepackage[textsize=small]{todonotes}       
\usepackage{graphicx}
\usepackage[numeric,initials,nobysame]{amsrefs}
\usepackage{upref,setspace}
\usepackage{tikz}
\usetikzlibrary{calc,intersections,through,backgrounds,shapes.geometric}
\usetikzlibrary{graphs}

\usepackage{caption}
\usepackage{subcaption}
\usepackage{color}

\usepackage{fourier}

\usepackage{enumerate}
\newenvironment{enumeratei}{\begin{enumerate}[\upshape (i)]}{\end{enumerate}}
\newenvironment{enumeratea}{\begin{enumerate}[\upshape (a)]}{\end{enumerate}}

\newenvironment{enumerateA}{\begin{enumerate}[\upshape (A)]}{\end{enumerate}}

\usepackage{paralist}

\newenvironment{inparaenuma}{\begin{inparaenum}[\upshape \bfseries (a) ]}{\end{inparaenum}}

\newenvironment{inparaenumi}{\begin{inparaenum}[\upshape  (i) ]}{\end{inparaenum}}

\newenvironment{inparaenumibf}{\begin{inparaenum}[\upshape \bfseries (i) ]}{\end{inparaenum}}

\numberwithin{figure}{section}
\numberwithin{table}{section}

\newtheorem{thm}{Theorem}[section]
\newtheorem{theorem}{Theorem}[section]
\newtheorem{lem}[thm]{Lemma}

\newtheorem{prop}[thm]{Proposition}
\newtheorem{defn}[thm]{Definition}

\newtheorem{ass}[thm]{Assumption}
\newtheorem*{ass*}{Assumption}
\newtheorem*{theorem*}{Theorem}

\newtheorem{conj}[thm]{Conjecture}
\newtheorem{lemma}[thm]{Lemma}

\newtheorem{construction}[thm]{Construction}

\newtheorem*{thm3.9*}{Theorem 3.9*}

\theoremstyle{definition}
\newtheorem{rem}{Remark}
\newtheorem{obs}{Observation}
\numberwithin{obs}{section}

\newcommand{\ind}{\1}
\newcommand{\eps}{\varepsilon}

\newcommand{\set}[1]{\left\{#1\right\}}

\newcommand{\equald}{\stackrel{\mathrm{d}}{=}}
\newcommand{\probc}{\stackrel{\mathrm{P}}{\longrightarrow}}

\newcommand{\weakc}{\stackrel{\mathrm{d}}{\longrightarrow}}
\newcommand{\convas}{\stackrel{\mathrm{a.s.}}{\longrightarrow}}

\newcommand{\len}{\mathrm{len}}

\newcommand{\cb}{\mathrm{CB}}
\newcommand{\cbd}{\mathrm{CBD}}
\newcommand{\connects}[1]{\xleftrightarrow{#1}} 
\newcommand{\connectsn}{\xleftrightarrow{n^{\eta}/\lambda^{1+\pza}}} 
\def\qed{ \hfill $\blacksquare$}

\newcommand{\cmnd}{{\mathbb{CM}}^{\vd}}
\newcommand{\cgnd}{\mathbb{US}^{\vd}}

\newcommand*\xbar[1]{%
	\hbox{%
		\vbox{%
			\hrule height 0.5pt 
			\kern0.5ex
			\hbox{%
				\kern-0.5em
				\ensuremath{#1}%
				\kern-0.1em
			}%
		}%
	}%
}

\newcommand{\cA}{\mathcal{A}}\newcommand{\cB}{\mathcal{B}}\newcommand{\cC}{\mathcal{C}}
\newcommand{\cD}{\mathcal{D}}\newcommand{\cE}{\mathcal{E}}
\newcommand{\cG}{\mathcal{G}}\newcommand{\cH}{\mathcal{H}}
\newcommand{\cL}{\mathcal{L}}
\newcommand{\cM}{\mathcal{M}}\newcommand{\cN}{\mathcal{N}}

\newcommand{\cT}{\mathcal{T}}
\newcommand{\cV}{\mathcal{V}}\newcommand{\cW}{\mathcal{W}}
\newcommand{\cZ}{\mathcal{Z}}

\newcommand{\vd}{\mathbf{d}}

\newcommand{\vp}{\mathbf{p}}
\newcommand{\vt}{\mathbf{t}}

\newcommand{\mvD}{\boldsymbol{D}}
\newcommand{\mvG}{\boldsymbol{G}}

\newcommand{\mvM}{\boldsymbol{M}}

\newcommand{\mvS}{\boldsymbol{S}}\newcommand{\mvU}{\boldsymbol{U}}

\newcommand{\mvp}{\boldsymbol{p}}
\newcommand{\mvq}{\boldsymbol{q}}

\newcommand{\mvw}{\boldsymbol{w}}\newcommand{\mvx}{\boldsymbol{x}}\newcommand{\mvy}{\boldsymbol{y}}

\newcommand{\mvtheta}{\boldsymbol{\theta}}

\newcommand{\mvmu}{\boldsymbol{\mu}}

\newcommand{\fA}{\mathfrak{A}}\newcommand{\fB}{\mathfrak{B}}\newcommand{\fC}{\mathfrak{C}}
\newcommand{\fD}{\mathfrak{D}}\newcommand{\fE}{\mathfrak{E}}\newcommand{\fF}{\mathfrak{F}}

\newcommand{\fM}{\mathfrak{M}}

\newcommand{\fS}{\mathfrak{S}}\newcommand{\fT}{\mathfrak{T}}

\newcommand{\bE}{\mathbb{E}}
\newcommand{\bG}{\mathbb{G}}

\newcommand{\bP}{\mathbb{P}}\newcommand{\bR}{\mathbb{R}}
\newcommand{\bT}{\mathbb{T}}

\newcommand{\bZ}{\mathbb{Z}}

\newcommand{\dI}{\mathpzc{I}}

\DeclareFontFamily{OT1}{pzc}{}
\DeclareFontShape{OT1}{pzc}{m}{it}{<-> s * [1.10] pzcmi7t}{}
\DeclareMathAlphabet{\mathpzc}{OT1}{pzc}{m}{it}

\newcommand{\pzI}{\mathpzc{I}}

\newcommand{\pzS}{\mathpzc{S}}

\newcommand{\pzZ}{\mathpzc{Z}}

\newcommand{\pza}{\mathpzc{h}} 

\newcommand{\pzh}{\mathpzc{h}}

\newcommand{\pzt}{\mathpzc{t}}

\newcommand{\rH}{\mathrm{H}}

\newcommand{\rP}{\mathrm{P}}

\newcommand{\rr}{\mathrm{r}}

\newcommand{\sM}{\mathscr{M}}

\newcommand{\sW}{\mathcal{W}}

\newcommand{\sP}{\mathpzc{P}}

\DeclareMathOperator{\E}{\mathbb{E}}
\DeclareMathOperator{\pr}{\mathbb{P}}

\DeclareMathOperator{\GH}{GH}
\DeclareMathOperator{\GHP}{GHP}
\DeclareMathOperator{\diam}{diam}
\DeclareMathOperator{\height}{ht}
\DeclareMathOperator{\ord}{ord}
\DeclareMathOperator{\dist}{dist}

\DeclareMathOperator{\MC}{MC}
\DeclareMathOperator{\conn}{Conn}

\DeclareMathOperator{\core}{Core}
\DeclareMathOperator{\avail}{avail}

\DeclareMathOperator{\con}{con}

\DeclareMathOperator{\udim}{\underline{dim}}
\DeclareMathOperator{\odim}{\overline{dim}}
\DeclareMathOperator{\res}{res}

\DeclareMathOperator{\surplus}{sp}
\DeclareMathOperator{\spls}{sp}
\DeclareMathOperator{\partition}{ptn}

\newcommand{\sss}{\scriptscriptstyle}
\newcommand{\erdos}{Erd\H{o}s-R\'enyi }

\newcommand{\convd}{\stackrel{d}{\longrightarrow}}

\definecolor{aqua}{rgb}{0.0, 1.0, 1.0}

\def\1{{\mathchoice {1\mskip-4mu\mathrm l}      
		{1\mskip-4mu\mathrm l}
		{1\mskip-4.5mu\mathrm l} {1\mskip-5mu\mathrm l}}}

\newcommand{\e}{{\mathrm{e}}}

\usepackage{accents}
\newcommand\thickbar[1]{\accentset{\rule{.5em}{1pt}}{#1}}
\newcommand\tbar[1]{\accentset{\rule{.5em}{1pt}}{#1}}

\newcommand{\tbarGn}{\thickbar{\mvG}_n}
\newcommand{\barGn}{\tbar{G}_n}

\newcommand{\ch}[1]{\textcolor{black}{{#1}}}
\newcommand{\op}{o_{\sss \mathrm{P}}}

\newcommand{\falln}{\thickbar{\forall}n}
\newcommand{\fallnt}{{\ooalign{\hfil$\forall$\hfil\cr\kern.08em-\hfil}}}

\newcommand{\mst}{\mathscr{M}}

\usepackage{frcursive}
\newenvironment{frcseries}{\fontfamily{frc}\selectfont}{}
\newcommand{\textfrc}[1]{{\frcseries#1}}
\newcommand{\mathfrc}[1]{\text{\textfrc{#1}}}

\newcommand{\cs}{\mathfrc{s}}

\newcommand{\crv}{v}

\newcommand{\Poi}{\mathrm{Poi}}
\newcommand{\Unif}{\mathrm{Unif}}
\newcommand{\Bern}{\mathrm{Bernoulli}}
\newcommand{\EXP}{\mathrm{Exp}}

\newcommand{\stod}{\preceq_{\sf sd}}

\newcommand{\pmtr}{{\mathpzc{Pmtr}}}

\newcommand{\mtbp}{\mathrm{MTBP}}
\newcommand{\tlbp}{\mathrm{L}3\mathrm{BP}}
\newcommand{\twolbp}{\mathrm{L}2\mathrm{BP}}

\definecolor{aqua}{rgb}{0.0, 1.0, 1.0}
\definecolor{webbrown}{rgb}{.6,0,0}
\definecolor{pinegreen}{rgb}{0.0, 0.47, 0.44}
\definecolor{ultramarineblue}{rgb}{0.25, 0.4, 0.96}
\definecolor{jrnl}{rgb}{0.0, 0.5, 0.0}
\definecolor{lincolngreen}{rgb}{0.11, 0.35, 0.02}
\definecolor{green(html/cssgreen)}{rgb}{0.0, 0.5, 0.0}
\definecolor{airforceblue}{rgb}{0.36, 0.54, 0.66}
\definecolor{azure}{rgb}{0.0, 0.5, 1.0}
\definecolor{bleudefrance}{rgb}{0.19, 0.55, 0.91}
\definecolor{cobalt}{rgb}{0.0, 0.28, 0.67}

\hypersetup{colorlinks=true, linktocpage=true, pdfstartpage=1, pdfstartview=FitV,
	breaklinks=true, pdfpagemode=UseNone, pageanchor=true, pdfpagemode=UseOutlines,
	plainpages=false, bookmarksnumbered, bookmarksopen=true, bookmarksopenlevel=1,
	hypertexnames=true, pdfhighlight=/O,
	urlcolor=webbrown, linkcolor=cobalt, citecolor=jrnl
}

\usepackage[normalem]{ulem}

\begin{document}
	
\title[MST in the heavy-tailed regime]{
Geometry of the minimal spanning tree in the heavy-tailed regime: new universality classes}

	\date{}
	
	\author[Bhamidi]{Shankar Bhamidi$^1$}
	\address{\hskip-15pt $^1$Department of Statistics and Operations Research, University of North Carolina, Chapel Hill}
	\author[Sen]{Sanchayan Sen$^2$}
	\address{\hskip-15pt $^2$Department of Mathematics, Indian Institute of Science}
	\email{bhamidi@email.unc.edu, sanchayan.sen1@gmail.com}

	\maketitle
	\begin{abstract}
		A \ch{well-known open problem} on the behavior of optimal paths in random graphs in the strong disorder regime, formulated by statistical physicists, and supported by a large amount of numerical evidence over the last decade \cite{BraBulCohHavSta03,wu2006transport,braunstein2007optimal,chen2006universal} is as follows: for a large class of random graph models with degree exponent $\tau\in (3,4)$, distances in the minimal spanning tree (MST) on the giant component in the supercritical regime scale like $n^{(\tau-3)/(\tau-1)}$. 
	    \ch{The aim of this paper is to make progress towards a proof of this conjecture.}

		We consider a supercritical inhomogeneous random graph model with degree exponent $\tau\in(3, 4)$ that is closely related to Aldous's multiplicative coalescent, and show that the MST constructed by assigning i.i.d. continuous weights to the edges in its giant component, endowed with the tree distance scaled by $n^{-(\tau-3)/(\tau-1)}$, converges in distribution with respect to the Gromov-Hausdorff topology to a random compact real tree. 
		Further, almost surely, every point in this limiting space either has degree one (leaf), or two, or infinity (hub), both the set of leaves and the set of hubs are dense in this space, and the Minkowski dimension of this space equals $(\tau-1)/(\tau-3)$.

		The multiplicative coalescent, in an asymptotic sense, describes the evolution of the component sizes of various near-critical random graph processes.
		We expect the limiting spaces in this paper to be the candidates for the scaling limit of the MST constructed for a wide array of other heavy-tailed random graph models. 
		
	\end{abstract}
	
	
	\tableofcontents
	
\section{Introduction }\label{sec:intro}
Consider a finite, connected, and weighted graph $(V,E,b)$, where $(V,E)$ is the underlying graph and $b:E\to [0,\infty)$ is the weight function.
A spanning tree of $(V,E)$ is a tree that is a subgraph of $(V,E)$ with vertex set $V$.
A minimal spanning tree (MST) $T$ of $(V,E,b)$ satisfies
\begin{align}\label{eqn:def-mst}
\sum_{e\in T}b(e)=\min\bigg\{\sum_{e\in T'}b(e):\ T'\text{ is a spanning tree of }(V,E)\bigg\}.
\end{align}
The MST is one of the most studied functionals in combinatorial optimization. 
Studying this object when the edge weights and potentially the underlying graph is random has stimulated an enormous body of work in probabilistic combinatorics and geometric probability.
The papers  
\cite{beardwood-halton-hammersley, aldous-steele, alexanderI, avram-bertsimas, steele, kesten-lee, alexander, chatterjee-sen, penroseI, penroseII, penroseIII, frieze1985value, beveridge-frieze-mcdiarmid, frieze-mcdiarmid, frieze-ruszink-thoma, penrose1998random, aldous1990random, janson1995minimal, petegabor} 
and the references therein give a non-exhaustive account of the enormous literature on the probabilistic study of MSTs.

We are interested in the global geometric properties of the MST, e.g., the diameter and the typical distance.   
In the early 2000s, several major conjectures were made about the intrinsic geometry of MSTs in the statistical physics community where this object arises as models of disordered networks,
but until very recently, there were few rigorous mathematical results on this problem.
Consider a finite connected graph $G=(V, E)$, and assign costs $\exp(\beta\eps_e)$ to the edges $e\in E$ where $\beta>0$.
An optimal path $P$ in $G$ between $u, v\in V$ minimizes the total cost $\sum_{e\in P'}\exp(\beta\eps_e)$ among all paths $P'$ between $u$ and $v$. 
This model interpolates between the first passage percolation regime (the weak disorder regime) and the minimal spanning tree regime (strong disorder regime).
Assuming that $\eps_e$, $e\in E$, are pairwise distinct, it is easy to see that for sufficiently large values of $\beta$, the optimal path between any two vertices $u$ and $v$ is the path $P$ that minimizes the maximal edge weight $\max_{e\in P'}\eps_e$ among all paths $P'$ connecting $u$ and $v$ in $G$.
It is well-known (see Lemma \ref{lem:mst-minimax-criterion}) that this is simply the path connecting $u$ and $v$ in the MST of $G$ constructed using the weights $\eps_e$, $e\in E$.
Thus, the number of edges or the hopcount $l_{opt}$ in the optimal path between two typical vertices in the presence of strong disorder is simply the length of the path in the MST connecting two typical vertices.

Motivated by the availability of data on a host of real-world networks as well as the impact of complex networks in our daily lives, the last few years have witnessed an explosion in the formulation and study of mathematical models of networks. 
These models try to capture properties of networks observed empirically. 
Of relevance to us is the heavy-tailed nature of the empirical degree distribution. 
To study this phenomenon, a plethora of random graph models have been proposed that have heavy-tailed degree distributions with some degree exponent $\tau \in (1,\infty]$. 
A precise definition of the degree exponent will not be needed in the sequel, so we instead refer the reader to \cite{Hofs15, Hofs17, durrett-rg-book} for a detailed discussion on the random graph models now available to practitioners.

Coming back to disordered networks, in the 2000s, statistical physicists predicted  \cite{BraBulCohHavSta03,wu2006transport,braunstein2007optimal,chen2006universal} 
that if the underlying graph posseses a heavy-tailed degree distribution with exponent $\tau$, then in the presence of strong disorder, $l_{opt}$ exhibits the following scaling behavior:
\begin{align}\label{eqn:0}
l_{opt}
\sim
\left\{
\begin{array}{l}
n^{1/3}\, ,\ \text{ if }\tau>4,\\
n^{(\tau-3)/(\tau-1)}\, ,\ \text{ if }\tau\in(3, 4)\, ,
\end{array}
\right.
\end{align}
where $n$ denotes the number of vertices in the underlying graph, and further, such behavior should be universal, i.e., in principle should apply to a wide array of random graph models.

The above conjecture is related to the universality of the intrinsic geometry of the MST, proving which in full generality remains open to date.
Only recently, there has been some progress in the $\tau>4$ regime.
In \cite{AddBroGolMie13}, it was shown that the MST of the complete graph on $n$ vertices constructed using i.i.d. $\Unif[0, 1]$ edge weights, endowed with the tree distance scaled by $n^{-1/3}$ and the uniform probability measure on the vertices, converges in distribution to a random compact $\bR$-tree $\cM$.
Further, almost surely, $\cM$ is binary and the Minkowski dimension of $\cM$ equals $3$.
This result can be seen as a confirmation of \eqref{eqn:0} for the case $\tau=\infty$.
The limiting space $\cM$ is expected to be the scaling limit of the MST of a variety of random discrete structures including random graphs with degree exponent $\tau>4$.
A first step in this broader program of establishing universality of $\cM$ was taken in \cite{addarioberry-sen}, where it was shown that the scaling limit of the MST of the random $3$-regular (simple) graph as well as the $3$-regular configuration model, with the tree distance scaled by $n^{-1/3}$, is $6^{1/3}\cdot\cM$.

On the other hand, the problem has stayed completely open in the $\tau\in(3, 4)$ regime.
The aim of this work is to study the MST in this regime and prove the existence of the scaling limit of the MST.
We will consider a \ch{supercritical} inhomogeneous random graph (IRG) model that corresponds to the rank-$1$ case of the general class of IRGs studied in \cite{BJR07}.
We will show that under certain assumptions, the MST on the giant component of this random graph, endowed with the tree distance scaled by $n^{-(\tau-3)/(\tau-1)}$, converges in distribution with respect to the Gromov-Hausdorff topology to a random compact $\bR$-tree.
We will then study the topological properties of the scaling limit.
We will show that almost surely, every point in this limiting space either has degree one (leaf), or two, or infinity (hub), and that both the set of leaves and the set of hubs are dense in this space.
Further, the Minkowski dimension of this space equals $(\tau-1)/(\tau-3)$ almost surely.
Note the contrasting characteristics of such a space in comparison to the space $\cM$.

The inhomogeneous random graph model considered in this paper is of special interest as it is closely related to the multiplicative coalescent \cite{aldous-crit, aldous-limic}.
The evolution of a large class of dynamic random graph models around the point of phase transition can be well-approximated \cite{SBSSXW-universal, SB-SD-vdH-SS} by the multiplicative coalescent, and we expect that this fact can be leveraged to establish universality of the scaling limits obtained in this paper. 
We defer further discussion regarding the general program to establish universality to Section~\ref{sec:disc}.
Our main result proves convergence of the MST viewed as a metric space,
\ch{which in particular implies distributional convergence of the diameter of the MST rescaled by $n^{-(\tau-3)/(\tau-1)}$.
This presents some rigorous evidence supporting the prediction in \eqref{eqn:0}, although the asymptotic behavior of the typical distance claimed in \eqref{eqn:0} does not follow from this alone.
The scaling limit of the typical distance can be deduced if one were able to establish distributional convergence of the MST viewed as a metric measure space.}
This strengthening of our result can be achieved if one additional estimate is proved.
This will also be discussed in Section~\ref{sec:disc}.

\subsection{Organization of the paper}
A reasonable amount of notation regarding notions of convergence of metric space-valued random variables as well as $\bR$-trees is required. 
To quickly get to the main result, we first define the random graph models of interest in Section \ref{sec:models}, and then describe the main result in Section \ref{sec:res}. 
Various definitions and preparatory results are then given in Section \ref{sec:defn}. 
In Section \ref{sec:different-model}, we show that it is equivalent to work with a modified random graph model; this will make many of our calculations easier.
A key step in the proof of our main result is relating the MST in the supercritical graph to the MST in the critical window; this is accomplished in Section \ref{sec:7}.
Then we complete the proof of the main result in Section \ref{sec:proof-main-thm}.
In Section \ref{sec:disc} we discuss possible extensions and how the scaling limit obtained here can be shown to be universal.
We also discuss how our main result can be strengthened to convergence with respect to the Gromov-Haussdorff-Prokhorov topology.
A few technical proofs are relegated to Appendix~\ref{sec:appendix}.

\section{Random graph models}\label{sec:models}
We start by describing the class of random graph models of interest for this paper. 
Start with the vertex set $[n]:=\set{1,2,\ldots, n}$, and suppose each vertex $i\in [n]$ has a weight $w_i^{\sss(n)}\geq 0$ attached to it; intuitively this measures the propensity or attractiveness of this vertex in the formation of links. We assume that the vertices are labeled so that $w_1^{\sss(n)}\geq w_2^{\sss(n)} \geq \cdots\geq w_n^{\sss(n)}$. 
Write $\mvw^{\sss(n)}=(w_1^{\sss(n)},\hdots, w_n^{\sss(n)})$, and let
\begin{equation}
\label{eqn:rank-one-connection}
	\thickbar{\mvq}_{ij}:= \thickbar{\mvq}_{ij}(\mvw^{\sss(n)}) = \frac{w_i^{\sss(n)} w_j^{\sss(n)}}{L_n} \wedge 1\, , \qquad 
	1\leq i \neq j\leq n\, , 
\end{equation}
where $L_n$ is the total weight given by
\begin{equation}
\label{eqn:ln-def}
	L_n:= \sum_{i\in [n]} w_i^{\sss(n)}.
\end{equation}
Now construct a random graph on $[n]$ by placing an edge between $i$ and $j$ independently for each $i< j\in [n]$ with probability $\thickbar{\mvq}_{ij}$.
This random graph model corresponds to the rank-$1$ case of the general class of IRGs studied by Bollob{\'a}s, Janson, and Riordan \cite{BJR07}.
We will denote this random graph by $\tbarGn$. 
A closely related model \cite{NorRei06,BJR07} is obtained by using the edge connection probabilities
\begin{equation}
\label{eqn:nr-connection}
{\mvq}_{ij}={\mvq}_{ij}(\mvw^{\sss(n)}):= 1-\exp(-w_i^{\sss(n)} w_j^{\sss(n)}/L_n)\, ,
\end{equation}
and will be denoted by $\mvG_n$.
In the regime of interest for this paper, as shown in \cite{janson-equiv}, this model is equivalent to the Chung-Lu model \cite{CL-connected,CL-distance,CL-distances-2,CL-book} and the Britton-Deijfen-Martin-L{\"o}f model \cite{britton-deijfen-martinlof}.  
Probabilists will be more familiar with the above random graph via its connection to one of the most famous stochastic coalescent models--the multiplicative coalescent. We describe only a special case and refer the interested reader to  \cite{aldous-crit, aldous-limic, aldous1999deterministic, bertoin-coagulation} for more general constructions of this specific model as well as other coalescent models.

\begin{figure}
	\centering
	\includegraphics[trim=4cm 3.2cm 3.5cm 3cm, clip=true, angle=0, scale=.3]{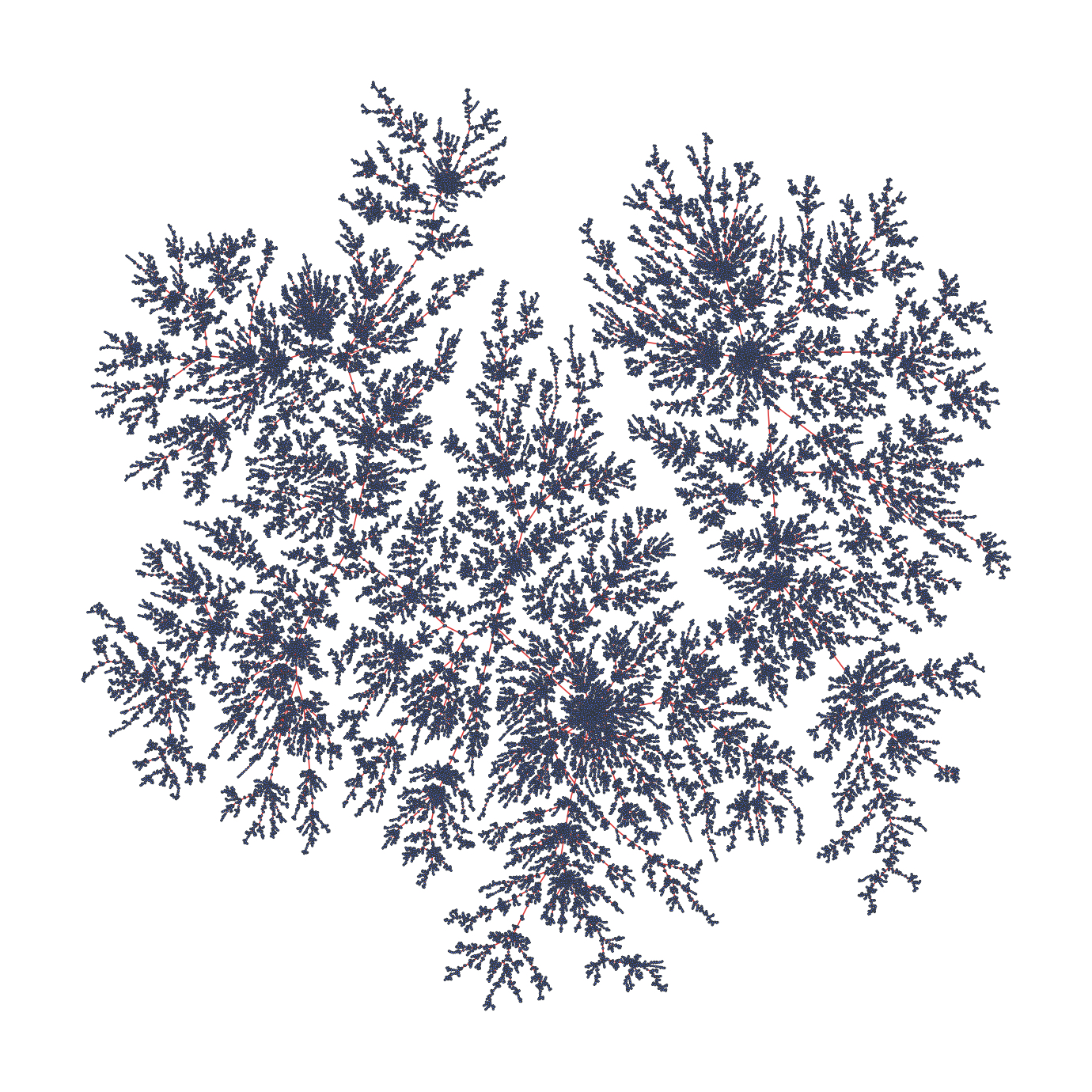}
	\captionof{figure}{The MST on the component of the vertex $1$ in $\overline\mvG_n$ with $w_i=3(n/i)^{\alpha}$, where $\tau=3.05$ and $n=80000$.}
	\label{fig:mst-3.05}
\end{figure}

\begin{figure}
	\centering
	\includegraphics[trim=3cm 3.2cm 3cm 3cm, clip=true, angle=0, scale=.3]{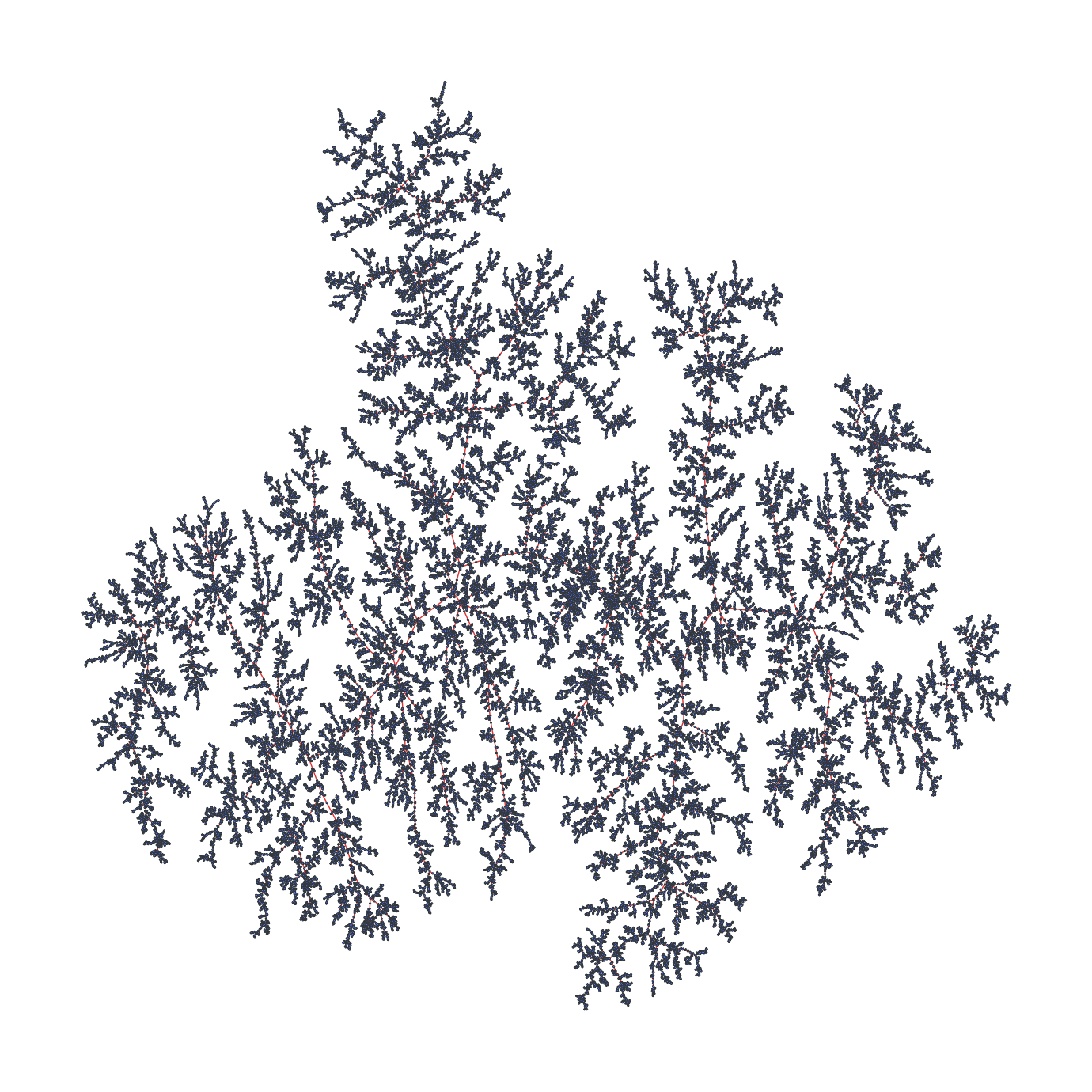}
	\captionof{figure}{The MST on the component of the vertex $1$ in $\overline\mvG_n$ with $w_i=3(n/i)^{\alpha}$, where $\tau=3.95$ and $n=80000$.}
	\label{fig:mst-3.95}
\end{figure}

\begin{defn}[Finite state multiplicative coalescent]
	\label{def:cgvx-def}
	Fix a finite vertex set $\cV$ and a collection of nonnegative vertex weights $\mvx = \big(x_v\, ;\, v\in \cV\big)$. Start at time zero with each vertex $v\in \cV$ being in a separate cluster of size $x_v$. At any time $t\geq 0$ we will have a collection of connected clusters with the weight of a cluster $\cC$ given by $\sW_{\mvx}(\cC) = \sum_{v\in \cC} x_v$.  Consider the continuous time Markov chain where two existing distinct clusters $\cC_a$ and $\cC_b$ merge at rate $\sW_{\mvx}(\cC_a)\cdot \sW_{\mvx}(\cC_b) $ into a single cluster of size $\sW_{\mvx}(\cC_a)+ \sW_{\mvx}(\cC_b)$. Write $\big(\MC\big((\cV, \mvx), t\big);\, t\geq 0\big)$ for this process. 
\end{defn}
Associated with the above dynamics is the following random graph related to $\mvG_n$. 
\begin{defn}[Random graph $\cG\big((\cV, \mvx), t\big)$]\label{def:1}
Consider a finite set $\cV$, nonnegative weights $\big(x_v\, ;\, v\in \cV\big)$, and $t\geq 0$. 
Let $\cG\big((\cV, \mvx), t\big)$ be  the random graph on vertex set $\cV$ obtained by independently placing edges between pairs of vertices $u,v\in \cV$ with probability
$ 1-\exp(- tx_v x_u)$.
\end{defn}

The following lemma is easy to check from the above description of the dynamics. 
\begin{lemma}\label{lem:mc-rg-graph}
	For any $t\geq 0$, the ordered sequence of weights of the connected components of $\cG((\cV, \mvx), t)$ has the same distribution as the ordered sequence of weights of clusters of $\MC((\cV, \mvx), t)$. 
\end{lemma}

The random graph $\cG$ is essentially the same as $\mvG_n$ with the parameters expressed in a different way.
This connection between the random graph $\mvG_n$ and the multiplicative coalescent will play a major role in our proofs.

We now specify how the vertex weights $\mvw^{\sss(n)}$ are chosen.
For the rest of this paper, we will work with a fixed exponent $\tau\in (3, 4)$.
We will use the following notation for constants associated to this exponent:
\begin{equation}\label{eqn:const-def}
\alpha := \frac{1}{\tau-1}\, , \qquad \rho := \frac{\tau-2}{\tau-1}\, ,\ \ \text{ and }\ \  \eta:= \frac{\tau-3}{\tau-1}\, . 
\end{equation}

\begin{ass}\label{ass:wts}
Consider a sequence $\big(\mvw^{\sss(n)}\, ;\, n\geq 1\big)$ of weight sequences, where $\mvw^{\sss(n)}= \big(w_1^{\sss(n)}, w_2^{\sss(n)}, \ldots, w_n^{\sss(n)}\big)$ with $w_1^{\sss(n)}\geq w_2^{\sss(n)}\geq\ldots\geq w_n^{\sss(n)}> 0$ . 
Then we assume the following:
\begin{enumeratei}
\item 
{\bf Supercriticality condition:}
Let $L_n$ be as in \eqref{eqn:ln-def}. 
Then
\[
\nu:=
\liminf_{n\to\infty}\ \frac{\sum_{j=1}^n (w_j^{\sss(n)})^2}{L_n}
>1\, .
\]
\item 
For each  $ i\geq 1$, there exists $\theta_i^*> 0$ such that 
\[
\lim_{n\to\infty} 
\frac{w_i^{\sss(n)}}{n^{\alpha}}
\cdot
\bigg(
\frac{n}{\sum_{j=1}^n (w_j^{\sss(n)})^2}
\bigg)^{1/2} 
= 
\theta_i^*.
\]
\item 
There exist constants $A_1, A_2 \in (0,\infty)$ such that for all $n\geq 2$ and $1\leq i\leq n/2$,
\[
A_1\left(\frac{n}{i}\right)^{\alpha} 
\leq 
w_i^{\sss(n)}
\leq 
A_2\left(\frac{n}{i}\right)^{\alpha}.
\]
\item
For all $n\geq 2$, 
$w_n^{\sss (n)}\geq A_1 \big(\log n\big)^{\frac{3}{2}}\cdot n^{-\frac{\eta}{4}}$\, .		
\end{enumeratei}
\end{ass}

Assumption~\ref{ass:wts}~(i) ensures supercriticality of the random graph model, and Assumption~\ref{ass:wts}~(ii) corresponds to the condition in \cite[Display (19)]{aldous-limic}.
It should be possible to relax Assumption~\ref{ass:wts}~(iii) and (iv) with a more intricate analysis; however, we do not pursue it here.
We will discuss this briefly in Section~\ref{sec:disc}.

Note that Assumption \ref{ass:wts} implies that for all $n\geq 2$ and $n/2\leq i\leq n$,
\begin{align}\label{eqn:up-bnd}
w_i^{\sss(n)}
\leq 
w_{n/2}^{\sss(n)}
\leq
2^{\alpha}A_2
\leq
2^{\alpha}A_2\left(\frac{n}{i}\right)^{\alpha}.
\end{align}
Write
\begin{align}\label{eqn:parameters}
\pmtr:=\big(A_1, A_2, \tau, \nu\big)
\end{align}
for the parameters in Assumption \ref{ass:wts}.
The following lemma is easy to verify and gives two natural settings that give rise to weights satisfying the above condition. The proof is omitted.

\begin{lemma}\label{lem:weight-choice}
Suppose $F$ is a cumulative distribution function (cdf) with support in $[0, \infty)$ such that for some $\tau \in (3,4)$, $\beta_{F}>4/\eta$, and $c_F\in (0, \infty)$, 
\begin{align*}
\limsup_{x\downarrow 0} \big(x^{-\beta_F}F(x)\big)<\infty , \qquad \, 
\lim_{x\to\infty} x^{\tau-1} [1-F(x)]= c_{\sss F}\, ,\ \text{ and }\ \
\int_0^{\infty} x^2 F(dx)> \int_0^{\infty} x F(dx)\, . 
\end{align*}
	Consider $\mvw^{\sss(n)}=\big(w_i^{\sss(n)} ,\,  1\leq i\leq n\big)$ obtained via one of the following:
	\begin{enumeratea}
		\item Let $w_i^{\sss(n)} = [1-F]^{-1}\big(i/(n+1)\big) $ for $1\leq i\leq n$. 
		\item Let $W_1, \ldots, W_n$ be i.i.d. random variables with cdf $F$, and let $w_1^{\sss(n)}\geq\ldots\geq w_n^{\sss(n)}$ be the corresponding ordered values.
	\end{enumeratea}
Then under {\upshape (a)}, $\big(\mvw^{\sss(n)},\, n\geq 1\big)$ satisfies Assumption \ref{ass:wts}.
Under {\upshape (b)}, we can construct $\mvw^{\sss(n)}$, $n\geq 1$, on the same probability space such that Assumption \ref{ass:wts} is satisfied almost surely in this space with a deterministic $\nu$, and random $A_1$, $A_2$, and $\theta_i^*$, $i\geq 1$.
\end{lemma}

Write $\cC(1, \mvG_n)$ (resp. $\cC(1, \tbarGn)$) for the component of vertex $1$ in the graph $\mvG_n$ (resp. $\tbarGn$).
For any component $\cC$ of $\mvG_n$ or $\tbarGn$, write $\cW(\cC)=\sum_{i\in\cC}w_i^{\sss(n)}$.
The following result states some basic properties of these two random graphs. 
\ch{The notations $O_P(\cdot)$ and $\Theta_P(\cdot)$ used in the sequel are as explained at the end of Section~\ref{sec:notation}.}

\begin{prop}
	\label{prop:prop-of-rg}
	The following hold under Assumption~\ref{ass:wts}:
	\begin{enumeratea}
		\item We have, $\sW\big(\cC(1, \mvG_n)\big) = \Theta_P(n)$, and
		\begin{align}\label{eqn:999}
		\max\big\{\sW(\cC)\, :\, \cC\text{ component in }\mvG_n\text{ and }\cC\neq \cC(1, \mvG_n)\big\} 
		 = O_P(\log{n})\, .
		\end{align}
		We call $\cC(1, \mvG_n)$ the giant component in $\mvG_n$. 
		Further, an analogous result holds for the random graph $\tbarGn$.
		\item By \cite[Theorem 3.13]{BJR07}, if further the empirical distribution $n^{-1} \sum_{i\in [n]} \delta \big\{w_i^{\sss(n)}\big\}$ converges to a cdf with tail exponent $\tau \in (3,4)$, then the degree distribution of $\mvG_n$ (resp. $\tbarGn$) converges in probability to a deterministic distribution with tail exponent $\tau$.  
	\end{enumeratea}
\end{prop}
The claim that $\sW\big(\cC(1, \mvG_n)\big)=\Theta_P(n)$ follows from Proposition \ref{prop:lem-ss-3}, \eqref{eqn:ss-8}, and Lemma \ref{lem:minus-coup}.
The proof of \eqref{eqn:999} is similar to that of Lemma \ref{lem:ss-lem-12}.
Proposition \ref{prop:prop-of-rg} (b) is not important for this study, but it gives some justification as to why these models are used to understand real world systems. 
We now define the central object of interest in this paper. 

\begin{defn}[Minimal spanning tree]\label{def:mst}
Let $\mvU=\big(U_{ij}\, ;\, 1\leq i < j\leq n\big)$ be a collection of i.i.d. $\Unif\, [0,1]$ random variables, and let $U_{ji}=U_{ij}$ for $1\leq i < j\leq n$. 
Let $\mvM^n$ denote the minimal spanning tree on the giant component of $\mvG_n$ using the edge weights $\big(U_{ij}\,;\, \{i, j\} \text{ is an edge in } \cC(1, \mvG_n)\big)$. 
Define $\tbar{\mvM}^n$ in analogous way.
\end{defn}

We make a convention here that we will follow throughout this paper:
\begin{quote}
When a finite connected graph $H$ is viewed as a \ch{metric} space, the underlying set will be the collection of vertices in $H$ joined by line segments of unit length that represent the edges in $H$, and the distance between two points will be the minimum of the lengths of paths connecting the two points.
For any metric space $(X, d)$ and $a>0$, $a\cdot X$ will denote the metric space $(X, a\cdot d)$, i.e, the space where the distance is scaled by $a$.
\end{quote}
Using this convention, we can view $\mvM^n$ and $\tbar{\mvM}^n$ as (random) compact metric spaces.

\section{Main result}
\label{sec:res}
Recall that the lower and upper box counting dimensions of a compact metric space $X$ are given by
\[\udim(X):= \liminf_{\delta\downarrow 0}\, \frac{\log{[\cN(X,\delta)]}}{\log(1/\delta)},
\quad\text{and}\quad
\odim(X):= \limsup_{\delta\downarrow 0}\, \frac{\log{[\cN(X,\delta)]}}{\log(1/\delta)}\]
respectively, where $\cN(X,\delta)$ is the minimum number of closed balls with radius $\delta$ required to cover $X$. 
If $\udim(X)=\odim(X)$, then the box-counting dimension or the Minkowski dimension of $X$ exists and equals this common value. 
For an $\bR$-graph $X$ and any point $x\in X$,  let $\deg(x\, ;\, X)$ denote the degree of $x$ in $X$; see Section \ref{sec:r-tree} for the relevant definitions. 
We write
\begin{equation}\label{eqn:leaf-hub-def}
	\cL(X):=\set{x\in X: \deg(x\, ;\, X) =1 }\, ,\ \ \text{ and }\ \ \  
	\cH(X) := \set{x\in X: \deg(x\, ;\, X) =\infty }
\end{equation}
for the set of leaves and the set of hubs in $X$ respectively.
Recall the notation from \eqref{eqn:const-def}.

\begin{thm}\label{thm:mst-convg}
Under Assumption \ref{ass:wts} on the weight sequence, there exists a random compact $\bR$-tree $\mst^{\mvtheta^*}$ whose law depends only on $\mvtheta^*:= (\theta_i^*\, ;\, i\geq 1)$ such that
\begin{align}\label{eqn:666}
n^{-\eta}\cdot \mvM^n \convd \mst^{\mvtheta^*} , \ \ \text{ as }\ \ \ n\to\infty\, ,
\end{align}
with respect to the Gromov-Hausdorff topology.
Further, almost surely,
\begin{enumeratea}
\item\label{it:a} 
$\deg\big(x\,;\, \mst^{\mvtheta^*}\big) \in \big\{1,2,\infty\big\}$ for all $x\in \mst^{\mvtheta^*}$;
\item\label{it:b} 
both the set of leaves $\cL(\mst^{\mvtheta^*})$ and the set of hubs $\cH(\mst^{\mvtheta^*})$ are dense in $\mst^{\mvtheta^*}$; and
\item\label{it:c} 
the Minkowski dimension of $\mst^{\mvtheta^*}$ satisfies
\[
\dim\big(\mst^{\mvtheta^*}\big) = \frac{1}{\eta}=\frac{\tau-1}{\tau -3}\, .
\]
\end{enumeratea}
Moreover, \eqref{eqn:666} continues to hold under Assumption \ref{ass:wts} if we replace $\mvM^n$ by $\tbar{\mvM}^n$.
\end{thm}

It should be possible to lift the Gromov-Hausdorff convergence in \eqref{eqn:666} to convergence with respect to the Gromov-Hausdorff-Prokhorov topology.
For this, it would suffice to prove an additional technical condition.
We will discuss this further in Section \ref{sec:disc}.

\section{Definitions and preliminary results}
\label{sec:defn}

\subsection{Notation}\label{sec:notation}
Throughout this paper, $C, C'$ etc. will denote constants that depend only on $\pmtr$ as defined in \eqref{eqn:parameters}, and their values may change from line to line. 
Special constants will be indexed by the relevant equations, e.g., $C_{\ref{eqn:55}}$ etc., and their values will depend only on $\pmtr$ unless specified otherwise. 
If a constant depends on any parameter other than $\pmtr$, then that will be explicitly mentioned when the constant is first introduced.
For example, $C_{\ref{eqn:555}}=C_{\ref{eqn:555}}(\Delta)$ appearing in the statement of Theorem \ref{thm:mst-crit-gh} below depends on $\pmtr$ and $\Delta$.

A claim holds `for all large $n$' or `$\falln$' will imply that there exists $n_0\geq 1$ depending only on $\pmtr$ such that the claim holds for all $n\geq n_0$. 
If the threshold depends on any other parameter, then that will be explicitly mentioned, and it will be assumed, without explicit mention, that such a threshold is chosen bigger than all thresholds involving $n$ that depend only on $\pmtr$ and were previously introduced in the proof.
For example, $n_{\ref{eqn:9A}} = n_{\ref{eqn:9A}}(\kappa)$ appearing around \eqref{eqn:9A} depends on $\pmtr$ and $\kappa$, and although not explicitly mentioned, is chosen so that it is bigger than the threshold involving $n$ above which \eqref{eqn:ss-7a} and \eqref{eqn:ss-8} hold.

Similarly, all thresholds involving $\lambda$ or $\eps$ will depend only on $\pmtr$ unless specified otherwise when the threshold is first introduced.

A relation of the form $a\asymp b$ will mean that there exist $C, C'>0$ depending only on $\pmtr$ such that $Ca\leq b\leq C' b$.
Here, $a, b$ could be elements of two sequences, or two functions defined on the same domain. 
As an example, consider the following claim made around \eqref{eqn:ss-6} below:
``\ldots for all large $n$,
$i_\lambda(u) \asymp u^{1/\alpha}\ \ \text{ for }\ \ (\lambda, u)\in\dI^{\sss(n),2}$."
This statement can be rewritten as follows: 
There exist $C, C'>0$, and $n_0\geq 1$ depending only on $\pmtr$ such that 
$Cu^{1/\alpha}\leq i_\lambda(u)\leq C' u^{1/\alpha}$
for all $n\geq n_0$ and $(\lambda, u)\in\dI^{\sss(n),2}$.

For a set $S$, we use $|S|$ or $\#S$ to denote the number of elements in $S$. 
For any graph $H$, we write $V(H)$ and $E(H)$ for the set of vertices and the set of edges of $H$ respectively.
We write $|H|$ for the number of vertices in $H$, i.e., $|H|=|V(H)|$.
For any finite connected graph $H=(V, E)$, we write $\spls(H)$ for the number of surplus edges in $H$, i.e.,
$\spls(H):=|E|-|V|+1$.
For a finite (not necessarily connected) graph $H$, we let
\begin{align}\label{eqn:234}
\max\surplus(H)
=
\max\big\{\surplus(\cC)\, :\, \cC\text{ connected component of }H \big\}\, .
\end{align}
Similarly, we write $\diam(H)$ for the maximum of the diameters (with respect to graph distance) of all the components in $H$.
We will write $\mathrm{LP}(H)$ to denote the length of the longest self-avoiding path in $H$.
For a rooted tree $T$, $\height(T)$ will denote the height of $T$, i.e., the distance to the farthest leaf from the root of $T$.

For a graph $H=(V, E)$ and $v\in V$, we write $\cC(v, H)$ to denote the component of $v$ in $H$. 
For $V'\subseteq V$, we write $H\setminus V'$ to denote the restriction of $H$ to the vertex set $V\setminus V'$.
If the vertices in $H$ have weights $\mvy = (y_v;\, v\in V)$ associated to them, then for any subgraph $H_0$ in $H$, we write $\sW_{\mvy}(H_0) = \sum_{v\in V(H_0)} y_v$ for the weight or mass of $H_0$ as measured by the prescribed vertex weights. 
When the weight sequence is $\mvw$, we omit the subscript and simply write $\sW(\cdot)$.

For two real valued random variables $X, Y$ we write $X\stod Y$ for the stochastic domination relation between the distributions of these random variables.  We use $\Poi(\cdot)$, $\Unif(\cdot)$ and $\Bern(\cdot)$ to respectively denote Poisson, uniform, and Bernoulli distributions with parameters in $\cdot$ that will be specified in the setting of interest.

For a non-negative function $n\mapsto g(n)$,
we write $f(n)=O(g(n))$ when $|f(n)|/g(n)$ is uniformly bounded, and
$f(n)=o(g(n))$ when $\lim_{n\rightarrow \infty} f(n)/g(n)=0$.
Furthermore, we write $f(n)=\Theta(g(n))$ if $f(n)=O(g(n))$ and $g(n)=O(f(n))$. 
We let $\weakc$, $\probc$, and $\convas$ respectively denote convergence in distribution, convergence in probability, and almost sure convergence.
\ch{For a sequence of random variables
$(X_n;\, n\geq 1)$ and a sequence of positive real numbers $(b_n;\ n\geq 1)$, we write 
$X_n=O_P(b_n)$ if the sequence of laws of $(X_n/ b_n;\, n\geq 1)$ is tight, and 
$X_n=\op(b_n)$ if $X_n/b_n\weakc 0$ as $n\rightarrow\infty$.
Further, we write $X_n=\Theta_P(b_n)$ if $X_n = O_P(b_n)$ and $1/X_n = O_P(1/b_n)$.}
For simplicity, we will freely omit ceilings and floors; this will not affect the argument.

\subsection{Convergence of metric spaces}\label{sec:metric-convg}
The following five topologies will be relevant to us:
(i) the Hausdorff topology on closed subsets of a compact metric space,
(ii) the Gromov-Hausdorff (GH) topology on $\fS_{\GH}$--the isometry equivalence classes of compact metric spaces,
(iii) the marked Gromov-Hausdorff topology on the isometry equivalence classes of triples of the form $\big(X, C, d\big)$, where $(X, d)$ is a compact metric space and $C\subseteq X$ is closed, 
(iv) the Gromov-Hausdorff-Prokhorov (GHP) topology on the isometry equivalence classes of compact metric measure spaces, and
(v) the Gromov-weak topology on the isometry equivalence classes of metric measure spaces.
The Hausdorff distance and the GHP distance will be denoted by $d_{\rH}(\,\cdot\, ,\,\cdot\,)$ and $d_{\GHP}(\,\cdot\, ,\,\cdot\,)$ respectively.
We refer the reader to \cite{metric-geometry-book} for background on the topologies in (i) and (ii).
For the topologies in (iii) and (iv), we will primarily follow \cite{miermont2009tessellations} and \cite{abraham-delmas-hoscheit, AddBroGolMie13} respectively.
For the Gromov-weak topology, we refer the reader to \cite{AthLorWin14, Winter-gromov-weak}.
The definitions and results related to the topologies in (i)--(iv) needed in this paper can be found in one place in \cite[Section 3.2]{addarioberry-sen}.

\subsection{$\bR$-trees and $\bR$-graphs} \label{sec:r-tree}
For any metric space $(X, d)$, a geodesic between $x_1, x_2\in X$ is an isomeric embedding $f:[0, d(x_1,x_2)]\to X$ such that $f(0)=x_1$ and $f\big(d(x_1,x_2)\big)=x_2$.
$(X, d)$ is a geodesic space if there is a geodesic between any two points in $X$.
An embedded cycle in $X$ is a subset of $X$ that is a homeomorphic image of the unit circle $S^1$.
\begin{defn}[Real trees \cite{legall-survey,evans-book}]\label{def:real-tree}
A compact geodesic metric space $(X,d)$ is called a real tree or $\bR$-tree if it has no embedded cycles.
\end{defn}
For a metric space $(X, d)$, $x\in X$, and $\eps>0$, let 
$B(x,\eps\, ; X):=\big\{y\in X\, :\, d(y, x)\leq \eps\}$.
\ch{We next recall some definitions and constructs from \cite{AddBroGolMie13}.
We refer the reader to \cite[Section~2.3]{AddBroGolMie13} for a detailed treatment.}

\begin{defn}[$\bR$-graphs \cite{AddBroGolMie13}]\label{def:r-graph}
A compact geodesic metric space $(X,d)$ is called an $\bR$-graph if for every $x\in X$, there exists $\eps>0$ such that $\big(B(x,\eps\, ; X), d\rvert_{B(x,\eps\, ;\, X)}\big)$ is an $\bR$-tree.
	
The core of an $\bR$-graph $(X,d)$, denoted by $\core(X)$, is the union of all the simple arcs having both endpoints in embedded cycles of $X$. If it is non-empty, then $(\core(X), d)$ is an $\bR$-graph with no leaves.
We define $\conn(X)$ to be the set of all $x\in X$ such that $x$ belongs to an embedded cycle in $X$.
\end{defn}

Clearly, $\conn(X)\subseteq\core(X)$.
By \cite[Theorem 2.7]{AddBroGolMie13}, if $X$ is an $\bR$-graph with a non-empty core, then $(\core(X), d)$ can be represented as $(k(X), e(X), l)$, where $(k(X), e(X))$ is a finite connected multigraph in which all vertices have degree at least $3$ and $l: e(X)\to (0,\infty)$ gives the edge lengths of this multigraph.
We denote by $\surplus(X)$ the number of surplus edges in $(k(X), e(X))$.
On any $\bR$-graph $(X,d)$ there exists a unique $\sigma$-finite Borel measure $\len$ called the length measure such that if $x_1, x_2\in X$ and $[x_1,x_2]$ is a geodesic path between $x_1$ and $x_2$ then $\len\big([x_1,x_2]\big)=d(x_1, x_2)$.
Note that
\begin{align}\label{eqn:def-L}
\sum_{e\in e(X)} l(e)=\len\big(\core(X)\big).
\end{align}
Clearly, $\len\big(\conn(X)\big)\leq\len\big(\core(X)\big)<\infty$.
If $\conn(X)\neq\emptyset$ (in which case $\len\big(\conn(X)\big)>0$), we write $\len_{\conn(X)}$ for the restriction of the length measure to $\conn(X)$ normalized to be a probability measure, i.e.,
\begin{align*}
\len_{\conn(X)}(\cdot)=\frac{\len(\cdot)}{\len(\conn(X))}~.
\end{align*}

For an $\bR$-graph $(X, d)$ and $x\in X$, choose $\eps>0$ such that $B(x,\eps\, ; X)$ is an $\bR$-tree, and define the degree of $x$ as 
\[
\deg(x\, ; X):=\big|\big\{\text{connected components of }B(x,\eps\, ; X)\setminus\{x\}\big\}\big|\, .
\]
Note that the value of $\deg(x\, ; X)$ is independent of the choice of $\eps$.

Since any finite connected graph, viewed as a metric space, is an $\bR$-graph,
the above definitions make sense for any finite connected graph $H$.
Note the difference between $e(H)$ defined above and $E(H)$-the set of edges in $H$.
Note also that in this case, the graph theoretic $2$-core of $H$, viewed as a metric space, coincides with the space $\core(H)$ as defined above, and
$\len(\core(H))$ equals the number of edges in the graph theoretic $2$-core of $H$.

\subsection{Some properties of MSTs}\label{sec:mst-properties}
Suppose $H=(V, E, b)$ is a weighted, connected, and labeled graph.
Assume that $b(e)\neq b(e')$ whenever $e\neq e'$.
We now state a useful property of the MST.




\begin{lemma}[Minimax paths property]\label{lem:mst-minimax-criterion}
	Let $H=(V,E,b)$ be as above.
	Then the MST $T$ of $H$ is unique.
	Further, $T$ has the following property:
	Any path $(x_0,\hdots,x_n)$ with $x_i\in V$ and $\{x_i,x_{i+1}\}\in E(T)$
	satisfies
	\[\max_{1\leq i\leq n}\ b\big(\{x_i,x_{i+1}\}\big)\leq \max_{1\leq j\leq m}\ b\big(\{x_j',x_{j+1}'\}\big)\]
	for any path $(x_0',\hdots,x_m')$ with $\{x_j',x_{j+1}'\}\in E$
	and $x_0=x_0'$ and $x_n=x_m'$.
	In words, the maximum edge weight in the path in the MST connecting two given vertices is smallest among all paths in $G$ connecting those two vertices.
	
	Moreover, $T$ is the only spanning tree of $G$ with the above property.
\end{lemma}
The above lemma is just a restatement of \cite[Lemma 2]{kesten-lee};
see also \cite[Proposition 2.1]{alexander-forest}.
We record the following useful observations that follow directly from Lemma~\ref{lem:mst-minimax-criterion}:
\begin{obs}\label{observation:only-rank-matters}
	The MST can be constructed just from the ranks of the different edge weights.
	This fact is not needed in the sequel, but it shows that the laws of $\mvM^n$ and $\tbar\mvM^n$ will remain unchanged if we used a set of exchangeable and pairwise distinct edge weights instead of $\mvU$ as in Definition~\ref{def:mst}.
\end{obs}
\begin{obs}\label{observation:percolation}
	Let $H=(V, E, b)$ be a connected and labeled graph with pairwise distinct edge weights.
	Let $u\in[0,\infty)$ and $\cC$ be a component of the graph $G^{u}=(V, E^u)$,
	where $E^u\subseteq E$ contains only those edges $e$ for which $b(e)\leq u$.
	Then the restriction of the MST of $(V, E, b)$ to $\cC$ is the MST of $\big(V(\cC), E(\cC), b|_{E(\cC)}\big)$.
	This fact is extremely useful as it can be used to connect the structure of the MST to the geometry of components of the graph under percolation.
\end{obs}

\subsection{Cycle-breaking}\label{sec:cycle-breaking}
\ch{In this section we recall two procedures from \cite{AddBroGolMie13} that can be applied to $\bR$-graphs and combinatorial graphs.
We refer the reader to \cite[Sections~3.1 and 3.2]{AddBroGolMie13} for a detailed treatment.}
Recall the notation $k(X)$, $e(X)$, $(l(e),\ e\in e(X))$, and $\mathrm{sp}(X)$ introduced below Definition~\ref{def:r-graph}.

\begin{defn}[Cycle-breaking ($\cb$), \ch{\cite[Section~3.2]{AddBroGolMie13}}]\label{def:algo-cb}
Let $X$ be an $\bR$-graph. If $X$ has no embedded cycles, then set $\cb(X)=X$.
Otherwise, sample $x\in X$ using the measure $\len_{\conn(X)}$. 
Endow $X\setminus \{x\}$ with the intrinsic metric: the distance between two points is the minimum of the lengths of paths in $X\setminus \{x\}$ that connect the two points.
Set $\cb(X)$ to be the completion of $X\setminus \{x\}$ with respect to the intrinsic metric.
(Thus, $\cb(X)$ is also an $\bR$-graph.)
	
For $k\geq 2$, we inductively define $\cb^k(X)$ to be the space $\cb\big(\cb^{k-1}(X)\big)$.
(Thus, at the $k$-th step, if $\cb^{k-1}(X)$ has an embedded cycle, then we are using the measure to $\len_{\conn(\cb^{k-1}(X))}$ to sample a point.)
	
Note that $\cb^k(X)=\cb^{\mathrm{sp}(X)}(X)$ for all $k\geq\mathrm{sp}(X)$, i.e., the spaces $\cb^k(X)$ remain the same after all cycles have been cut open.
We denote this final space (which is a real tree) by $\cb^\infty(X)$.
\end{defn}

Next we define a cycle-breaking process for discrete graphs.

\begin{defn}[Cycle-breaking for discrete graphs ($\cbd$), \ch{\cite[Section~3.1]{AddBroGolMie13}}]\label{def:algo-cbd}
Let $H=(V, E)$ be a finite connected graph.
Sample $e\in E$ uniformly. 
If $(V, E\setminus \{e\})$ is connected, set 
$\cbd(H)=(V, E\setminus \{e\})$.
Otherwise, set $\cbd(H)=H$.
Inductively set $\cbd^{k+1}(H)=\cbd\big(\cbd^k(H)\big)$, $k\geq 1$.

Almost surely, the graphs $\cbd^k(H)$ are the same (and are all trees) for all large values of $k$.
We denote this tree by $\cbd^\infty(H)$.
\end{defn}

Suppose $H$ is a finite connected graph.
Let $f_1,\ldots,f_s$ be the edges of $H$ that get removed in the process $\big(\cbd^k(H), k\geq 1\big)$. Clearly, $s=\mathrm{sp}(H)$.
For $1\leq i\leq s$, let $y_i$ be a uniformly sampled point on $f_i$.
It is easy to see that viewing $H$ as an $\bR$-graph, the completion of the space $H\setminus\{y_1,\ldots,y_s\}$ with repect to the intrinsic metric has the same distribution as $\cb^\infty(H)$.
In this coupling, $\cbd^\infty(H)$ is a subspace of $\cb^\infty(H)$, and
\begin{align}\label{eqn:cbd-cb-close}
d_{\rH}\big(\cbd^\infty(H),\ \cb^\infty(H)\big)\leq 1\, .
\end{align}
We now state a lemma that connects cycle-breaking to MSTs.

\begin{lem}\label{lem:cycle-breaking-gives-mst}
Suppose $H$ is a finite connected graph.
Then $\cbd^\infty(H)$ has the same law as the MST of $H$ constructed by assigning exchangeable pairwise distinct weights to the edges in $H$.
\end{lem}

Lemma \ref{lem:cycle-breaking-gives-mst} follows easily from Lemma \ref{lem:mst-minimax-criterion}.
A proof can be found in \cite[Proposition 3.5]{AddBroGolMie13}.
For $r\in (0,1)$ define $\cA_r$ to be the set of all $\bR$-graphs $X$ that satisfy
\begin{align}\label{eqn:65}
\mathrm{sp}(X)+\len(\core(X))\leq 1/r,\ \text{ and }\ \min_{e\in e(X)}\len(e)\geq r\, .
\end{align}
The following theorem will allow us to prove convergence of MSTs from GHP convergence of the underlying graphs.

\begin{thm}\label{thm:mst-from-ghp-convergence}
Fix $r\in (0,1)$.
Suppose $(X,d)$ and $(X_n, d_n)$, $n\geq 1$, are $\bR$-graphs in $\cA_r$ such that 
$(X_n, d_n)\to (X,d)$ as $n\to\infty$ w.r.t. GH topology.
Further, suppose for each $n\geq 1$, $(X_n, d_n)$ is isometric to $\eps_n\cdot H_n$, where $H_n$, $n\geq 1$, are finite connected graphs and $\eps_n\to 0$.
Then as $n\to\infty$, $\eps_n\cdot\cbd^{\infty}(H_n)\weakc\cb^{\infty}(X)$ w.r.t. GH topology.
\end{thm}
Theorem \ref{thm:mst-from-ghp-convergence} follows from \cite[Theorem 3.3]{AddBroGolMie13} and \eqref{eqn:cbd-cb-close}.

\section{A slightly different model}\label{sec:different-model}
Recall the definition of $L_n$ from \eqref{eqn:ln-def}, and let
\begin{align}\label{eqn:444}
\ell_n := \sum_{i=2}^n w_i^{\sss(n)}=L_n-w_1^{\sss(n)}\, .
\end{align}
Analogous to \eqref{eqn:rank-one-connection} and \eqref{eqn:nr-connection}, define
\begin{align}\label{eqn:300}
\bar q_{ij}=\bar q_{ij}(\mvw^{\sss (n)})=1\wedge \big(w_i^{\sss (n)}w_j^{\sss (n)}/\ell_n\big)\, ,
\ \text{ and }\ \ 
q_{ij}=q_{ij}(\mvw^{\sss (n)})=1-\exp\big(-w_i^{\sss (n)}w_j^{\sss (n)}/\ell_n\big)
\end{align}
for $1\leq i\neq j\leq n$.
Note that under Assumption \ref{ass:wts}, for all large $n$,
$\bar q_{ij}=w_i^{\sss (n)}w_j^{\sss (n)}/\ell_n$ for all $1\leq i\neq j\leq n$.
Let $\barGn$ (resp. $G_n$) be the random graph on $[n]$ obtained by placing an edge between $i$ and $j$ independently for each $i< j\in [n]$ with probability $\bar q_{ij}$ (resp. $q_{ij}$).

\begin{lemma}\label{lem:minus-coup}
Under Assumption \ref{ass:wts}, 
there exists a coupling of $\tbarGn$ and $\barGn$ such that $\pr\big(\tbarGn\neq\barGn\big) \to 0$ as $n\to\infty$. 
Similar assertions hold for the pairs $\big(\mvG_n\, , G_n\big)$, and $\big(\barGn\, , G_n\big)$.
\end{lemma}

\noindent{\bf Proof:}
We will prove the assertion for $\tbarGn$ and $\barGn$; a similar argument works for the other pairs.
Under Assumption \ref{ass:wts}, $\max_{i, j}\big(w_i^{\sss (n)}w_j^{\sss (n)}/\ell_n\big)=O(n^{-\eta})$.
By \cite[Corollary 2.12]{janson-equiv}, it is enough to prove that
$\sum_{i\neq j} {(\thickbar{\mvq}_{ij} - \overline q_{ij})^2 }/\overline q_{ij} \to 0$, as $n\to\infty$.
Now, under Assumption \ref{ass:wts}, for all large $n$,
\[
\sum_{i\neq j} \frac{(\thickbar{\mvq}_{ij} - \overline q_{ij})^2}{\overline q_{ij}} 
= 
\sum_{i\neq j} \frac{w_i w_j}{\ell_n}\bigg(1-\frac{\ell_n}{L_n}\bigg)^2 
\leq
\frac{L_n^2}{\ell_n}\cdot \bigg(\frac{w_1}{L_n}\bigg)^2
\leq Cn^{-\eta}\, .
\]
This completes the proof.
\qed

\medskip

Similar to Definition \ref{def:mst}, let $M^n$ denote the minimal spanning tree on $\cC(1, G_n)$ using the edge weights $\big(U_{ij}\,;\, \{i, j\}\in E\big(\cC(1, G_n)\big)\big)$.
Define $\tbar M^n$ as the MST on $\cC(1, \barGn)$ in an analogous way.
In view of Lemma \ref{lem:minus-coup}, it is enough to prove Theorem \ref{thm:mst-convg} for $M^n$ and $\tbar M^n$. 
Hence, from now on we will only work with the random graphs $G_n$ and $\barGn$.
This will make some of the computations simpler.

\section{Relating the MST with the components in the critical window}\label{sec:7}

We will now need some notation.  
Recall that $\ell_n=\sum_{i=2}^n w_i^{\sss (n)}$, and define
\begin{equation}
\label{eqn:sigm2-nu-def}
\sigma_2^{\sss(n)} = n^{-1}\sum_{i=2}^n (w_i^{\sss(n)})^2\, , \ \ \text{ and }\ \ 
\nu_n:= n\sigma_2^{\sss(n)}/\ell_n\, .
\end{equation}
(Note that once again we have omitted the weight of vertex $1$ in the formulae.) 
Using Assumption \ref{ass:wts}, it is easy to see that $\liminf_{n\to\infty}\nu_n=\nu>1$.
Without loss of generality, we can assume the following:
\begin{ass}\label{ass:wlog}
For all $n\geq 1$, $\nu_n\geq \nu':=1+(\nu-1)/2$.
\end{ass}

For $\lambda\geq 0$, define 
\begin{equation}
\label{eqn:pnlam-def}
p^n_{\lambda}:= \bigg(1+\frac{\lambda}{n^{\eta}}\bigg)\frac{1}{\nu_n}\, .
\end{equation}
Let $\mvU$ be as in Definition \ref{def:mst}.
For any subgraph $H$ of $\barGn$, we write 
\begin{equation}
\label{eqn:uni-def}
\mvU\big|_{H}:= \big(U_{ij}\, ;\, \{i, j\} \in E(H)\big)\, .
\end{equation}
In words, $\mvU\big|_{H}$ is the collection of edge weights corresponding to the edges in $H$. 

\begin{defn}[The graph $\barGn(\lambda)$]\label{def:graph-vgnp}
Fix $\lambda\geq 0$. 
Let $\barGn(\lambda)$ be the subgraph of $\barGn$ with edge set
\[
\big\{\{i, j\} \, :\,  \{i, j\}  \in E(\barGn) \text{ and } U_{ij} \leq p^n_{\lambda}\big\}\, .
\]
For $v\in[n]$, let $\tbar\cC_v^n(\lambda)=\cC\big(v, \barGn(\lambda)\big)$.
Let $\tbar M_{\lambda}^n$ denote the MST on $\tbar\cC_1^n(\lambda)$ constructed using $\mvU\big|_{\tbar\cC_1^n(\lambda)}$. 
\end{defn}

Thus,
$\tbar M^n=\tbar M^n_{\lambda}$ for all $\lambda\geq (\nu_n-1)n^{\eta}$.
Note that if $p^n_{\lambda}\in[0, 1]$ and $w_i^{\sss (n)}w_j^{\sss (n)}/\ell_n\leq 1$ for all $i\neq j\in[n]$, then $\barGn(\lambda)$ is a random graph on $[n]$ with independent edge connection probabilities $\big(p^n_{\lambda} w_i^{\sss(n)} w_j^{\sss(n)}/\ell_n\ ;\, 1\leq i< j\leq n\big)$.
Note also that by Observation \ref{observation:percolation},
$\tbar M_{\lambda_1}^n$ is a subtree of $\tbar M_{\lambda_2}^n$ whenever $0\leq\lambda_1\leq\lambda_2$.
In particular, $\tbar M_{\lambda}^n$ is a subtree of $\tbar M^n$ for all $\lambda\geq 0$.

\begin{defn}[The graph $G_n(\lambda)$]	\label{def:graph-vgexp}
Fix $\lambda\geq 0$. 
Let $G_n(\lambda)$ be the random graph on $[n]$ obtained by placing an edge between $i$ and $j$ independently for each $i\neq j\in [n]$ with probability 
$1-\exp\big(-p^n_{\lambda}w_i^{\sss (n)}w_j^{\sss (n)}/\ell_n\big)$
for $1\leq i<j\leq n$.
Let $\cC_1^n(\lambda)=\cC\big(1, G_n(\lambda)\big)$. 
\end{defn}

We are now ready to state the result that connects $\tbar M^n$ to $\tbar M^n_{\lambda}$.

\begin{theorem}\label{thm:mst-crit-gh}
Under Assumptions \ref{ass:wts} and \ref{ass:wlog}, for every $\Delta\in (0,1/2]$, there exist $\lambda_{\ref{eqn:555}}=\lambda_{\ref{eqn:555}}(\Delta)\geq 1$, $n_{\ref{eqn:555}}=n_{\ref{eqn:555}}(\Delta)\geq 2$, and $C_{\ref{eqn:555}}=C_{\ref{eqn:555}}(\Delta)>0$ such that for $n\geq n_{\ref{eqn:555}}$ and $\lambda \in [\lambda_{\ref{eqn:555}}, (\nu_n-1)n^{\eta}]$,
\begin{align}\label{eqn:555}
\pr\bigg(d_{\rH}\big(\tbar M^n_{\lambda}\, ,\, \tbar M^n\big) 
\geq 
\frac{n^\eta}{\lambda^{1-\Delta}} + \frac{n^{\eta}}{(\log n)^{1/6}}\bigg)
\leq \frac{C_{\ref{eqn:555}}}{\sqrt{\lambda}}\, . 
\end{align}
\end{theorem}

\begin{rem}
Note that $\cW\big(\cC(1, \barGn)\big) = \Theta_p(n)$, while it turns out (see \eqref{eqn:66}, \eqref{eqn:ss-8}, and Lemma~\ref{lem:1})
that for large (and fixed) $\lambda$, $\cW(\tbar\cC^n_1(\lambda)) = \Theta_P\big(\lambda^{1/(\tau-3)}n^{\rho}\big) = o_P(n)$. 
Thus, the above result shows that despite this major gap in their respective masses, $n^{-\eta}\cdot\tbar M^n_{\lambda}$ approximates $n^{-\eta}\cdot\tbar M^n$ quite well when $\lambda$ is large.
However, as we will see later (see \eqref{eqn:75} and \cite[Theorem 1.2\,(c)]{SB-vdH-SS-PTRF}), if $\lambda$ is kept fixed and $n\to\infty$,
$n^{-\eta}\cdot\tbar M^n_{\lambda}$ converges in distribution in the GH sense to a limiting compact $\bR$-tree whose Minkowski dimension is $(\tau-2)/(\tau-3)$, whereas $\dim(\sM^{\mvtheta^*})=(\tau-1)/(\tau-3)$.
Thus, the intrinsic geometry of $\tbar M^n_{\lambda}$ has features significantly different from that of $\tbar M^n$. 
Similar results were proved in \cite{AddBroGolMie13} in the context of the MST of the complete graph. 
\end{rem}

Theorem \ref{thm:mst-crit-gh} follows upon combining the next two propositions.
\begin{prop}\label{prop:ss-10}
Under Assumptions \ref{ass:wts} and \ref{ass:wlog}, 
there exists $\eps_{\ref{eqn:333}}\in (0, \nu' -1)$, 
and for every 
$\Delta \in (0,1/2]$, there exist $\lambda_{\ref{eqn:555}}=\lambda_{\ref{eqn:555}}(\Delta)\geq 1$, $n_{\ref{eqn:333}}=n_{\ref{eqn:333}}(\Delta)\geq 1$, and $C_{\ref{eqn:333}}=C_{\ref{eqn:333}}(\Delta) >0$ such that for $n\geq n_{\ref{eqn:333}}$ and $\lambda \in [\lambda_{\ref{eqn:555}}, \eps_{\ref{eqn:333}} n^\eta]$, 
\begin{align}\label{eqn:333}
\pr\bigg(d_{\rH}\big(\tbar M^n_{\lambda}\, ,\,  \tbar M^n_{\eps_{\ref{eqn:333}}n^{\eta}}\big) 
\geq 
\frac{n^\eta}{\lambda^{1-\Delta}} \bigg) \leq \frac{C_{\ref{eqn:333}}}{\sqrt{\lambda}}\, .
\end{align}
\end{prop}

\begin{prop}\label{prop:ss-11}
Under Assumptions \ref{ass:wts} and \ref{ass:wlog}, for all large $n$,
\begin{align*}
\pr\bigg(
d_{\rH}\big(\tbar M^n_{\eps_{\ref{eqn:333}}n^{\eta}}\, ,\, \tbar M^n\big) 
\geq 
\big(\log n\big)^{-1/6}n^{\eta} 
\bigg) 
\leq 
n^{-1}.
\end{align*}
\end{prop}

The rest of Section~\ref{sec:7} is devoted to the proofs of Propositions~\ref{prop:ss-10} and \ref{prop:ss-11}.
In Section~\ref{sec:proof-strategy} below, we first explain the rationale behind proving these two results, and then the proofs are completed in several steps in the following sections.
From now on, all results will be proved under Assumptions~\ref{ass:wts} and \ref{ass:wlog}, and we will not mention this explicitly.

\subsection{The general strategy}\label{sec:proof-strategy}
The scaling limit of the MST of the complete graph viewed as a metric measure space was established in \cite{AddBroGolMie13} relying on the results of \cite{addarioberry-broutin-reed,BBG-12,BBG-limit-prop-11}.
This was an important breakthrough, and the limiting space is--quoting the authors of \cite{AddBroGolMie13}--``one of the first scaling limits to be identified for any problem from combinatorial optimisation."
This proof has four key ingredients:
\begin{inparaenumibf}
\item\label{it:i}
One of them is deriving the critical scaling limit of the \erdos random graph.
The scaling limit of the maximal components of the \erdos random graph inside the critcal window was established in \cite{BBG-12} with respect to the GH topology.
This was strengthened to convergence with respect to the GHP topology in \cite{AddBroGolMie13}.
\item\label{it:ii}
Consider the MST $\fM_n$ on the complete graph $K_n$ on $n$ vertices constructed using i.i.d. 
Uniform$[0, 1]$ edge weights $U_{ij}$, $1\leq i<j\leq n$.
Now, for $\lambda>0$, consider the subgraph of $K_n$ with vertex set $[n]$ and edge set
$\big\{
\{i, j\}\, :\, 1\leq i<j\leq n,\ U_{ij}\leq n^{-1}+\lambda n^{-4/3}
\big\}$; 
let $\fM_{n,\lambda}$ denote the restriction of $\fM_n$ to the maximal component in this graph.
An important step in the proof is getting a tail bound on the Hausdorff distance between  $\fM_{n,\lambda}$ and $\fM_n$ for large fixed $\lambda$.
This bound was obtained in \cite{addarioberry-broutin-reed}.
In words, this result states that $\fM_{n,\lambda}$ is quite close to $\fM_n$ in the GH sense if $\lambda$ is sufficiently large, and thus, the structure of $\fM_n$ viewed as a metric space is essentially determined in the late stages of the critical window.
(Note that the authors of \cite{addarioberry-broutin-reed} prove their results in the slightly different setting of a random graph process evolving through the addition of edges in discrete time.
However, this result translates to the setting mentioned above in a straighforward way.)
\item\label{it:iii}
The third ingredient is proving a tail bound on the maximal number of vertices in the trees obtained by removing the edges of $\fM_{n,\lambda}$ from $\fM_n$ for large $\lambda$.
This was established in \cite[Lemma~4.11]{AddBroGolMie13}. 
\item\label{it:iv}
Finally, certain topological properties of the scaling limit were established in \cite{AddBroGolMie13}.
This included showing that the Minkowski dimension of the limiting space is $3$ almost surely; the proof of this result also made use of the results in \cite{BBG-limit-prop-11}.
\end{inparaenumibf}

The scaling limit of $\fM_n$ with respect to the GH topology was established in \cite{AddBroGolMie13} by combining the critical scaling limit of the \erdos random graph mentioned in \textbf{(\ref{it:i})} with the tail bound in \textbf{(\ref{it:ii})} via the results of \cite[Section~3]{AddBroGolMie13}.
This can be strengthened to GHP convergence by using the bound in \textbf{(\ref{it:iii})}  and using the fact that the number of vertices in the trees obtained by removing the edges of $\fM_{n,\lambda}$ from $\fM_n$ satisfy a certain exchangeability propoerty; see the proof of \cite[Proposition~4.8]{AddBroGolMie13} and \cite[Lemma~6.19 and Remark~5]{addarioberry-sen}.

In the context of the multiplicative coalescent in the regime of interest in this paper, the critical metric scaling limit of the maximal components was obtained in \cite{SB-vdH-SS-PTRF,broutin2020limits} (we will need  slightly tweaked versions of these results in our proof as will be discussed in  Section~\ref{sec:existence-scaling-limit} below). 
The missing ingredient for proving the convergence in \eqref{eqn:666} was the analogue of the result in 
\textbf{(\ref{it:ii})} mentioned above in the present setting.
Theorem~\ref{thm:mst-crit-gh} provides this tail bound.
Once we have this bound, it can be combined with the critical scaling limit in a manner similar to \cite{AddBroGolMie13} to deduce the claimed GH convergence in \eqref{eqn:666}.
The techniques used in the proof of Proposition~\ref{prop:ss-10} to analyze the graph $G_n(\lambda)$ will also be useful in Sections~\ref{sec:upper} and \ref{sec:lower-bound-on-dimenension} where we prove the claimed Minkowski dimension in Theorem~\ref{thm:mst-convg}~\eqref{it:c}.

Now, the tail bound mentioned in \textbf{(\ref{it:ii})} above was established in \cite{addarioberry-broutin-reed} in two stages: 
\begin{inparaenuma}
\item\label{it:a-bf}
First, a bound on the Hausdorff distance between  $\fM_{n, n^{1/3}(\log n)^{-1}}$ and $\fM_n$ (going from the barely supercritical regime to the purely supercritical regime) is proved using a variation of Prim's algorithm \cite{prim-algorithm}; this approach used in \cite{addarioberry-broutin-reed} is explained in Algorithm~1 in Section~\ref{sec:45}.
\item\label{it:b-bf}
Next, the Hausdorff distance between $\fM_{n,\lambda}$ and  $\fM_{n, n^{1/3} (\log n)^{-1}}$ (from the critical window to the barely supercritical regime) is bounded for fixed $\lambda>0$.
This step makes use of \cite[Theorem~7]{luczak1990component} which, in words, says that the identity of the maximal component in the \erdos random graph process gets fixed in the late critical window.
The desired bound is then established by considering time points in a suitably chosen geometric progression in the interval $[\lambda, n^{1/3}(\log n)^{-1}]$, and estimating the Hausdorff distance between the MSTs of the maximal components at consecutive time points in this progression.
This technique of proving a property of a random graph process by considering time points where the consecutive points are neither too close nor too far away was previously used in  
\cite[Sections~4 and 6]{luczak1990component} where it was termed the ``scanning method."
\end{inparaenuma}


\ch{Proposition~\ref{prop:ss-11} stated above gives a bound on the Hausdorff distance between 
$\tbar M^n$ and $\tbar M^n_{\eps_{\ref{eqn:333}}n^{\eta}}$--the MST on the component of the vertex $1$ in the purely supercritical regime where the graph is only slightly supercritical.
This gives a result analogous to the one derived in \cite{addarioberry-broutin-reed} (mentioned in \textbf{(\ref{it:a-bf})} above) for the model of interest in this paper.
In our setting, applying Algorithm~1 directly would not yield the desired bound; rather we have to use a modification of this approach (explained in Algorithm~2 in Section~\ref{sec:45}).
Proposition~\ref{prop:ss-10} proves the complementary bound by connecting 
$\tbar M^n_{\lambda}$ to  $\tbar M^n_{\eps_{\ref{eqn:333}}n^{\eta}}$.
In the proof of Proposition~\ref{prop:ss-10}, we use the scanning method as in \cite{addarioberry-broutin-reed,luczak1990component}.
Here, the bulk of the work lies in choosing an appropriate geometric progression to which the scanning method can be applied, and bounding the Hausdorff distance between the MSTs at successive time points in this progression. 
This requires several new techniques.
In particular, we rely on a new method of getting tail bounds on heights of branching processes recently developed in \cite{addario2019most}, concentration inequalities proved in \cite{addario2017probabilistic}, and results on the relation between connected components of inhomogeneous rank-$1$ graphs and $\vp$-trees derived in \cite{SB-vdH-SS-PTRF,SBSSXW14}.
}

\ch{In Section~\ref{sec:size-bias-const} below, we will define two processes called `breadth-first walks' that will help us analyze the random graph models of interest.
In Section~\ref{sec:proofs-mst-crit}, we use a concentration inequality for the suprema of the centered breadth-first walk (which relies on a similar inequality derived in \cite{addario2017probabilistic}) to prove bounds on the lower tails of the total weight of the component of $1$ in $G_n(\lambda)$ as well as the sum of the squares of the vertex weights in the component of $1$ in $G_n(\lambda)$.
These bounds hold for a range of values of $\lambda$ where $G_n(\lambda)$ passes from the late critcal window to the purely supercritical regime.
We also show that when $\lambda$ is in this range, the probability that many of the high-weight vertices are contained in the component of $1$ is lower bounded by an appropriately chosen function of $\lambda$.
In particular, these results will allow us to show that for $\lambda$ in this range, if the component of $1$ is removed from $G_n(\lambda)$, then the rest of the graph is sufficiently subcritical.
In Section~\ref{sec:54}, we obtain a tail bound on the height of a branching process closely related to the random graph $G_n(\lambda)$.
Here, we use a technique developed in \cite{addario2019most}.
Now, as mentioned above, we aim to use the scanning method to prove the claimed bound in Proposition~\ref{prop:ss-10}.
With this in mind and building on the results of Sections~\ref{sec:size-bias-const}, \ref{sec:proofs-mst-crit}, and \ref{sec:54}, 
we achieve the following in Sections~\ref{sec:diam} and \ref{sec:surplus}:
(1) 
We define $\delta_1>0$ such that the scanning method can be applied to a geometric progression of time points where the  common ratio is $(1+\delta_1/2)$.
We also define a random graph $H_n(\lambda, \delta_1)$ such that 
$d_{\rH}\big(
\tbar M^n_{\lambda}\, ,\, \tbar M^n_{\lambda(1+\delta_1/2)}
\big) 
$
can be stochastically bounded in terms of the longest self-avoiding path in $H_n(\lambda, \delta_1)$.
(2)
We prove a tail bound on the diameter of $H_n(\lambda, \delta_1)$.
This is done by establishing height bounds for a random tree that has (potentially) three layers, each of which resembles a multitype branching process.
The offspring distributions in these three layers and the depths of the different layers are chosen in a suitable way to obtain the desired bound.
(3)
We obtain a lower bound for the probability that each component of $H_n(\lambda, \delta_1)$ is either a tree or is unicyclic.
Here, we make use of a construction of a connected component of a rank-$1$ inhomogeneous random graph using $\vp$-trees \cite{SB-vdH-SS-PTRF,SBSSXW14}.
We then use these results and apply the scanning method to complete the proof of Proposition~\ref{prop:ss-10} in Section~\ref{sec:scanning}.
Finally, the proof of Proposition~\ref{prop:ss-11} is given in Section~\ref{sec:45}.
As mentioned above, here we use a modification of the approach used in the proof of \cite[Lemma~4]{addarioberry-broutin-reed}.
}

\subsection{An exploration process}\label{sec:size-bias-const}
Fix $\lambda\geq 0$.
It will be useful in the proof to express $G_n(\lambda)$ in a reparametrized form. 
Define
\begin{equation}
\label{eqn:xin-def}
x_i^{\sss(n)}:= \frac{w_i^{\sss(n)}}{n^\rho \big(\sigma_2^{\sss(n)}\big)^{1/2}}\, ,\ \ \text{ and }\ \ \ \theta_{i,\lambda}^{\sss(n)} := \bigg(1+\frac{\lambda}{n^\eta}\bigg) \frac{w_i^{\sss(n)}}{n^\alpha \big(\sigma_2^{\sss(n)}\big)^{1/2}}\, , \qquad 1\leq i\leq n. 
\end{equation}
Write $\mvx^{\sss(n)}=\big(x_i^{\sss(n)}\, ,\, i\in[n]\big)$, and
$\sigma_2(\mvx^{\sss(n)}):=\sum_{i=2}^n (x_i^{\sss(n)})^2$. 
(Note that similar to \eqref{eqn:sigm2-nu-def}, $x_1^{\sss(n)}$ is not included in the sum.)
Then
\begin{equation}
\label{eqn:ss-1}
\sigma_2(\mvx^{\sss(n)}) = n^{-\eta}\, , \ \ \text{ and }\ \ 
\big(\lambda + n^{\eta}\big) x_j^{\sss(n)} = \theta_{j,\lambda}^{\sss(n)} \ \ \text{ for }\ \  j\in [n]\, . 
\end{equation}
Further, by Assumption \ref{ass:wts}, for each $i\geq 1$,
\begin{equation}
\label{eqn:ss-2}
\big(\sigma_2(\mvx^{\sss(n)})\big)^{-1}x_i^{\sss(n)}
=n^{\eta}\cdot x_i^{\sss(n)}
\to \theta_i^*\, ,\ \ \text{ as }\ \ n\to\infty\, .
\end{equation}
Now note that for $2\leq i<j\leq n$, 
\begin{equation}
\label{eqn:711}
p^n_{\lambda}\bigg(\frac{w_i^{\sss(n)} w_j^{\sss(n)}}{\ell_n}\bigg) 
= 
\bigg(\lambda + \frac{1}{\sigma_2(\mvx^{\sss(n)})}\bigg)x_i^{\sss(n)} x_j^{\sss(n)} \, ,
\end{equation}
and consequently, 
\begin{align}\label{eqn:712}
G_n(\lambda)
\equald
\cG_n(\lambda):=
\cG\bigg(\big([n], \mvx^{\sss(n)}\big)\, ,\,  \lambda+\big(\sigma_2(\mvx^{\sss(n)})\big)^{-1}\bigg)\, ,
\end{align}
where the latter is as in Definition \ref{def:1}.

Let $\cG_n^-(\lambda)=\cG_n(\lambda)\setminus [1]$.
A useful tool in the study of the random graph $\cG_n^-(\lambda)$ is the ``breadth-first walk" process $\big(Z^{n, -}_{\lambda}(u), u\geq 0\big)$ associated with a breadth-first exploration of the random graph $\cG_n^-(\lambda)$, which we describe next.
This is very much in the spirit of \cite{aldous-limic}, although our breadth-first walk (defined in \eqref{eqn:7777}) is slightly different from the one considered in \cite{aldous-limic}, as it is easier to analyze.

For $2\leq j\leq n$, let the size of vertex $j$ be $x_j^{\sss (n)}$.
Choose $v(1)$ from $[n]\setminus\{1\}$ in a size-biased way, i.e., $\pr\big(v(1)=v\big)\propto x_v^{\sss(n)}$, $v\in [n]\setminus\{1\}$.
Explore the component of $v(1)$ in $\cG_n^-(\lambda)$ in a breadth-first fashion. 
Let $h^1$ be the height of the breadth-first tree, and for $0\leq i\leq h^1$, let $\text{Gen}^1_i$ be the set of vertices in the $i$-th generation of the breadth-first tree.
For $t\geq 2$, having explored the component of $v(1),\ldots, v(t-1)$, 
choose $v(t)$ in a size-biased way from the remaining vertices, 
explore its component in $\cG_n^-(\lambda)$ in a breadth-first fashion, 
let $h^t$ be the height of its breadth-first tree, and for $0\leq i\leq h^t$, let $\text{Gen}^t_i$ be the set of vertices in the $i$-th generation of the breadth-first tree.
Stop when all vertices $j\in [n]\setminus\{1\}$ have been found.

Using properties of exponential random variables, it can be easily checked that the above collection of random variables can be constructed in the following way: 
Recall the relations from \eqref{eqn:ss-1}.
Let $\xi_j^n$, $2\leq j\leq n$, be independent random variables such that
\begin{align}\label{eqn:713}
\xi_{j, \lambda}^n
\sim
\EXP\big(\big( \lambda+ \big(\sigma_2(\mvx^{\sss(n)})\big)^{-1}\big) x_j^{\sss(n)}\big)
\, ,\
\ch{\text{ or equivalently, }\
\xi_{j, \lambda}^n
\sim
\EXP\big(\theta_{j,\lambda}^{\sss(n)}\big)\, .
}
\end{align}
To simplify notation we will write $\xi_j^n$ instead of $\xi_{j, \lambda}^n$.
Let $v(1)$ be such that $\xi_{v(1)}^n=\min\big\{\xi_j^n\, :\, 2\leq j\leq n\big\}$, and set
$
\text{Gen}^1_0=\big\{v(1)\big\}.
$
Inductively define
\[
\text{Gen}^1_i
=
\bigg\{ j\in [n]\setminus\{1\}\, :\, 
\xi_j^n\in
\bigg[\xi_{v(1)}^n+\sum_{k=0}^{i-2}\sum_{v\in\text{Gen}^1_k}x_v^{\sss (n)}\, ,\ 
\xi_{v(1)}^n+\sum_{k=0}^{i-1}\sum_{v\in\text{Gen}^1_k}x_v^{\sss (n)}\bigg]
\bigg\}\text{ for } 1\leq i\leq h^1+1\, ,
\]
where 
\[
h^1+1=\min\, \big\{ i\geq 1\, :\, \text{Gen}^1_i=\emptyset\big\}\, .
\]
For $t\geq 2$, let $v(t)$ be such that 
\[
\xi_{v(t)}^n=
\min\bigg\{\xi_j^n\, :\, 2\leq j\leq n\ \text{ and }\ 
j\notin \bigcup_{k=1}^{t-1}\bigcup_{i=0}^{h^k}\text{Gen}^k_i\bigg\}\, .
\] 
Set
$\text{Gen}^t_0=\big\{v(t)\big\}$,
and define $\text{Gen}^t_i$, $1\leq i\leq h^t+1$, in a manner analogous to the case $t=1$.
Stop when all vertices $j\in [n]\setminus\{1\}$ have been found.

We define the breadth-first walk as
\begin{align}\label{eqn:7777}
Z^{n, -}_{\lambda}(u) 
= 
- u+ \sum_{2\leq j\leq n} x_j^{\sss(n)} \ind\big\{\xi_j^n \leq u\big\}\, ,\ \ u\geq 0\, .
\end{align}
The correspondence described above allows one to prove various properties of $\cG_n^-(\lambda)$ by studying the process $Z^{n, -}_{\lambda}$. 
Here we make note of an elementary property of $Z^{n, -}_{\lambda}$ that will be useful to us:
Suppose $\cG_n^-(\lambda)$ and $\xi_j^n$, $2\leq j\leq n$, are coupled by means of the correspondence described above. Then in this coupling, for $1\leq t\leq m$, 
\[
Z^{n, -}_{\lambda}\big(\xi_{v(t)}^n-\big)
= 
Z^{n, -}_{\lambda}\bigg(\xi_{v(t)}^n
+\cW_{\mvx^{\sss (n)}}\bigg(\cC\big(v_{(t)}, \cG_n^-(\lambda)\big)\bigg)\bigg)
=
-\xi_{v(t)}^n+\sum_{k=1}^{t-1}\cW_{\mvx^{\sss (n)}}\bigg(\cC\big(v_{(k)}, \cG_n^-(\lambda)\big)\bigg) 
<0\, ,
\]
where $m$ denotes the number of components in $\cG_n^-(\lambda)$.
This leads to the following:
\begin{lemma}\label{lem:domin-ald-limic}
In the above coupling, if for some $0<u_1<u_2$, $Z^{n, -}_{\lambda}(u)>0$ for $u\in [u_1, u_2]$, then there exists a component $\cC^{\star}$ in $\cG_n^-(\lambda)$ such that 
$\cW_{\mvx^{\sss (n)}}\big(\cC^{\star}\big)\geq u_2-u_1$.
\end{lemma}

Note that we can explore $\cG_n(\lambda)$ starting from the vertex $1$ in a manner similar to the exploration of $\cG_n^-(\lambda)$.
In this case, we define the breadth-first walk as
\begin{align}\label{eqn:8888}
Z^{\sss n, (1)}_{\lambda}(u) 
= 
x_1^{\sss(n)}- u+ \sum_{2\leq j\leq n} x_j^{\sss(n)} \ind\big\{\xi_j^n \leq u\big\}\, ,\ \ u\geq 0\, ,
\end{align}
where $\xi_j^n$, $2\leq j\leq n$, are as in \eqref{eqn:713}.
(Here we append a ``$(1)$'' to specify that the process starts at vertex $1$.)
As in the case of $\cG_n^-(\lambda)$, there is a natural coupling between $\cG_n(\lambda)$ and $\xi_j^n$, $2\leq j\leq n$.
\begin{lemma}\label{lem:domin-ald-limic-1}
	In the above coupling, 
	\[
	\cW_{\mvx^{\sss (n)}}\big(\cC\big(1, \cG_n(\lambda)\big)\big)
	=
	\inf\, \big\{u\geq 0\, :\,  Z^{\sss n, (1)}_{\lambda}(u)=0\big\}\, .
	\]
\end{lemma}

Let us get back to the process $Z^{n, -}_{\lambda}(\cdot)$.
Using \eqref{eqn:713} and the second identity in \eqref{eqn:ss-1}, it can be directly checked that for any $u\geq 0$,
\begin{equation}
\label{eqn:ss-3}
\E\big[n^{\eta}Z^{n, -}_{\lambda}(u)\big] = \lambda u - \sum_{j=2}^n \frac{\theta_{j,\lambda}^{\sss(n)}}{\big(1+\lambda n^{-\eta}\big)}\left(u\theta_{j,\lambda}^{\sss(n)} + \exp(-u\theta_{j,\lambda}^{\sss(n)}) -1\right) =: \Phi_{\lambda}^{\sss(n)}(u). 
\end{equation}
Note that $\Phi_{\lambda}^{\sss(n)}(\cdot)$ is a strictly concave function (in $u$).
In particular, for any $\lambda>0$, $\Phi_{\lambda}^{\sss(n)}(\cdot)$ has a unique positive zero, which we will denote by $s^{\sss(n)}(\lambda)$.
Define 
\begin{equation}
\label{eqn:ss-4}
\varphi_{\lambda}^{\sss(n)}(u) := \sum_{j=2}^n \big(\theta_{j,\lambda}^{\sss(n)}\big)^2\frac{\left(u \theta_{j,\lambda}^{\sss(n)} + \exp(- u\theta_{j,\lambda}^{\sss(n)}) -1\right)}{u \theta_{j, \lambda}^{\sss(n)}}
\end{equation}
for $u\geq 0$ (the value at $u=0$ is understood to be the limit of $\varphi_{\lambda}^{\sss(n)}(u)$ as $u\downarrow 0$), so that 
\begin{equation}
\label{eqn:ss-5}
\Phi_{\lambda}^{\sss(n)}(u) = \lambda u  - u \varphi_{\lambda}^{\sss(n)}(u)\big/\big(1+\lambda n^{-\eta}\big)\, .
\end{equation}

From now on, we will drop the dependence on $n$ in the superscripts and simply write $\theta_{j,\lambda}$, $x_j$, $\mvx$, and $w_j$ to ease notation. 
Recall the constants $A_1, A_2$ from Assumption \ref{ass:wts}.  Consider the interval 
\begin{equation}
\label{eqn:inter-def}
\dI^{\sss(n)}:= 
\bigg[\frac{2\cdot\big(\sigma_2^{\sss(n)}\big)^{1/2}}{A_1}\, ,\ \frac{n^{\alpha}\big(\sigma_2^{\sss(n)}\big)^{1/2}}{A_2 2^{\alpha+1}}\bigg]\, .
\end{equation}
For $u \in \dI^{\sss(n)}$, define $i_\lambda(u) = \min\set{i\geq 1: \theta_{i,\lambda} u <1}$. 
From the definition of $\dI^{\sss(n)}$ and Assumption \ref{ass:wts}~(iii), it follows that $3\leq i_{\lambda}(u) \leq n/2$ for $u \in \dI^{\sss(n)}$. 
Writing 
\begin{align}\label{eqn:111}
\dI^{\sss(n),2} = \big\{(\lambda, u)\, :\, 0\leq \lambda \leq  n^{\eta}/10\, ,\  u \in \dI^{\sss(n)} \big\}\, ,
\end{align}
we have, for all large $n$,
\begin{equation}
\label{eqn:ss-6}
i_\lambda(u) \asymp u^{1/\alpha}\ \ \text{ for }\ \ (\lambda, u)\in\dI^{\sss(n),2} \, .
\end{equation}

Define the function $g(s):= s+e^{-s} -1$ for $s\geq 0$. Since $ g(s) \asymp s^2$ on $[0,1]$ while $g(s) \asymp s$ on $[1,\infty)$, using Assumption \ref{ass:wts} we have, $\falln$ and for all $(\lambda, u) \in \dI^{\sss(n),2}$, 
\begin{align}\label{eqn:ss-88}
u\varphi_{\lambda}^{\sss(n)}(u) 
&
\asymp 
\sum_{j=2}^{i_{\lambda}(u) -1} \theta_{j,\lambda}\big(u\theta_{j,\lambda}\big) 
+ 
\sum_{j=i_\lambda(u)}^n \theta_{j,\lambda}\big(u\theta_{j,\lambda}\big)^2 
\asymp 
u \sum_{j=2}^{i_\lambda(u) -1} \frac{1}{j^{2\alpha}} 
+ 
u^2 \sum_{j=i_\lambda(u)}^n\big(\theta_{j,\lambda}\big)^3\, .
\end{align}
Using \eqref{eqn:ss-6}, we see that $\falln$ and for all $(\lambda, u) \in \dI^{\sss(n),2}$, 
\begin{align}\label{eqn:ss-88-a}
u\sum_{j=2}^{i_{\lambda}(u) -1} \frac{1}{j^{2\alpha}}
\asymp
u\cdot\big(i_{\lambda}(u)\big)^{1-2\alpha} 
\asymp
u\cdot u^{(1-2\alpha)/\alpha}
\asymp
u^{\tau-2}\, ,
\end{align}
whereas Assumption \ref{ass:wts} (iii) and \eqref{eqn:up-bnd} yield
\begin{gather}
0
\leq 
u^2 \sum_{j=i_{\lambda}(u)}^n\big(\theta_{j,\lambda}\big)^3
\leq
Cu^2 
\sum_{j=i_{\lambda}(u)}^n \frac{1}{j^{3\alpha}}
\leq
C'u^2\big(i_{\lambda}(u)\big)^{1-3\alpha} 
\leq
C'' u^{\tau-2}\, .\label{eqn:ss-88-b}
\end{gather}
Combining \eqref{eqn:ss-88}, \eqref{eqn:ss-88-a}, and \eqref{eqn:ss-88-b}, we get, $\falln$ and for all $(\lambda, u) \in \dI^{\sss(n),2}$, 
\begin{align}\label{eqn:ss-7}
u\varphi_{\lambda}^{\sss(n)}(u) 
\asymp 
u^{\tau -2}\, .
\end{align}
Now let us switch back to $\Phi_{\lambda}^{\sss(n)}(\cdot)$. 
Using the asymptotics for $\varphi_{\lambda}^{\sss(n)}(\cdot)$ in \eqref{eqn:ss-7}, choose $\lambda_{\ref{eqn:ss-7a}} \geq 1$ large enough and $\eps_{\ref{eqn:ss-7a}}  \in \big(0,\frac{1}{10}\wedge(\nu' -1)\big)$ small 
so that $\falln$,
\[
\Phi_{\lambda_{\ref{eqn:ss-7a}}}^{\sss (n)}
\bigg( \frac{2\big(\sigma_2^{\sss(n)}\big)^{1/2}}{A_1}\bigg)>0\, , \ \ 
\text{ and }\ \ 
\Phi_{\eps_{\ref{eqn:ss-7a}}n^{\eta}}^{\sss (n)}
\bigg(\frac{n^{\alpha}\big(\sigma_2^{\sss(n)}\big)^{1/2}}{A_2 2^{\alpha+1}}\bigg)<0\, ,
\]
and consequently, the unique positive zero $s^{\sss(n)}(\lambda)$ of $\Phi_{\lambda}^{\sss(n)}(\cdot)$ satisfies, $\falln$,
\begin{equation}
\label{eqn:ss-7a}
\qquad s^{\sss(n)}(\lambda) \in \dI^{\sss(n)} 
\qquad \text{ for }\ \ \lambda \in [\lambda_{\ref{eqn:ss-7a}} ,\, \eps_{\ref{eqn:ss-7a}}  n^\eta]. 
\end{equation}
(Note that  
$\eps_{\ref{eqn:ss-7a}}  \in \big(0,\frac{1}{10}\wedge(\nu'-1)\big)$ 
implies that $\falln$, $\eps_{\ref{eqn:ss-7a}}  < (\nu_n-1)$.)
Thus, \eqref{eqn:ss-7}, \eqref{eqn:ss-5}, and the relation $\Phi_{\lambda}^{\sss (n)}\big(s^{\sss(n)}(\lambda)\big)=0$ implies, $\falln$,
\begin{equation}
\label{eqn:ss-8}
s^{\sss(n)}(\lambda) \asymp \lambda^{1/(\tau-3)} 
\qquad \text{ for }\  \lambda \in [\lambda_{\ref{eqn:ss-7a}} ,\, \eps_{\ref{eqn:ss-7a}}  n^\eta]. 
\end{equation}

\subsection{The component of the vertex $1$ in $G_n(\lambda)$}
\label{sec:proofs-mst-crit}
Recall the notation $\cC_1^n(\lambda)$ from Definition~\ref{def:graph-vgexp}.
In this section we \ch{will} study properties of $\cC_1^n(\lambda)$.
We start with a lower bound on $\cW\big(\cC_1^n(\lambda)\big)$.

\begin{prop}\label{prop:lem-ss-3}
There exists $\kappa_0\in (0,1/4)$ 
such that the following holds for all $\kappa\in (0,\kappa_0]$: there exists $n_{\ref{eqn:66}} = n_{\ref{eqn:66}}(\kappa)$ and $\lambda_{\ref{eqn:66}}=\lambda_{\ref{eqn:66}}(\kappa) \geq \lambda_{\ref{eqn:ss-7a}}$ such that for all $n\geq n_{\ref{eqn:66}}$ and $\lambda \in [\lambda_{\ref{eqn:66}} ,\, \eps_{\ref{eqn:ss-7a}} n^\eta]$,
\begin{align}\label{eqn:66}
\pr\bigg(\cW_{\mvx}\big(\cC\big(1, \cG_n(\lambda)\big)\big) 
\geq 
(1-2\kappa) s^{\sss(n)}(\lambda)\bigg) \geq 
1-\exp\big(-C\lambda^{1/(\tau-3)}\big)\, ,
\end{align}
and consequently, 
$\pr\bigg(\cW\big(\cC_1^n(\lambda)\big) 
\geq 
(1-2\kappa) s^{\sss(n)}(\lambda) \big(\sigma_2^{\sss(n)}\big)^{1/2} n^{\rho}\bigg)
\geq
1-\exp\big(-C\lambda^{1/(\tau-3)}\big)$.
 \end{prop}

We will make use of the next two results in the proof of Proposition \ref{prop:lem-ss-3}. 
Recall the independent exponential random variables $\xi^n_j$, $2\leq j\leq n$, from \eqref{eqn:713}.
\begin{lemma}[{\cite[Lemma 6.3]{addario2017probabilistic}}]\label{lem:ss-lem-3A}
There exists a constant $C_{\ref{eqn:99}}>0$ such that for all $s\geq 1/C_{\ref{eqn:99}}$, $y\geq 1/C_{\ref{eqn:99}}$, for all large $n$, and $\lambda \in [1,n^{\eta}]$,
\begin{align}\label{eqn:99}
\pr\bigg(
\sup_{u\leq s}\bigg|\sum_{j=2}^n \theta_{j,\lambda}\bigg(\ind\big\{\xi^n_j\leq u\big\} 
- 
\pr\big(\xi^n_j\leq u\big)\bigg)\bigg| \geq ys^{\frac{\tau-3}{2}}\bigg) 
\leq 
\exp\big(-C_{\ref{eqn:99}} y\log\log{y}\big)\, .
\end{align}
\end{lemma}

\begin{rem}
The result in \cite[Lemma 6.3]{addario2017probabilistic} is given in a slightly different setting.
However, the key ingredient in its proof is the Klein-Rio bound \cite[Theorem 1.1]{klein-rio05}.
Applying \cite[Theorem 1.1]{klein-rio05}, the proof of \cite[Lemma 6.3]{addario2017probabilistic} boils down to establishing a uniform upper bound on the variance of a certain collection of functionals, and an upper bound on the expectation of the supremum of the said collection of functionals.
This is achieved in the proof of \cite[Lemma 6.3]{addario2017probabilistic} and in \cite[Lemma 6.4]{addario2017probabilistic}.
Now, using Assumption \ref{ass:wts}\,(iii) and \eqref{eqn:up-bnd}, those same arguments can be used to prove Lemma \ref{lem:ss-lem-3A}.
Further, examining the proofs of \cite[Lemma 6.3 and Lemma 6.4]{addario2017probabilistic} will reveal that the constant $C_{\ref{eqn:99}}$ can be chosen so that it depends only on $A_1$, $A_2$, and $\tau$.
We omit the proof of Lemma \ref{lem:ss-lem-3A} as no new idea is involved.
\end{rem}

The next lemma describes technical properties of the function $\varphi_\lambda^{\sss(n)}$ defined  in \eqref{eqn:ss-4}. 
For fixed $\kappa>0$, define analogous to the set $\dI^{\sss(n),2}$, 
\[
\dI^{\sss(n),2}_\kappa:= 
\bigg\{
(\lambda, u)\, :\, 
0\leq \lambda\leq  n^{\eta}/10\, ,\ 
u\in \bigg[\frac{2\cdot\big(\sigma_2^{\sss(n)}\big)^{1/2}}{A_1 \kappa}\, ,\ \frac{n^{\alpha}\cdot\big(\sigma_2^{\sss(n)}\big)^{1/2}}{A_2 2^{\alpha+1}}\bigg]\, 
\bigg\}\, .
\]

\begin{lemma}\label{lem:varphi-prop}
There exist $C_{\ref{eqn:55}}>0$ and $\kappa_0 \in (0,1/4)$ such that for all large $n$, for all $\kappa \in (0,\kappa_0]$, and $(\lambda, u) \in \dI_\kappa^{\sss(n),2}$,
\begin{align}\label{eqn:55}
\frac{\varphi_\lambda^{\sss(n)}(\kappa u)}{\varphi_\lambda^{\sss(n)}(u)} \leq \frac{\varphi_\lambda^{\sss(n)}\big((1-\kappa)u\big)}{\varphi_\lambda^{\sss(n)}(u)} \leq \frac{1}{1+C_{\ref{eqn:55}}\kappa}\, .
\end{align}
\end{lemma}

\noindent {\bf Proof of Proposition \ref{prop:lem-ss-3} assuming Lemma \ref{lem:varphi-prop}:} 
Recall the function $\Phi_{\lambda}^{\sss(n)}$ from \eqref{eqn:ss-3} and its connection to $\varphi_\lambda^{\sss(n)}$ from \eqref{eqn:ss-5}.
By \eqref{eqn:ss-7a} and \eqref{eqn:ss-8}, for any $\kappa\in(0, \kappa_0]$, there exist $\lambda_{\ref{eqn:9A}} = \lambda_{\ref{eqn:9A}}(\kappa)\geq \lambda_{\ref{eqn:ss-7a}}$ and $n_{\ref{eqn:9A}} = n_{\ref{eqn:9A}}(\kappa)$ such that 
$\big(\lambda, s^{\sss(n)}(\lambda) \big)\in \dI^{\sss(n), 2}_\kappa$ for all $\lambda \in [\lambda_{\ref{eqn:9A}},\, \eps_{\ref{eqn:ss-7a}} n^\eta]$ and $n\geq n_{\ref{eqn:9A}}$. 
Then for $n\geq n_{\ref{eqn:9A}}$ and $\lambda \in [\lambda_{\ref{eqn:9A}},\, \eps_{\ref{eqn:ss-7a}} n^{\eta}]$, 
\begin{align}
	\Phi_\lambda^{\sss(n)}\big((1-\kappa) s^{\sss(n)}(\lambda)\big) 
	&\geq (1-\kappa) s^{\sss(n)}(\lambda)\bigg(\lambda - \frac{1}{1+\lambda n^{-\eta}}\cdot \frac{\varphi_\lambda^{\sss(n)}(s^{\sss(n)}(\lambda))}{1+C_{\ref{eqn:55}}\kappa}\bigg)\notag \\
	&= (1-\kappa) s^{\sss(n)}(\lambda)\cdot\lambda\cdot \bigg(\frac{\kappa C_{\ref{eqn:55}}}{1+\kappa C_{\ref{eqn:55}}}\bigg) \geq C_{\ref{eqn:9A}}\kappa \lambda^{\frac{\tau-2}{\tau-3}}\, , \label{eqn:9A}
	\end{align}
	where the first inequality uses \eqref{eqn:ss-5} and Lemma \ref{lem:varphi-prop},
	the second step follows from the definition of $s^{\sss(n)}(\lambda)$ which implies that $\varphi_\lambda^{\sss(n)}(s^{\sss(n)}(\lambda)) = \lambda \cdot (1+\lambda n^{-\eta})$, and the final step follows from \eqref{eqn:ss-8}. 
	Similarly, for all $n\geq n_{\ref{eqn:9A}}$ and $\lambda \in [\lambda_{\ref{eqn:9A}},\, \eps_{\ref{eqn:ss-7a}}n^{\eta}]$, 
	\begin{equation}
	\label{eqn:9B}
		\Phi_\lambda^{\sss(n)}\big(\kappa s^{\sss(n)}(\lambda)\big) \geq C_{\ref{eqn:9A}}\kappa^2 \lambda^{(\tau-2)/(\tau-3)}.  
	\end{equation}
Now recall the process $Z_{\lambda}^{n, -}$ from \eqref{eqn:7777}, and note that for any $u\geq 0$,
\begin{equation}
\label{eqn:9D}
n^{\eta}\bigg(Z_{\lambda}^{n, -}(u) - \E\big[Z_{\lambda}^{n, -}(u)\big]\bigg) 
= 
\sum_{j=2}^n \frac{\theta_{j,\lambda}}{(1+\lambda n^{-\eta})}
\bigg(\ind\big\{\xi^n_j\leq u\big\} - \pr\big(\xi^n_j \leq u\big)\bigg)\, .
\end{equation}
Thus, by Lemma \ref{lem:ss-lem-3A}, there exists $\lambda_{\ref{eqn:9C}}\geq\lambda_{\ref{eqn:ss-7a}}$ such that $\falln$ and $\lambda\in[\lambda_{\ref{eqn:9C}},\, \eps_{\ref{eqn:ss-7a}}n^{\eta} ]$,
\begin{equation}
	\label{eqn:9C}
		\pr\bigg(\sup_{u\leq s^{\sss(n)}(\lambda)}
		\big(\lambda+n^{\eta}\big)\bigg|Z_{\lambda}^{n, -}(u) - \E\big[Z_{\lambda}^{n, -}(u)\big]\bigg|
		\geq 
		\lambda^{\frac{1}{\tau-3}}\big(s^{\sss(n)}(\lambda)\big)^{\frac{\tau-3}{2}}\bigg)\leq \exp\big(-C\lambda^{\frac{1}{\tau-3}}\log\log \lambda^{\frac{1}{\tau-3}}\big). 
	\end{equation}
Using \eqref{eqn:ss-8}, 
$\falln$ and $\lambda\in[\lambda_{\ref{eqn:9C}},\, \eps_{\ref{eqn:ss-7a}}n^{\eta} ]$,
$\lambda^{\frac{1}{\tau -3}} \big(s^{\sss(n)}(\lambda)\big)^{\frac{\tau-3}{2}} 
\leq C \lambda^{\frac{\tau-1}{2(\tau-3)}}$. 
Now using \eqref{eqn:9A}, \eqref{eqn:9B}, and \eqref{eqn:9C} together with the concavity of $\Phi_{\lambda}^{\sss(n)}$ and the fact that $(\tau-1)/(2(\tau-3) )< (\tau-2)/(\tau-3)$, 
we can find
$\lambda_{\ref{eqn:338}} = \lambda_{\ref{eqn:338}}(\kappa)\geq\lambda_{\ref{eqn:ss-7a}}$ and $n_{\ref{eqn:338}} = n_{\ref{eqn:338}}(\kappa)$ such that for all $n\geq n_{\ref{eqn:338}}$ and $\lambda \in [\lambda_{\ref{eqn:338}},\, \eps_{\ref{eqn:ss-7a}} n^\eta]$,
\begin{equation}
\label{eqn:338}
\pr\bigg(
Z_{\lambda}^{n, -}(u)>0\ \text{ for every }\ 
u\in\big[ \kappa s^{\sss(n)}(\lambda), (1-\kappa) s^{\sss(n)}(\lambda) \big] 
\bigg)
\geq
1-\exp\big(-C\lambda^{\frac{1}{\tau-3}}\big)\, .
\end{equation}

By Lemma \ref{lem:domin-ald-limic}, on the event in \eqref{eqn:338}, there exists a component $\cC^{\star}$ in $\cG_n^-(\lambda)$ with 
$\cW_{\mvx}(\cC^{\star})\geq (1-2\kappa)s^{\sss(n)}(\lambda) $. 
Now we can generate $\cG_n(\lambda)$ by first generating $\cG_n^-(\lambda)$, and then independently sampling the edges from vertex $1$ to the vertex set $[n]\setminus\{1\}$. 
Thus, for all $n\geq n_{\ref{eqn:338}}$ and $\lambda \in [\lambda_{\ref{eqn:338}},\, \eps_{\ref{eqn:ss-7a}} n^\eta]$,
\begin{align}
\pr\bigg(\sW_{\mvx}\big(\cC\big(1, \cG_n(\lambda)\big)\big)
\leq (1-2\kappa) s^{\sss(n)}(\lambda)\bigg) 
\leq \exp\bigg(-C\lambda^{\frac{1}{\tau-3}}\bigg)
+ 
\exp\bigg(-\big(\lambda+n^{\eta}\big)x_1\cdot (1-2\kappa)s^{\sss(n)}(\lambda)\bigg)\, , \notag
\end{align}
where the second term is an upper bound for the probability that vertex $1$ is {\bf not} connected to the component $\cC^{\star}$ of mass at least $(1-2\kappa) s^{\sss(n)}(\lambda)$. 
Now an application of \eqref{eqn:ss-8} and \eqref{eqn:ss-2} completes the proof of Proposition \ref{prop:lem-ss-3}. 
\qed

\medskip

Let us now turn to the proof of Lemma \ref{lem:varphi-prop}.

\medskip

\noindent {\bf Proof of Lemma \ref{lem:varphi-prop}:}  
Recall the definition of $\dI^{\sss(n),2}$ from \eqref{eqn:111}.
Let $f:[0, \infty)\to [0, 1]$ be given by $f(u) := 1-e^{-u} - ue^{-u}$.  
Note that $f(u)\asymp 1$ for $u\geq 1$ while $f(u) \asymp u^2$ for $u\in [0,1]$. 
Thus, $\falln$ and for $(\lambda, u)\in \dI^{\sss(n),2}$, 
\begin{align}\label{eqn:10}
\sum_{j=2}^n \theta_{j,\lambda}
\big(1-e^{-u\theta_{j,\lambda}} - u\theta_{j,\lambda}e^{-u\theta_{j,\lambda}}\big)
\asymp 
\sum_{j=2}^{i_{\lambda}(u) -1}\theta_{j,\lambda}
+
\sum_{j=i_{\lambda}(u)}^n \theta_{j,\lambda}\big(u\theta_{j,\lambda}\big)^2
\asymp 
u^{\tau-2},
\end{align} 
where the last step uses arguments similar to the ones leading to \eqref{eqn:ss-7}. 
Now for any $u\geq 0$, $\kappa\in (0,1)$, and $s\in [1-\kappa, 1]$, 
\begin{align}\label{eqn:100}
\big|\sum_{j=2}^n \theta_{j,\lambda}
\big(u\theta_{j,\lambda}e^{-us\theta_{j,\lambda}} - u\theta_{j,\lambda}e^{-u\theta_{j,\lambda}}\big)\big|
\leq 
\sum_{j=2}^n 
u\theta_{j,\lambda}^2\big(1-e^{-\kappa u \theta_{j,\lambda}}\big)\, .
\end{align}
For all large $n$ and $(\lambda, \kappa u) \in \dI^{\sss(n),2}$, 
\begin{align}\label{eqn:10A}
&
\sum_{j=2}^n \theta_{j,\lambda}^2 u \big(1-e^{-\kappa u\theta_{j,\lambda}}\big) 
\asymp 
\sum_{j=2}^{i_\lambda(\kappa u)-1} \theta_{j,\lambda}^2u 
+ 
\sum_{j=i_\lambda(\kappa u)}^n \theta_{j,\lambda}^2 u \big(\kappa u\theta_{j,\lambda}\big) 
\notag \\ 
&\hskip30pt
\asymp 
u \big(i_\lambda(\kappa u)\big)^{1-2\alpha} 
+ 
\kappa u^2\sum_{j=i_{\lambda}(\kappa u)}^{n}\big(\theta_{j,\lambda}\big)^3
\asymp 
\kappa^{\tau -3} u^{\tau -2}\, ,
\end{align}
where we have used Assumption \ref{ass:wts}\,(iii), \eqref{eqn:up-bnd}, and \eqref{eqn:ss-6}.
Combining \eqref{eqn:10}, \eqref{eqn:100}, and \eqref{eqn:10A}, we see that there exists $ \kappa_0\in (0,1/4)$ small such that $\falln$ and for all $\kappa \in (0,\kappa_0]$,
\begin{equation}
\label{eqn:10B}
	\sum_{j=2}^n \theta_{j,\lambda} 
	\big(1-e^{-u\theta_{j,\lambda}} - u\theta_{j,\lambda}e^{-us\theta_{j,\lambda}}\big) 
	\asymp u^{\tau -2}
	\ \ \ \text{ for }\ \ \ (\lambda, u)\in\dI_\kappa^{\sss(n), 2}
\end{equation}
uniformly over $s\in [1-\kappa, 1]$. 
Hence, $\falln$, for any $\kappa \in (0,\kappa_0]$, and $(\lambda, u) \in \dI^{\sss(n),2}_\kappa$, 
\begin{align}
\frac{\varphi_\lambda^{\sss(n)}(u)}{\varphi_\lambda^{\sss(n)}\big((1-\kappa)u\big)} -1 
&=
\frac{\sum_{j=2}^n \theta_{j,\lambda}\big(\kappa\big(1-e^{-u\theta_{j,\lambda}}\big) + e^{-u\theta_{j,\lambda}} - e^{-(1-\kappa)u\theta_{j,\lambda}}\big)}{(1-\kappa)u\cdot\varphi_\lambda^{\sss(n)}\big((1-\kappa)u\big)} \notag \\
&\geq 
\frac{1}{C u^{\tau -2}}\cdot
\bigg(\sum_{j=2}^n \theta_{j,\lambda}\bigg(\kappa\big(1-e^{-u\theta_{j,\lambda}}\big) + e^{-u\theta_{j,\lambda}} - e^{-(1-\kappa)u\theta_{j,\lambda}}\bigg)\bigg)\notag\\
&=
\frac{1}{C u^{\tau -2}}\cdot
\bigg(
\sum_{j=2}^n \theta_{j,\lambda} \int_{1-\kappa}^1 \big(1-e^{-u\theta_{j,\lambda}} - u\theta_{j,\lambda}e^{-u\theta_{j,\lambda}s}\big)ds\bigg)
\geq 
C_{\ref{eqn:55}}\kappa\, ,\notag
\end{align}
where the second step uses \eqref{eqn:ss-7}, and the last step uses \eqref{eqn:10B}. 
It thus follows that $\falln$, for all $\kappa\in (0,\kappa_0]$ and $(\lambda,u)\in \dI_\kappa^{\sss(n),2}$, 
\[
\frac{\varphi_\lambda^{\sss(n)}(\kappa u)}{\varphi_\lambda^{\sss(n)}(u)} 
\leq 
\frac{\varphi_\lambda^{\sss(n)}\big((1-\kappa) u\big)}{\varphi_\lambda^{\sss(n)}(u)} 
\leq 
\frac{1}{1+C_{\ref{eqn:55}}\kappa}\, ,
\]
where the first inequality follows since $\varphi_\lambda^{\sss(n)}(\cdot)$ is an increasing function. This completes the proof of Lemma \ref{lem:varphi-prop}. \qed 

\medskip

Next, we study the sum of squares of weights in the component of vertex $1$, as well as the inclusion of maximal weight vertices within this component.

\begin{prop}\label{prop:lem-ss-4}
There exist $\delta_0>0$ and $\lambda_{\ref{it:lem-ss-4}} \geq \lambda_{\ref{eqn:ss-7a}}$ such that the following hold for all large $n$ and for $\lambda \in [\lambda_{\ref{it:lem-ss-4}},\, \eps_{\ref{eqn:ss-7a}} n^\eta]$:
\begin{gather}
\pr\bigg(\sum_{j\in \cC(1, \cG_n(\lambda)),\, j\neq 1 } \theta_{j,\lambda}^2 \leq (1+\delta_0)\lambda\bigg)\leq \exp(-C\lambda)\, ,\label{it:lem-ss-4} \ \ \ \
\text{ and }\\
\pr\bigg(j\notin \cC\big(1, \cG_n(\lambda)\big)\ \ \text{ for some }\ \
1\leq j\leq \frac{\lambda^{1/\eta}}{\log^3\lambda}\bigg) 
\leq 
\exp\big(-C\log^{3\alpha}\lambda\big)\, .\label{it:lem:ss-8} 
\end{gather}	
\end{prop}

\noindent{\bf Proof of \eqref{it:lem-ss-4}:}
Define $\psi_\lambda^{\sss(n)}:[0, \infty)\to [0, \infty)$ by
\[
\psi_\lambda^{\sss(n)}(u):=\sum_{j=2}^n \theta_{j,\lambda}^2\big(1-e^{-\theta_{j,\lambda} u}\big)\, .
\] 
Since $(1-e^{-s})\geq (s-1+e^{-s})/s$ for $s\in (0,\infty)$, we have $\psi_\lambda^{\sss(n)}(u) \geq \varphi_\lambda^{\sss(n)}(u)$ for $u\geq 0$. 
Further, $\falln$ and for $(\lambda, u) \in \dI^{\sss(n), 2}$, 
\[
u\big(\psi_\lambda^{\sss(n)}(u) - \varphi_\lambda^{\sss(n)}(u)\big) 
= 
\sum_{j=2}^n\theta_{j,\lambda}
\big(1-e^{-u\theta_{j,\lambda}} - u\theta_{j,\lambda}e^{-u\theta_{j,\lambda}}\big)
\asymp 
u^{\tau-2} 
\asymp 
u\varphi_\lambda^{\sss(n)}(u)\, ,
\]
where the penultimate step uses \eqref{eqn:10}, and the last step uses \eqref{eqn:ss-7}. 
Hence, there exists $\delta_0 >0$ such that $\falln$ and for $(\lambda, u)\in \dI^{\sss(n), 2}$, $\psi_\lambda^{\sss(n)}(u)\geq (1+3\delta_0) \varphi_\lambda^{\sss(n)}(u)$. 
Combined with \eqref{eqn:ss-7a}, we see that $\falln$ and $\lambda\in [\lambda_{\ref{eqn:ss-7a}},\, \eps_{\ref{eqn:ss-7a}} n^{\eta}]$, 
\begin{equation}\label{eqn:11}
\psi_\lambda^{\sss(n)}\big(s^{\sss(n)}(\lambda)\big) 
\geq 
(1+3\delta_0)\varphi_\lambda^{\sss(n)}\big(s^{\sss(n)}(\lambda)\big)\, .
\end{equation}

Let $\kappa_0$ be as in Lemma \ref{lem:varphi-prop}, and choose $\kappa_1 \in (0,\kappa_0]$ small so that 
\begin{equation}\label{eqn:12}
(1-2\kappa_1)(1+3\delta_0)\geq (1+2\delta_0)\, .
\end{equation}
Recall the process $Z_{\lambda}^{\sss n, (1)}$ from \eqref{eqn:8888}.
By Proposition \ref{prop:lem-ss-3} and Lemma \ref{lem:domin-ald-limic-1}, writing $n_{\ref{eqn:12A}}=n_{\ref{eqn:66}}(\kappa_1)$ and $\lambda_{\ref{eqn:12A}}= \lambda_{\ref{eqn:66}}(\kappa_1)$, 
we have, for all 
$n\geq n_{\ref{eqn:12A}}$ and $\lambda\in [\lambda_{\ref{eqn:12A}},\, \eps_{\ref{eqn:ss-7a}} n^\eta]$, 
\begin{equation}\label{eqn:12A}
\pr\big(Z_{\lambda}^{\sss n, (1)}(\cdot)\ \text{ hits zero before } (1-2\kappa_1)s^{\sss(n)}(\lambda) \big) 
\leq 
\exp\big(-C\lambda^{1/(\tau-3)}\big) \, .
\end{equation}
Consider the coupling between $\cG_n(\lambda)$ and the random variables $\xi_j^n$, $2\leq j\leq n$, as mentioned below \eqref{eqn:8888}.
Then note that in this coupling, on the complement of the event in \eqref{eqn:12A}, 
$j\in\cC\big(1, \cG_n(\lambda)\big)$ whenever $\xi_j^n \leq (1-2\kappa_1)s^{\sss(n)}(\lambda)$, and consequently,
\begin{align}
\sum_{ j\in \cC(1, \cG_n(\lambda)),\, j\neq 1} \theta_{j,\lambda}^2 
&\geq 
\sum_{j=2}^n \theta_{j,\lambda}^2\cdot \ind\big\{\xi_j^n \leq (1-2\kappa_1)s^{\sss(n)}(\lambda)\big\}
\notag \\
&=  
\sum_{j=2}^n \theta_{j,\lambda}^2 
\bigg(\ind\big\{\xi_j^n \leq (1-2\kappa_1)s^{\sss(n)}(\lambda)\big\} - 
\pr\big(\xi_j^n \leq (1-2\kappa_1)s^{\sss(n)}(\lambda)\big)\bigg) \notag \\
& \qquad 
+ \psi_\lambda^{\sss(n)}\big((1-2\kappa_1)s^{\sss(n)}(\lambda)\big):= 
\fT_1^{\sss(n)} + \fT_2^{\sss(n)}. \notag
\end{align}

To lower bound $\fT_2^{\sss(n)}$, note that 
$1-e^{-su} \geq s(1-e^{-u})$ for all $s\in (0,1)$ and $u\geq 0$. 
Thus, 
$\psi_\lambda^{\sss(n)}\big((1-2\kappa_1)s^{\sss(n)}(\lambda)\big) 
\geq 
(1-2\kappa_1)\psi_\lambda^{\sss(n)}\big(s^{\sss(n)}(\lambda)\big)$.  
Thus, $\falln$ and $\lambda\in [\lambda_{\ref{eqn:ss-7a}},\, \eps_{\ref{eqn:ss-7a}} n^{\eta}]$,
\begin{equation}\label{eqn:13}
\psi_\lambda^{\sss(n)}\big((1-2\kappa_1)s^{\sss(n)}(\lambda)\big) 
\geq 
(1-2\kappa_1)(1+3\delta_0) \varphi_\lambda^{\sss(n)}\big(s^{\sss(n)}(\lambda)\big) 
\geq 
(1+2\delta_0) \lambda\, ,
\end{equation} 
where the first step uses \eqref{eqn:11}, and the second step uses \eqref{eqn:12} and the definition of $s^{\sss (n)}(\lambda)$.

Turning to $\fT_1^{\sss(n)}$, an application of the bounded difference inequality shows that $\falln$ and $\lambda \in [1,\, \eps_{\ref{eqn:ss-7a}} n^\eta]$, 
\begin{equation}\label{eqn:14}
\pr\big(\fT_1^{\sss(n)} \leq -\delta_0 \lambda\big) \leq \exp\big(-C\lambda^2\big)\, .
\end{equation}
	
Combining \eqref{eqn:12A}, \eqref{eqn:13}, and \eqref{eqn:14}, it follows that there exists $\lambda_{\ref{eqn:22}} \geq \lambda_{\ref{eqn:12A}}$ such that $\falln$ and $\lambda\in [\lambda_{\ref{eqn:22}},\, \eps_{\ref{eqn:ss-7a}} n^{\eta}]$, 
\begin{align}\label{eqn:22}
\pr\bigg(\sum_{ j\in \cC(1, \cG_n(\lambda)),\, j\neq 1} \theta_{j,\lambda}^2 
\leq 
(1+\delta_0)\lambda\bigg) 
\leq 
\exp\big(-C\lambda^{2\wedge \frac{1}{\tau-3}}\big) 
\leq 
\exp\big(-C \lambda\big)\, .
\end{align}
This completes the proof of \eqref{it:lem-ss-4}. 
\qed

\medskip

\noindent{\bf Proof of \eqref{it:lem:ss-8}:}
Write $\sum_\ast$ for $\sum_{j=2}^{\frac{\lambda^{1/\eta}}{\log^3\lambda}} $, and note that $\falln$ and $\lambda \in [\lambda_{\ref{eqn:ss-7a}},\,  \eps_{\ref{eqn:ss-7a}} n^\eta]$,
\begin{align}\label{eqn:32}
&\sum\displaystyle_\ast
\pr\big(\xi^n_j \geq  (1-2\kappa_1) s^{\sss(n)}(\lambda)\big) 
= 
\sum\displaystyle_\ast \exp\big(-\theta_{j,\lambda}(1-2\kappa_1)s^{\sss(n)}(\lambda)\big) \notag\\
&\hskip20pt
\leq 
\sum\displaystyle_\ast 
\exp\bigg(-\frac{C}{j^\alpha}\lambda^{\frac{1}{\tau-3}}\bigg) 
\leq 
\frac{\lambda^{1/\eta}}{\log^3\lambda} \exp\big(-C\log^{3\alpha}\lambda\big) 
\leq 
C' \exp\big(-C'' \log^{3\alpha}\lambda\big)\, , 
\end{align}
where the second step uses \eqref{eqn:ss-8}, and the last step uses the fact $3\alpha > 1$. 

Now consider again the coupling between $\cG_n(\lambda)$ and the random variables $\xi_j^n$, $2\leq j\leq n$, as mentioned below \eqref{eqn:8888}.
As already noted below \eqref{eqn:12A}, in this coupling, on the complement of the event in \eqref{eqn:12A}, 
$j\in\cC\big(1, \cG_n(\lambda)\big)$ whenever $\xi_j^n \leq (1-2\kappa_1)s^{\sss(n)}(\lambda)$.
Thus, combining \eqref{eqn:32} and \eqref{eqn:12A} completes the proof. 
\qed

\subsection{Height bounds for a branching process}\label{sec:54}

Recall the definition of $\nu_n$ from \eqref{eqn:sigm2-nu-def}, and define 
$v_i = w_i/\nu_n$ for $2\leq i\leq n$. 
Let $V_n$ be a random variable with distribution
\begin{align}\label{def:V-n}
\frac{1}{\sum_{i=2}^n v_i}\sum_{i=2}^n v_i\cdot\delta\{v_i\}
=\frac{1}{\ell_n}\sum_{i=2}^n w_i\cdot\delta\{v_i\}\, .
\end{align}
Let $\Poi(V_n)$ denote a random variable that conditional on $V_n$ is distributed as a Poisson random variable with mean $V_n$. 
We will need the following property of the random variable $\Poi(V_n)$.
\begin{lemma}\label{lem:5}
There exist $ C'_{\ref{eqn:2323}}> C_{\ref{eqn:2323}} >0$ such that for all $n\geq 1$,
\begin{align}\label{eqn:2323}
C'_{\ref{eqn:2323}}/u^{\tau-2}\geq \pr\big(\Poi(V_n) \geq u\big) \geq C_{\ref{eqn:2323}}/u^{\tau-2} \qquad \text{ for } \qquad  1\leq u\leq \crv_2/2.
\end{align}
\end{lemma}
The proof of Lemma \ref{lem:5} is given in Section \ref{sec:proof-lem-5}. 
The main result of this section, stated in the next proposition, describes height asymptotics of a branching process $T_n$ with offspring distribution $\Poi(V_n)$. 
Note that $\E\big[\Poi(V_n)\big] = 1$, so that $T_n$ is finite almost surely. 

\begin{prop}\label{prop:lem-6-bp-height}
There exists $\gamma_* > 0$ and $C>0$ such that for $\gamma \in (0,\gamma_*]$ and all $n\geq 1$,
\[
\pr\big(\height(T_n) \geq \gamma n^{\eta}\big) \leq \frac{C}{n^{\alpha} \gamma^{1/(\tau-3)}}\, . 
\]
\end{prop}  

The proof of this result uses a technique recently developed in \cite{addario2019most}.
For a random variable $X$, define L\'evy's concentration function as
\[
Q(X, u):= \sup_{x\in \bR}\, \pr\big(x\leq X\leq x+u\big)\, ,\ \ \ u>0 . 
\]
Let $\big(X_i\, ;\, i\geq 1\big)$ be an i.i.d. sequence with 
$X_1 \equald \Poi(V_n)-1$.
Let $\big(S_k\, ;\, k\geq 0\big)$ be a random walk with $S_k-S_{k-1}=X_k$, $k\geq 1$, 
and for $x\in\bZ$, write $P_x$ for the probability distribution of this random walk when $S_0 = x$.  
By \cite[Theorem 3.1]{esseen1968concentration} there exists a (universal) constant $C>0$ such that for any $x\in\bZ$, under $P_x$,
\begin{equation}\label{eqn:21}
Q(S_k, u) \leq 
\frac{Cu}{\sqrt{k\cdot\E\big[\big(X_1 - X_2\big)^2\ind\big\{\big|X_1 - X_2\big| \leq u\big\}\big]}}\, ,
\ \ \ k\geq 1\, .
\end{equation}
We start with the following lemma. 
\begin{lemma}\label{lem:eqn-22}
Let $Y = X_1 - X_2$. Let $\beta_* >0$ be such that $\falln$, one has $\beta_*n^\alpha \leq v_2/2 -1$. Then $\falln$ and for any $u\in\big[1,\,  \beta_* n^\alpha\big]$,
\[\E\big[Y^2\ind\big\{|Y|\leq u\big\}\big] \geq C u^{4-\tau}\, .\]
\end{lemma}
\noindent{\bf Proof:}
Choose $\Delta\in (0, 1)$ such that $C_{\ref{eqn:2323}}=2C'_{\ref{eqn:2323}}\Delta^{\tau-2}$.
Then an application of Lemma \ref{lem:5} shows that $\falln$, for any $u\in [2/\Delta,\, \beta_*n^{\alpha}]$ and $y\in [1, u\Delta]$,
\begin{align}\label{eqn:232}
\pr\big(\Poi(V_n)\in [y,\, u]\big)
&=
\pr\big(\Poi(V_n)\geq y\big)-\pr\big(\Poi(V_n)> u\big)
\geq
\frac{C_{\ref{eqn:2323}}}{y^{\tau-2}}-\frac{C'_{\ref{eqn:2323}}}{u^{\tau-2}}\notag\\
&\geq
\frac{C_{\ref{eqn:2323}}}{(u\Delta)^{\tau-2}}-\frac{C'_{\ref{eqn:2323}}}{u^{\tau-2}}
=\frac{C'_{\ref{eqn:2323}}}{u^{\tau-2}}\, .
\end{align}
Hence, $\falln$ and for any $u\in [2/\Delta,\, \beta_*n^{\alpha}]$,
\begin{align}
\E\big[Y^2\ind\big\{|Y|\leq u\big\}\big] 
&=
2\int_0^u y \pr\big(u\geq |Y|\geq y\big) dy 
\geq 
2\pr\big(X_2 =-1\big)\int_1^{u\Delta} y \pr\big(u\geq X_1+1\geq y\big) dy \notag \\
&\geq 
2\pr\big(X_2 =-1\big) \int_1^{u\Delta}  \frac{C'_{\ref{eqn:2323}}y}{u^{\tau-2}} dy 
\geq 
C\pr\big(X_2 =-1\big) u^{4-\tau}, \label{eqn:1225}
\end{align}
where the penultimate step uses \eqref{eqn:232}.
Now, using Assumption \ref{ass:wts}, 
\begin{equation}\label{eqn:1226}
\pr\big(X_2 = -1\big) 
= 
\E\big[e^{-V_n}\big] 
= 
\frac{1}{\ell_n}\sum_{i=2}^n w_i e^{-w_i/\nu_n} 
\geq 
\frac{C}{n}\sum_{i=n/4}^{n/2} \left(\frac{n}{i}\right)^\alpha e^{-C'(\frac{n}{i})^\alpha} 
\geq C'' >0\, .
\end{equation}
This yields the desired result for $u\in [2/\Delta,\, \beta_*n^{\alpha}]$.
Now, an argument similar to the one used in \eqref{eqn:1226} will show that $\pr(X_1=0)\geq C>0$, which would in turn imply that for $u\in [1, 2/\Delta]$, 
\[
\E\big[Y^2\ind\big\{|Y|\leq u\big\}\big]
\geq 
\pr(X_1=0)\pr(X_2=-1)
\geq 
C'
\geq 
C' u^{4-\tau}(\Delta/2)^{4-\tau}\, .
\]
Using this last observation we extend the lower bound to the interval $u\in [1,\, \beta_*n^{\alpha}]$.
\qed	

\medskip

\noindent{\bf Proof of Proposition \ref{prop:lem-6-bp-height}:}
Consider the intervals $I_l = [2^{l-1}, 2^{l+2}) $ for $l\geq 1$. 
Define
\[
r_{nl}:= 
\min\bigg\{k\geq 1\, :\, \sup_{x\in I_l}\, P_x\bigg(\bigcap_{j=1}^k\set{S_j\in I_l} \bigg) \leq 1/2\bigg\}\, .
\]
By \eqref{eqn:21} and Lemma \ref{lem:eqn-22}, when $2^{l+2}\leq \beta_*n^\alpha$, we have $\falln$ and a large constant $C_{\ref{eqn:233}}>0$,
\begin{align}\label{eqn:233}
Q\big(S_{C_{\ref{eqn:233}} 2^{l(\tau-2)}}\, ,\, 2^{l+2}\big)\leq \frac{C 2^{l+2}}{\sqrt{C_{\ref{eqn:233}} 2^{l(\tau-2)} (2^{l+2})^{(4-\tau)}} } \leq 1/2,
\end{align}
Hence, $\falln$,
\begin{equation}\label{eqn:23}
r_{nl}\leq C_{\ref{eqn:233}} 2^{l(\tau-2)}, \qquad \text{ for } 1\leq l\leq m(\beta_*),
\end{equation}	
where
$m(\beta) :=\max \set{k\in \bZ_{>0}: 2^k \leq \lfloor \beta n^{\alpha} \rfloor}$ for $\beta>0$.

For the rest of this proof, we will work with the random walk $\big(S_k\, ;\, k\geq 0\big)$ started at $S_0=1$, i.e., we will work under the measure $P_1$. 
Let $\xi=\inf\,\big\{k\geq 0\, :\, S_k=0 \big\}$.
By \cite[Proposition 1.7]{addario2019most},
\begin{align}\label{eqn:61}
\height(T_n)\stod\,\, 3\sum_{k=0}^{\xi-1} \frac{1}{S_k}=: 3\cdot J_n(\xi)\, .
\end{align}
Following
\cite{addario2019most} we derive tail bounds for $J_n(\xi)$ by decomposing the trajectory of the random walk into various ``scales'' which we now define. 
Let $\zeta_0 = 0$ and $R_0 =1$. 
For $i\geq 0$, define the stopping times
$\zeta_{i+1} = \min\big\{t\geq \zeta_i\, :\, S_t \notin\big[2^{R_i-1}, 2^{R_i+2}\big)\big\}$, 
and let 
$R_{i+1} = \max\big\{l\, :\, S_{\zeta_{i}+1}\geq 2^l\big\}$. 
	
Next, for $0\leq k\leq \xi$, let $\Lambda(k)$ denote the scale of $\mvS$ at time $k$. 
Precisely, let $j=\max\big\{i\, :\, \zeta_i\leq k\big\}$ be the most recent epoch for a change in scale, and let $\Lambda(k) = R_j$. 
Finally, for $l\geq 1$, define
\[
J_{nl}(\xi):= \sum_{k=0}^{\xi -1}\frac{1}{S_k}\ind\set{\Lambda(k) = l}\, .
\]
Now consider $\beta\in [4n^{-\alpha}, \beta_{*}]$, and note that $\falln$ and for any such choice of $\beta$, $m(\beta) \geq 2$, where $m(\cdot)$ is as defined below \eqref{eqn:23}. 
Define
\[
b_l:= 18\bigg(m(\beta) -l+1 +2\log_2\big(m(\beta) -l+1\big)\bigg)\, , \qquad 1\leq l\leq m(\beta)\, .
\]
Then $\falln$,
\begin{align}\label{eqn:24}
\sum_{l=1}^{m(\beta)} \frac{b_l r_{nl}}{2^{l-1}} 
&\leq 
18 C_{\ref{eqn:233}}\sum_{l=1}^{m(\beta)}
\frac{1}{2^{l-1}}
\bigg(m(\beta) -l+1 +2\log_2\big(m(\beta) -l+1\big)\bigg)\cdot 2^{l(\tau-2)} \notag\\
&\leq
36 C_{\ref{eqn:233}} 2^{(m(\beta)+1)(\tau-3)}\sum_{j=1}^{\infty} \frac{(j+2\log_2 j)}{2^{j(\tau-3)}}\leq C_{\ref{eqn:24}} (\beta n^\alpha)^{\tau -3}\, ,
\end{align}
where the first step uses \eqref{eqn:23}, and in the second step we have made the substitution $j=m(\beta) -l+1 $.

Note that $\max_{0\leq k\leq \xi-1}\Lambda(k)\leq m(\beta)$ on the event 
$\big\{ \max_{0\leq k\leq \xi} S_k< \lfloor \beta n^\alpha \rfloor\big\}$.
Hence, $\falln$ and for $\beta \in [4n^{-\alpha},\beta_* ]$, 
\begin{align}
&
P_1\bigg(J_n(\xi) \geq C_{\ref{eqn:24}}(\beta n^\alpha)^{\tau-3}\bigg) 
\leq 
P_1\bigg(\max_{0\leq k\leq \xi} S_k\geq \lfloor \beta n^\alpha \rfloor\bigg) 
+ 
P_1\bigg( \sum_{l=1}^{m(\beta)} J_{nl}(\xi) \geq \sum_{l=1}^{m(\beta)} \frac{b_l r_{nl}}{2^{l-1}} \bigg) \notag\\
&\hskip60pt
\leq 
\frac{1}{\lfloor \beta n^\alpha \rfloor} 
+ 
\sum_{l=1}^{m(\beta)} P_1\bigg(J_{nl}(\xi) \geq \frac{b_l r_{nl}}{2^{l-1}} \bigg) 
\leq 
\frac{1}{\lfloor \beta n^\alpha \rfloor} 
+ 
\sum_{l=1}^{m(\beta)} \frac{1}{2^l}\cdot \frac{2}{2^{b_l/18}}\, ,\notag
\end{align} 
where the first step uses \eqref{eqn:24},
the second step follows from a simple application of the optional stopping theorem, 
and the third step follows from \cite[Theorem 3.6]{addario2019most}. 
Using the expression for $b_l$, the above bound yields, 
$\falln$ and for $\beta \in [4n^{-\alpha},\beta_* ]$,
\begin{equation}\label{eqn:405}
P_1\bigg(J_n(\xi) \geq C_{\ref{eqn:24}}(\beta n^\alpha)^{\tau-3}\bigg) 
\leq 
\frac{1}{\lfloor \beta n^\alpha \rfloor} 
+ 
\sum_{l=1}^{m(\beta)} \frac{2}{2^{m(\beta)+1} \big(m(\beta)-l+1\big)^2} 
\leq 
\frac{C}{\beta n^\alpha}\, .
\end{equation}
To reparametrize from $\beta$ to $\gamma$ as  in the statement of Proposition \ref{prop:lem-6-bp-height}, write 
$\gamma = 3 C_{\ref{eqn:24}}\beta^{\tau-3}$, and $\gamma_* = 3 C_{\ref{eqn:24}}\beta_*^{\tau-3}$.
Then \eqref{eqn:61} and \eqref{eqn:405} give the desired result for all large $n$ and for $\gamma \in [3 C_{\ref{eqn:24}}4^{\tau-3}n^{-\eta},\, \gamma_*]$. 
Now choose a larger constant to make the bound work for 
$\gamma \in (0,\, 3C_{\ref{eqn:24}}4^{\tau-3}n^{-\eta}]$ and for all $n\geq 1$. 
This completes the proof.
\qed

\subsection{Diameter outside the component of the vertex $1$}\label{sec:diam}
For $\lambda\geq 0$ and $\delta>0$, let $H_n(\lambda, \delta)$ be the random graph constructed in the following way:
Let the vertex set be 
\begin{align}\label{eqn:466}
\cV'_{\lambda}:=[n]\setminus V\big(\cC_1^n(\lambda)\big)\, ,
\end{align}
and place edges independently between $i, j\in\cV'_{\lambda}$ with probability
$1-\exp\big(-p^n_{(1+\delta)\lambda}w_i w_j/\ell_n\big)$.

\begin{prop}\label{prop:lem-ss-7} 
There exist $\delta_1>0$, $\eps_{\ref{eqn:94}}\in (0, \eps_{\ref{eqn:ss-7a}} ]$ and
$\lambda_{\ref{eqn:94}}\geq\lambda_{\ref{it:lem-ss-4}}$ such that for all large $n$, for $\lambda\in [\lambda_{\ref{eqn:94}}, \eps_{\ref{eqn:94}} n^\eta]$, and for every $\Delta\in (0,1/2]$,
\begin{align}\label{eqn:94}
\pr\bigg(\diam\big(H_n(\lambda, \delta_1)\big) \geq \frac{n^{\eta}}{\lambda^{1-\Delta}}\bigg) 
\leq 
C \lambda^{\frac{\tau+1}{\tau-3}}\exp\big(-C' \lambda^\Delta\big)\, .
\end{align}
\end{prop}

Let us specify here how we choose the thresholds.
Choose $\lambda_{\ref{eqn:94}}\geq\lambda_{\ref{it:lem-ss-4}}$, $\eps_{\ref{eqn:94}}\in (0, \eps_{\ref{eqn:ss-7a}} ]$,  $\varrho\in (0,1/2)$, and $\delta_1 >0$ such that 
\begin{equation}\label{eqn:944}
\varrho\cdot\big(\lambda_{\ref{eqn:94}}\big)^{-1/2}\leq\gamma_*\, ,\ \ \
\frac{(1+\delta_0)(1-2\varrho)}{(1+\eps_{\ref{eqn:94}})^2} - (1+\delta_1)\geq \frac{\delta_0}{2}\, , \ \ \
\text{ and } \ \ \
\eps_{\ref{eqn:94}}\leq \frac{\varrho}{10}\, ,
\end{equation}
where $\gamma_*$ is as in Proposition \ref{prop:lem-6-bp-height}, and
$\delta_0$ is as in Proposition \ref{prop:lem-ss-4}. 

Let $i_1 < i_2 < \cdots$ be the vertices in $\cV'_{\lambda}$. 
For $1\leq k\leq |\cV'_{\lambda}|$,  write $\cC_{\res}(i_k; \lambda)$ for the component of $i_k$ in 
$H_n(\lambda, \delta_1) \setminus\big\{i_1, \ldots, i_{k-1}\big\}$.
Write 
\begin{align}\label{eqn:45}
\pr_1(\cdot)=\pr\big(\cdot\, \big|\, \cC_1^n(\lambda)\big)\, ,\ \ \text{ and }\ \ 
\bE_1[\cdot]=\bE\big[\cdot\, \big|\, \cC_1^n(\lambda)\big]
\, .
\end{align}
Now, 
\begin{equation}\label{eqn:25}
\pr_1\bigg(\diam\big(H_n(\lambda, \delta_1)\big) \geq \frac{n^\eta}{\lambda^{1-\Delta}}\bigg) 
\leq 
\sum_k \pr_1\bigg(\diam\big(\cC_{\res}(i_k; \lambda)\big) \geq \frac{n^\eta}{\lambda^{1-\Delta}} \bigg)\, . 
\end{equation} 
For $1\leq k\leq |\cV'_{\lambda}|$, let
\begin{equation}\label{eqn:1046}
\cV'_{\lambda, k} := \cV'_{\lambda} \setminus \big\{i_1, i_2, \ldots, i_{k-1}\big\}. 
\end{equation}
Note that for any $s>0$, 
\begin{align}\label{eqn:88}
\Bern\big(1-e^{-s}\big)\stod\Poi(s)\, ,
\end{align} 
and consequently,
the breadth-first exploration tree of $\cC_{\res}(i_k; \lambda)$ starting from $i_k$ is upper bounded by a multitype branching process with state space $\cV'_{\lambda, k}$ in which 
the type of the root is $i_k$, and any vertex of type $i$ has 
$\Poi\big(p^n_{(1+\delta_1)\lambda}w_i w_j/\ell_n\big)$ many type $j$ children for $i, j\in\cV'_{\lambda, k}$. 
This leads us to the following definition which will be useful in the proof.
\begin{defn}\label{defn:mtbp}
For $\lambda\geq 0$, $D\subseteq [n]$, and $i\in D$, let ~$\mtbp_n^{\lambda, i}(D)$ be a multitype branching process tree with type space $D$ that is rooted at a vertex of type $i$, and in which a vertex of type $j$ has $\Poi\big(p^n_{(1+\delta_1)\lambda}w_j w_k/\ell_n\big)$ many children of type $k$ for each $j,k\in D$. 
Let $\mtbp_n^{i}(D)=\mtbp_n^{0, i}(D)$, and note that in this case a vertex of type $j$ has  $\Poi(w_j w_k/\nu_n \ell_n)$ many children of type $k$.  
\end{defn}
	 
The bounds in the next lemma will be crucial for dealing with $\diam\big(\cC_{\res}(i_k; \lambda)\big)$.

\begin{lemma}\label{lem:26}
For any $\lambda\geq 0$, $1\leq k\leq |\cV'_{\lambda}|$, and $h>0$,
\[
\pr_1\big(\diam\big(\cC_{\res}(i_k; \lambda)\big) \geq h \big) 
\leq 
\pr_1\big(\height\big(\mtbp_n^{\lambda, i_k}\big(\cV'_{\lambda, k}\big)\big)\geq h/2\big)\, .
\]
Further, for any $i\in D\subseteq [n]$, and $h\in\bZ_{>0}$,
\[
\pr\bigg(\height\big(\mtbp_n^{\lambda, i}(D) \big)\geq h \bigg)
\leq 
\exp\big((1+\delta_1)h\lambda n^{-\eta}\big) \cdot 
\pr\bigg(\height\big(\mtbp_n^i(D) \big) \geq h \bigg)\, .
\]
\end{lemma}

\noindent{\bf Proof:}
The first assertion follows from the discussion above Definition \ref{defn:mtbp}. 
For the second assertion, write $\zeta=1+(1+\delta_1)\lambda n^{-\eta}$, and note that $\mtbp_n^{i}(D)$ can be obtained as a subtree of $\mtbp_n^{\lambda, i}(D)$ by killing every child independently with probability $1-\zeta^{-1}$. 
Write $\cA$ for the event in which $\height\big(\mtbp_n^{\lambda, i}(D)\big)\geq h$, and no vertex in the leftmost path of length $h$ starting from the root in $\mtbp_n^{\lambda, i}(D)$ is killed. 
Then
\begin{align}
\pr(\cA) = \zeta^{-h} \pr\bigg(\height\big(\mtbp_n^{\lambda, i}(D)\big)\geq h\bigg)
\geq
\exp\big(-(1+\delta_1)h\lambda n^{-\eta}\big) 
\pr\bigg(\height\big(\mtbp_n^{\lambda, i}(D)\big)\geq h\bigg) .\notag
\end{align}
To finish the proof, note that $\cA$ implies $\height\big(\mtbp_n^i(D)\big)\geq h$.
\qed

\medskip

We next record two properties of the above branching process. 
Let $V_n$ be as defined around \eqref{def:V-n}.
	
\begin{lemma}\label{lem:mtbp-size-bias}
Fix $i \in [n]\setminus \{1\}$. 
Consider $\mtbp_n^i\big([n]\setminus \{1\}\big)$, and erase the types of all vertices.
Then this tree has the same distribution as a branching process tree where the root has $\Poi\big(w_i/\nu_n\big)$ many children, and every other vertex has $\Poi(V_n)$ many children.
\end{lemma}

This result was noted in the discussion above \cite[Proposition 3.2]{NorRei06}.
We will briefly include the proof.

\medskip

\noindent{\bf Proof of Lemma \ref{lem:mtbp-size-bias}:}
Note that $\mtbp_n^i\big([n]\setminus \{1\}\big)$ can be constructed by starting from the root and inductively continuing through the generations as follows:
to each vertex of type $j$ assign $\Poi(w_j/\nu_n)$ many children, and conditional on this step, declare the type of every child independently to be $k\in\{2,\ldots, n\}$ with probability $w_k/\ell_n$.
Thus, in the tree obtained by erasing the types of every vertex in $\mtbp_n^i\big([n]\setminus \{1\}\big)$, the root has $\Poi(w_i/\nu_n)$ many children, and every other vertex has $\Poi(w_Y/\nu_n)$ many children, where $\pr\big(Y=j\big)=w_j/\ell_n$, $j=2,\ldots, n$.
The proof is complete upon noting that $w_Y/\nu_n\equald V_n$.
\qed

\medskip

The next lemma follows easily from Definition \ref{defn:mtbp}. 

\begin{lemma}\label{lem:stod-mtbp}
Consider $i\in D\subseteq D' \subseteq [n]\setminus \{1\}$. 
Then $\mtbp_n^i(D)$ can be coupled with $\mtbp_n^i(D')$ so that the former is a subtree of the latter.
\end{lemma}
	
Given $\cV'_{\lambda}$, we will now construct a three layer branching process ($\tlbp$), which we will use to obtain tail bounds on $\height\big(\mtbp_n^{i_k}\big(\cV'_{\lambda, k}\big)\big)$.
Let $\varrho$ be as in \eqref{eqn:944}, and $\Delta\in(0, 1/2]$ be as in the statement of Proposition \ref{prop:lem-ss-7}.

\begin{defn}\label{def:three-layer-cake}
For $\lambda>0$ and $1\leq k\leq |\cV'_{\lambda}|$, consider the following (potentially) three layer process $\tlbp_n^{i_k}(\lambda)$: 
\begin{enumeratea}
\item {\bf Layer 1:} 
Start $\mtbp_n^{i_k}(\cV'_{\lambda})$, and run this process up to generation $(1-2\varrho)n^{\eta}/(2\lambda^{1-\Delta})$. 
Call this the first layer. 
If there is at least one vertex in generation $(1-2\varrho)n^{\eta}/(2\lambda^{1-\Delta})$, then we say that the first layer has been fully activated. 
\item {\bf Layer 2:} 
If the first layer is fully activated, then starting from every vertex $v$ in generation  $(1-2\varrho)n^{\eta}/(2\lambda^{1-\Delta})$, run independent $\mtbp_n^{\mathrm{type}(v)}\big([n]\setminus [k]\big)$ processes up to generation
$\varrho n^{\eta}/(2\lambda^{1-\Delta})$, 
where $\mathrm{type}(v)\in[n]$ denotes the type of the vertex $v$. 
Call this the second layer. 
If any of these branching processes survives up to generation $\varrho n^{\eta}/(2\lambda^{1-\Delta})$, then we say that the second layer has been fully activated; in this case, there is at least one vertex in generation $(1-\varrho)n^{\eta}/(2\lambda^{1-\Delta})$ of $~\tlbp_n^{i_k}(\lambda)$.
\item {\bf Layer 3:} 
If the second layer is fully activated, then starting from every vertex $v$ in generation  $(1-\varrho)n^{\eta}/(2\lambda^{1-\Delta})$, run independent $\mtbp_n^{\mathrm{type}(v)}\big([n]\setminus \{1\}\big)$ processes.
Call this the third layer. 
If any of these branching processes survives up to generation $\varrho n^{\eta}/(2\lambda^{1-\Delta})$, then we say that the third layer has been fully activated.
\end{enumeratea}
\end{defn}

Write $\fF_{n,\lambda}^{i_k}$ for the event that all three layers have been fully activated, and note that 
$
\fF_{n,\lambda}^{i_k}
=
\big\{\height\big(\tlbp_n^{i_k}(\lambda)\big) \geq n^{\eta}/(2\lambda^{1-\Delta})\big\}
$. 
Since $\cV'_{\lambda, k}\subseteq\cV'_{\lambda}$ and 
$\cV'_{\lambda, k}\subseteq [n]\setminus [k]\subseteq [n]\setminus \{1\}$, 
using Lemma \ref{lem:stod-mtbp}, 
$\height\big(\mtbp_n^{i_k}\big(\cV'_{\lambda, k}\big)\big)
\stod
\height\big(\tlbp_n^{i_k}\big(\lambda\big)\big)
$.
In particular, 
\begin{equation}\label{eqn:27}
\pr_1\bigg(\height\big(\mtbp_n^{i_k}\big(\cV'_{\lambda, k}\big)\big) \geq \frac{n^\eta}{2\lambda^{1-\Delta}}\bigg) 
\leq 
\pr_1\big(\fF_{n,\lambda}^{i_k}\big)\, . 
\end{equation}

For $\lambda>0$, define
$
\fE_{n,\lambda} := 
\big\{\sum_{j=2}^n \theta_{j,\lambda}^2\ind\big\{j\in \cC^n_1(\lambda)\big\} \geq (1+\delta_0)\lambda \big\}
$.
Note that by \eqref{it:lem-ss-4} and \eqref{eqn:712}, $\falln$ and for
$\lambda\in[\lambda_{\ref{it:lem-ss-4}},\, \eps_{\ref{eqn:ss-7a}} n^{\eta}]$,
\begin{align}\label{eqn:98}
\pr\big(\fE_{n,\lambda}^c\big)\leq\exp\big(-C\lambda\big)\, .
\end{align}
\begin{lemma}\label{lem:fnl-bound}
We have, for all $\Delta\in (0,1/2]$, for all large $n$, and for $\lambda\in [\lambda_{\ref{eqn:94}},\, \eps_{\ref{eqn:94}} n^\eta]$, on the event $\fE_{n,\lambda}$,
\[
\pr_1\big(\fF_{n,\lambda}^{i_k}\big) 
\leq 
C \lambda^{\frac{\tau+1}{\tau-3}} k^{-\frac{\tau+1}{\tau-1}}
\exp\bigg(-\, \frac{(1+\delta_0)(1-2\varrho)\lambda^\Delta}{2(1+ \eps_{\ref{eqn:94}})^2}\bigg)
\, ,   
\]
for $1\leq k\leq |\cV'_{\lambda}|$.
\end{lemma}

\noindent{\bf Proof:}
Fix $k\in\{1, \ldots, |\cV'_{\lambda}|\}$ and $\Delta\in(0, 1/2]$.
Write $\Gamma_{2\leadsto 3}$ for the number of vertices in generation $(1-\varrho)n^{\eta}/(2\lambda^{1-\Delta})$ of $\tlbp_n^{i_k}(\lambda)$ (i.e., vertices at the end of the second layer), whose subtree in the third layer has height at least $\varrho n^{\eta}/(2\lambda^{1-\Delta})$.  
Clearly,
\begin{equation}\label{eqn:515}
\pr_1\big(\fF_{n,\lambda}^{i_k}\big) 
\leq 			
\E_1\big[\Gamma_{2\leadsto 3}\big]\, .
\end{equation}
Hence, it is enough to find an upper bound for $\E_1\big[\Gamma_{2\leadsto 3}\big]$.

For $j\geq 1$, let  $\text{Gen}(j)$ be the set of vertices in the $j$-th generation of $\tlbp^{i_k}(\lambda)$, and let 
$
Y(j)
=
\sum_{v\in\text{Gen}(j)}w_{\text{type}(v)}
$ 
denote the sum of weights of vertices in the $j$-th generation.
Conditioning on the first two layers,  using Lemma \ref{lem:mtbp-size-bias} for the branching processes in the third layer, and then using Proposition \ref{prop:lem-6-bp-height} for the height of such branching processes, we get, $\falln$ and for all 
$\lambda\in [\lambda_{\ref{eqn:94}},\, \eps_{\ref{eqn:94}} n^\eta]$, 
\begin{equation}\label{eqn:gamma2-3}
\E_1\big[\Gamma_{2\leadsto 3}\big]
\leq 
\frac{C}{n^\alpha \big(\varrho/\lambda^{1-\Delta}\big)^{\frac{1}{\tau-3}}}\cdot \E_1\bigg[Y\bigg(\frac{(1-\varrho)n^{\eta}}{2\lambda^{1-\Delta}}\bigg)\bigg]\, .
\end{equation}
It can be checked by a direct computation that in the second layer, 
\begin{equation}\label{eqn:554}
\E_1\bigg[Y\bigg(\frac{(1-\varrho)n^{\eta}}{2\lambda^{1-\Delta}}\bigg)\bigg]
=
R_k
\E_1\bigg[Y\bigg(\frac{(1-\varrho)n^{\eta}}{2\lambda^{1-\Delta}}-1\bigg)\bigg]
=
\ldots
=
\big(R_k\big)^{\frac{\varrho n^{\eta}}{2\lambda^{1-\Delta}}}
\E_1\bigg[Y\bigg(\frac{(1-2\varrho)n^{\eta}}{2\lambda^{1-\Delta}}\bigg)\bigg]\, ,
\end{equation}
where $R_k=\sum_{j=k+1}^n w_j^2/(\nu_n \ell_n)$.
Similarly, in the first layer, 
\begin{equation*}
\E_1\bigg[Y\bigg(\frac{(1-2\varrho)n^{\eta}}{2\lambda^{1-\Delta}}\bigg)\bigg]
=
\bigg(\sum_{j\in\cV'_{\lambda}} \frac{w_j^2}{\nu_n \ell_n}\bigg) 
\E_1\bigg[Y\bigg(\frac{(1-2\varrho)n^{\eta}}{2\lambda^{1-\Delta}}-1\bigg)\bigg]
=\ldots=
\bigg(\sum_{j\in\cV'_{\lambda}} \frac{w_j^2}{\nu_n \ell_n}
\bigg)^{\frac{(1-2\varrho)n^{\eta}}{2\lambda^{1-\Delta}}}\cdot
w_{i_k}\, ,
\end{equation*}
which combined with \eqref{eqn:554} and \eqref{eqn:gamma2-3} yields,
$\falln$ and for all 
$\lambda\in [\lambda_{\ref{eqn:94}},\, \eps_{\ref{eqn:94}} n^\eta]$, 
\begin{align}
\E_1\big[\Gamma_{2\leadsto 3}\big]
\leq 
C\lambda^{\frac{1}{\tau-3}}\times 
\bigg(\frac{w_k}{n^\alpha}\bigg)\times 
\bigg(\sum_{j\in \cV'_{\lambda}}\frac{w_j^2}{\nu_n \ell_n}\bigg)^{\frac{(1-2\varrho)n^\eta}{2\lambda^{1-\Delta}}} \times
\bigg(\sum_{j=k+1}^n \frac{w_j^2}{\nu_n \ell_n}\bigg)^{\frac{\varrho n^\eta}{2\lambda^{1-\Delta}}}, \label{eqn:28}
\end{align}
where we have used the relation $w_{i_k}\leq w_k$. 
Writing $\wedge k$ for $k\wedge (n/2)$, we have
\begin{equation}\label{eqn:741}
\sum_{j=k+1}^n \frac{w_j^2}{\nu_n\ell_n} 
=1-\sum_{j=2}^k \frac{w_j^2}{\nu_n\ell_n}
\leq 
1-\frac{C}{n^{1-2\alpha}}\sum_{j=2}^{\wedge k}\frac{1}{j^{2\alpha}}
\leq
1-\frac{C'(\wedge k)^{\eta}}{n^{\eta}}
\leq
1-\frac{C'(k/2)^{\eta}}{n^{\eta}}
\, ,
\end{equation} 
where the first step uses the relation 
\begin{align}\label{eqn:741-a}
\sum_{j=2}^n w_j^2/(\nu_n\ell_n)=1\, ,
\end{align} 
and the second step uses Assumption \ref{ass:wts}.
Further, on the event $\fE_{n,\lambda}$,
\[
\sum_{ j\in \cC^n_1(\lambda),\, j\neq 1 } w_j^2 
\geq 
\frac{n^{2\alpha} \sigma_2^{\sss (n)}}{(1+\lambda n^{-\eta})^2}\cdot(1+\delta_0)\lambda
=
\frac{\nu_n\ell_n}{n^{\eta}(1+\lambda n^{-\eta})^2}\cdot(1+\delta_0)\lambda
\, .
\]
Thus, for $n\geq 2$ and $\lambda\in (0 ,\, \eps_{\ref{eqn:94}} n^\eta]$, on $\fE_{n,\lambda}$,
\begin{equation}\label{eqn:737}
\sum_{j\in \cV'_{\lambda}}\frac{w_j^2}{\nu_n \ell_n}
=
1-\sum_{ j\in \cC^n_1(\lambda),\, j\neq 1} \frac{w_j^2}{\nu_n \ell_n}
\leq 
1- \frac{(1+\delta_0)\lambda}{(1+\lambda n^{-\eta})^2 n^\eta} 
\leq 
1- \frac{(1+\delta_0)\lambda}{(1+\eps_{\ref{eqn:94}})^2 n^\eta}\, .
\end{equation}		 
Using \eqref{eqn:28}, \eqref{eqn:741}, \eqref{eqn:737}, and the inequality $1-u\leq e^{-u}$, we get, 
$\falln$ and for $\lambda\in [\lambda_{\ref{eqn:94}},\, \eps_{\ref{eqn:94}} n^\eta]$,
\begin{align}\label{eqn:288}
\E_1\big[\Gamma_{2\leadsto 3}\big]
\leq 
C\lambda^{\frac{1}{\tau-3}}\times 
k^{-\alpha}\times 
\exp\bigg(-\frac{(1+\delta_0)(1-2\varrho)\lambda^\Delta}{2(1+\eps_{\ref{eqn:94}})^2}\bigg)
\times
\exp\bigg(-\frac{C_{\ref{eqn:288}}k^{\eta}}{\lambda^{1-\Delta}}\bigg)
\end{align}
on the event $\fE_{n,\lambda}$.
Since $\sup_{u\geq 0} u^{\frac{\tau}{\tau-3}} e^{-u}< \infty$, and $\lambda_{\ref{eqn:94}}\geq\lambda_{\ref{eqn:ss-7a}}\geq 1$, 
\[
\exp\bigg(-\frac{C_{\ref{eqn:288}}k^{\eta}}{\lambda^{1-\Delta}}\bigg)
\leq
\exp\bigg(-\frac{C_{\ref{eqn:288}}k^{\eta}}{\lambda}\bigg)
\leq
C\bigg(\frac{\lambda}{k^{\eta}}\bigg)^{\frac{\tau}{\tau-3}}\, 
\]
for $\lambda\geq\lambda_{\ref{eqn:94}}$, which combined with \eqref{eqn:288} yields the desired result.
\qed

\medskip

\noindent {\bf Proof of Proposition \ref{prop:lem-ss-7}:}  
Combining Lemma \ref{lem:26}, \eqref{eqn:27}, and Lemma \ref{lem:fnl-bound}, we see that 
$\falln$, for all $\lambda\in [\lambda_{\ref{eqn:94}},\, \eps_{\ref{eqn:94}} n^\eta]$, $\Delta\in (0,1/2]$, 
and $1\leq k\leq |\cV'_{\lambda}|$, on the event $\fE_{n,\lambda}$,
\begin{align}\label{eqn:287}
\pr_1\bigg(\diam\big(\cC_{\res}(i_k; \lambda)\big) \geq \frac{n^\eta}{\lambda^{1-\Delta}} \bigg)
&
\leq 
C 
\lambda^{\frac{\tau+1}{\tau-3}}k^{-\frac{\tau+1}{\tau-1}}
\exp\bigg(-\, \frac{(1+\delta_0)(1-2\varrho)\lambda^\Delta}{2(1+ \eps_{\ref{eqn:94}})^2}\bigg)
\exp\bigg(\frac{(1+\delta_1)\lambda^{\Delta}}{2}\bigg) \notag\\
&
\leq 
C 
\lambda^{\frac{\tau+1}{\tau-3}}k^{-\frac{\tau+1}{\tau-1}}
\exp\bigg(-\frac{\delta_0 \lambda^{\Delta}}{4} \bigg) \, ,
\end{align}
where the last step uses the second relation in \eqref{eqn:944}.
Combining \eqref{eqn:287} with \eqref{eqn:25} and \eqref{eqn:98} completes the proof. 
\qed

\subsection{Maximum surplus outside the component of the vertex $1$}\label{sec:surplus}
Let $H_n(\lambda, \delta_1)$ and $\eps_{\ref{eqn:94}}$ be as in the setting of Proposition \ref{prop:lem-ss-7}.
Our aim in this section is to prove the following result.
\begin{prop}\label{prop:lem-ss-9}
There exists $\lambda_{\ref{eqn:96}}\geq \lambda_{\ref{it:lem-ss-4}}$ such that 
$\falln$ and $\lambda \in [\lambda_{\ref{eqn:96}},\, \eps_{\ref{eqn:94}} n^{\eta}]$,
\begin{align}\label{eqn:96}
\pr\bigg(\max\surplus\big(H_n(\lambda, \delta_1)\big)\geq 2\bigg) 
\leq 
C/\sqrt{\lambda}\, . 
\end{align}
\end{prop}

The proof of this proposition will require some results from \cite{SBSSXW14, SB-vdH-SS-PTRF}, which we will now recall briefly.
Fix a finite vertex set $\cV$ and write $\bG_{\cV}^{\con}$ for the space of all simple connected graphs with vertex set $\cV$.
For fixed $a > 0$ and probability mass function $\mvq = \big(q_v,\, v \in \cV\big)$, define a probability distribution $\pr_{\con}(\,\cdot\, ; \mvq, a, \cV)$ on $\bG_{\cV}^{\con}$ as follows: 
\begin{equation}
\label{eqn:pr-con-vp-a-cV-def}
\pr_{\con}\big(G; \mvq, a, \cV\big): = \frac{1}{Z(\mvq, a)} \prod_{(i,j)\in E(G)} \big(1-\exp(-aq_i q_j)\big) \prod_{(i,j)\notin E(G)} \exp(-aq_i q_j)\, ,\ \mbox{ for } G \in \bG_{\cV}^{\con},
\end{equation}
where $Z(\mvq,a)$ is a normalizing constant so that $\pr_{\con}(\bG_{\cV}^{\con}; \mvq, a, \cV)=1$.

For $t\geq 0$, consider the random graph $\cG\big(\big([n], \mvw^{(n)}\big),\, t\big)$ from Definition \ref{def:1}, and write 
$(\cC_i,\, i\geq 1)$ for its components in decreasing order of their masses. 
Let $\cV^{\sss(i)} := V(\cC_i)$ be the vertex set of $\cC_i$, $i \geq 1$, and note that 
$\big(\cV^{\sss(i)};\, i\geq 1\big)$ is a random partition of $[n]$. 

\begin{prop}[{\protect{\cite[Proposition 6.1]{SBSSXW14}}}]\label{prop:generate-nr-given-partition}
Conditional on the partition $\big(\cV^{\sss(i)};\,  i\geq 1\big)$, define
\begin{equation*}
\mvq^{\sss(i)} := \bigg( \frac{w_v}{\sum_{v \in \cV^{\sss(i)}}w_v } ;\, v \in \cV^{\sss(i)} \bigg)\, , \ \ 
\text{ and }\ \ 
a^{\sss(i)}:= t\bigg(\sum_{v\in \cV_{\sss(i)}} w_v\bigg)^2, \ \ \  i\geq 1.
\end{equation*}
For each fixed $i \geq 1$, let $G_i \in  \bG_{\cV^{\sss(i)}}^{\con}$ be a connected simple graph with vertex set $\cV^{\sss(i)}$. 
Then
\begin{equation*}
\pr\bigg(\cC_i = G_i \;\; \forall i \geq 1\ \big|\ \big(\cV^{\sss(i)};\ i\geq 1\big) \bigg) 
= 
\prod_{i\geq 1} \pr_{\con}\big(G_i\, ;\, \mvq^{\sss(i)}, a^{\sss(i)}, \cV^{\sss(i)}\big)\, .
\end{equation*}
\end{prop}
Thus, the random graph $\cG\big(\big([n], \mvw^{(n)}\big),\, t\big)$ can be generated in two stages:
{\bf (i) Stage I:} 
Generate the partition of the vertices into different components, i.e., generate $\big(\cV^{\sss(i)};\ i\geq 1\big)$.
{\bf (ii) Stage II:} 
Conditional on the partition, generate the internal structure of each component according to the law $\pr_{\con}\big(\, \cdot\, ;\, \mvq^{\sss(i)}, a^{\sss(i)}, \cV^{\sss(i)}\big)$ independently across different components.

In Proposition \ref{prop:SBSSXW} given below we will describe an algorithm to generate such connected components.
To state this result, we need some definitions.

For fixed $m \geq 1$, write $\bT_m$ for the set of all rooted trees with vertex set $[m]$.
Let $\bT_m^{\ord}$ be the collection of all plane trees with $m$ vertices where the vertices are labeled by elements of $[m]$.
Thus, an element of $\bT_m^{\ord}$ is a rooted tree with vertex set $[m]$ where the children of each vertex are arranged from right to left.
For $\vt\in\bT_m^{\ord}$ and $i\in[m]$, let $\sP(i; \vt)$ be the following set of vertices:
$j\in\sP(i; \vt)$ if and only if the parent of $j$ is a strict ancestor of $i$ in $\vt$, and $j$ lies on the right of the path connecting $i$ to the root of $\vt$.
Let 
\[
\sP(\vt):=
\big\{(i, j)\, :\, i\in[m]\, ,\ j\in\sP(i;\vt)  \big\}\, .
\]
For a probability mass function $\mvq=\big(q_1, \ldots, q_m\big)$ with $q_i>0$ for all $i\in[m]$ and $a>0$, define
\begin{equation}\label{eqn:ltpi-def}
L(\vt)=L(\vt; a, \mvq)
:= 
\prod_{\{i, j\}\in E(\vt)} \left[\frac{\exp(a q_i q_j)- 1}{aq_i q_j} \right] 
\exp\bigg(\sum_{(i, j) \in \sP(\vt)} a q_i q_j\bigg), \qquad \vt \in \bT_m^{\ord}.
\end{equation}

Consider $\vt\in\bT_m^{\ord}$, and suppose the vertices of $\vt$, arranged in a depth-first order, are $v(1),\ldots, v(m)$, with $v(1)$ being the root.
Define the function $f_{\vt}(\cdot)=f_{\vt}(\, \cdot\, ; a, \mvq)$ on $[0, 1]$ as follows:
\begin{align*}
f_{\vt}(s)=\sum_{j\in\sP(v(i); \vt)}q_j\, ,\ \ \ \text{ if }\ \ \ 
\sum_{k\leq i-1}q_{v(k)}\leq s<\sum_{k\leq i}q_{v(k)}\, ,
\end{align*}
and $f_{\vt}(1)=0$.
Clearly, 
\begin{align}\label{eqn:87}
\sum_{(i, j) \in \sP(\vt)} q_i q_j=\int_0^1 f_{\vt}(s)ds\, .
\end{align}
Now, since $1+s\leq e^s\leq 1+se^s$ for any $s\geq0$,
\begin{align*}
1
\leq
\prod_{\{i, j\}\in E(\vt)} \left[\frac{\exp(a q_i q_j)- 1}{aq_i q_j} \right] 
\leq
\prod_{\{i, j\}\in E(\vt)} \exp\big(a q_i q_j\big)
=
\exp\bigg(a\sum_{\{i, j\}\in E(\vt)}  q_i q_j\bigg)
\leq
\exp\big(aq_{\max}\big)\, ,
\end{align*}
where $q_{\max}=\max_{j\in [m]}q_j$.
Hence,
\begin{align}\label{eqn:91}
1
\leq
L(\vt)
\leq
\exp\big(aq_{\max}\big)\cdot
\exp\bigg(a\int_0^1 f_{\vt}(s)ds\bigg)\, .
\end{align}

Associated to the probability mass function $\mvq$ there is a random tree model called a $\mvq$-tree
\cite{pitman-camarri,pitman-random-mappings} which we now define. 
(This random tree is usually referred to as a $\mvp$-tree, but we instead use $\mvq$ to avoid confusion with $p^n_{\lambda}$ as defined in \eqref{eqn:pnlam-def}.)
For $\vt \in \bT_m$ and $v\in [m]$, write $d_v(\vt)$ for the number of children of $v$ in tree $\vt$. Then the law of the $\mvq$-tree, denoted by $P_{\text{tree}}$, is defined as follows:
\begin{equation}\label{eqn:p-tree-def}
P_{\text{tree}}(\vt) 
= 
P_{\text{tree}}(\vt; \mvq) = \prod_{v\in [m]} q_v^{d_v(\vt)}, \qquad \vt \in \bT_m.
\end{equation}
Generating a random $\mvq$-tree having distribution $P_{\text{tree}}$, and then assigning a uniform random order to the children of every vertex in this tree gives a random element in $\bT_m^{\ord}$ with law 
$P_{\ord}(\, \cdot\, ; \mvq)$ given by
\begin{equation}\label{eqn:ordered-p-tree-def}
P_{\ord}(\vt) 
= 
P_{\ord}(\vt; \mvq) = \prod_{v\in [m]} \frac{q_v^{d_v(\vt)}}{(d_v(\vt)) !}, \qquad \vt \in \bT_m^{\ord}.
\end{equation}
Using $L(\cdot)$ to tilt $P_{\ord}$ results in the following distribution:
\begin{equation}\label{eqn:tilt-ord-dist-def}
P_{\ord}^\star( \vt) 
:= 
P_{\ord}(\vt) \cdot \frac{L(\vt)}{\E[ L(T_{\mvq})]}\, , \qquad \vt \in \bT_m^{\ord}\, ,
\end{equation}
where $T_{\mvq}\sim P_{\ord}$.

\begin{prop}[{\cite[Proposition 7.4]{SBSSXW14}}]\label{prop:SBSSXW}
Fix $m\geq 1$, $a>0$, and a probability mass function $\mvq$ on $[m]$ with $\min_i q_i>0$.
Construct a random connected graph on $[m]$ as follows:
\begin{enumeratea}
\item 
First generate $T^{\star}_{\mvq}$ having distribution $P_{\ord}^\star$ as in \eqref{eqn:tilt-ord-dist-def}.
\item 
Conditional on $T^{\star}_{\mvq}$, 
add the edge $\{i, j\}$ with probability $1-\exp\big(-aq_i q_j\big)$
independently for $(i, j)\in\sP(T^{\star}_{\mvq})$.
\end{enumeratea}
Then the resulting random graph is distributed as $\pr_{\con}\big(\,\cdot\, ;\mvq, a, [m]\big)$ on $\bG_{[m]}^{\con}$.
\end{prop}

Note that in the above construction, the number of  surplus edges in the resulting graph equals the number of edges added in part (b) of the construction.
Using \eqref{eqn:88} and \eqref{eqn:87}, we see that the number of surplus edges in the random graph resulting in the above construction is stochastically dominated by $Y$, where 
\[
\big(Y\, \big|\, T^{\star}_{\mvq}\big)
\sim
\Poi\big(a\int_0^1 f_{T_{\mvq}^{\star}}(s)ds\big)\, .
\]
Hence, if $\cG_{\con}$ is distributed as $\pr_{\con}\big(\,\cdot\, ;\mvq, a, [m]\big)$, then by Proposition \ref{prop:SBSSXW},
\begin{align}\label{prop:surplus}
\pr\big(\surplus\big(\cG_{\con}\big) \geq 2\big) 
&\leq 
\pr\big(Y\geq 2\big)
\leq
a^2\cdot\bE\big[\|f_{T_{\mvq}^{\star}}\|_{\infty}^2\big]
=
a^2\cdot\frac{\bE\big[\|f_{T_{\mvq}}\|_{\infty}^2\cdot L(T_{\mvq})\big]}{\bE\big[L(T_{\mvq})\big]}
\notag\\&\leq
a^2\exp\big(aq_{\max}\big)\cdot
\bE\big[\|f_{T_{\mvq}}\|_{\infty}^2\cdot\exp\big(a\|f_{T_{\mvq}}\|_{\infty}\big)\big]\, ,
\end{align}
where $T_{\mvq}\sim P_{\ord}$.
Here, the second step uses the fact that $\pr\big(\Poi(s)\geq 2\big)\leq s^2$ for all $s>0$, and the final step uses \eqref{eqn:91}.
The next result gives a tail bound on $\|f_{T_{\mvq}}\|_{\infty}$.

\begin{lemma}\label{lem:tilt-bound}
There exists a universal constant $C > 0$ such that for every $m\geq 1$, probability mass function $\mvq$ on $[m]$ with $\min_i q_i>0$, and $x\geq \e$,
\[
\pr\left(\|f_{T_{\mvq}}\|_{\infty} \geq x ||\mvq||_2\right)\leq \exp\big(-C x\log(\log{x})\big)\, .
\]
\end{lemma}

Lemma \ref{lem:tilt-bound} follows by combining 
\cite[Lemma 7.9]{SBSSXW14} and \cite[Lemma 4.9]{SB-vdH-SS-PTRF}.
We are now ready to prove Proposition \ref{prop:lem-ss-9}.

\medskip

\noindent{\bf Proof of Proposition \ref{prop:lem-ss-9}:}
Let $\pr_1(\cdot)$ and $\E_1(\cdot)$ be as in \eqref{eqn:45}.
We will write $\sum_{1}$ for $\sum_{j\in\cV'_{\lambda}}$, where $\cV'_{\lambda}$ is as in \eqref{eqn:466}. 
Let $\fE_{n,\lambda}$ be as defined before \eqref{eqn:98}.
Then, for $n\geq 2$ and $\lambda\in (0 ,\, \eps_{\ref{eqn:94}} n^\eta]$, on the event $\fE_{n,\lambda}$,
\begin{align}\label{eqn:77}
p^n_{(1+\delta_1)\lambda}\cdot\sum\displaystyle_1 \ch{\frac{w_j^2}{\ell_n}}
&=
\bigg(1+\frac{(1+\delta_1)\lambda}{n^{\eta}}\bigg)
\sum\displaystyle_1 \frac{w_j^2}{\nu_n\ell_n}
\leq
\bigg(1+\frac{(1+\delta_1)\lambda}{n^{\eta}}\bigg)
\bigg(1-\frac{(1+\delta_0)\lambda}{(1+\eps_{\ref{eqn:94}})^2 n^\eta}\bigg)
\notag\\
&\leq
1+\frac{\lambda}{n^{\eta}}
\bigg(\big(1+\delta_1\big)-\frac{(1+\delta_0)}{(1+\eps_{\ref{eqn:94}})^2}\bigg)
\leq 
1-\frac{\lambda\delta_0}{2n^{\eta}}\, ,
\end{align}
where the second step uses \eqref{eqn:737}, and the final step uses the second relation in \eqref{eqn:944},
Consequently,
\begin{align}\label{eqn:33}
&\pr_1\bigg(
\sum_{\tiny \substack{\cC \text{component} \\ \text{of }H_n(\lambda, \delta_1)}} \hspace{-.2in} \big(\cW(\cC)\big)^2 \geq \frac{n^{2\rho}}{\sqrt{\lambda}} \bigg)	
\leq 
\frac{\sqrt{\lambda}}{n^{2\rho}}\cdot 
\E_1\bigg[\sum\displaystyle_1 w_v\cdot \cW\bigg(\cC\big(v, H_n(\lambda, \delta_1)\big)\bigg)\bigg]\\	
&\hskip40pt
\leq 
\frac{\sqrt{\lambda}}{n^{2\rho}}\cdot
\sum\displaystyle_1 w_v\cdot \bigg(\frac{w_v}{1-p^n_{(1+\delta_1)\lambda}\cdot\sum\displaystyle_1 \ch{w_j^2/\ell_n}}\bigg) 
\leq
\frac{C\sqrt{\lambda}}{n^{2\rho}}\cdot
\sum\displaystyle_1 w_v^2\cdot \bigg(\frac{n^{\eta}}{\lambda}\bigg) 
\leq\frac{C'}{\sqrt{\lambda}}\, ,\notag
\end{align} 
where the first step uses Markov's inequality, 
the second step uses \cite[Lemma 4.5]{addario2017probabilistic},
the third step uses \eqref{eqn:77},
and the final step uses the relation $\sum_{j=1}^n w_j^2=O(n)$.

For $\lambda \geq \lambda_{\ref{it:lem-ss-4}}$, define
\begin{equation}\label{eqn:34}
\fE_{n,\lambda}' 
= 
\fE_{n,\lambda} \cap \bigg\{j\in \cC^n_1(\lambda)\ \text{for all }\ 
1\leq j\leq \frac{\lambda^{1/\eta}}{\log^3\lambda}\bigg\}\, .
\end{equation}
For $\lambda>0$, let
\begin{equation}\label{eqn:35}
\fB_{n,\lambda}
:= 
\bigg\{
\sum_{\tiny \substack{\cC \text{component} \\ \text{of }H_n(\lambda, \delta_1)}} \hspace{-.2in} \big(\cW(\cC)\big)^2 < n^{2\rho}\lambda^{-1/2} 
\bigg\}\, ,
\end{equation}
and note that $\fB_{n,\lambda}$ is the complement of the event studied in \eqref{eqn:33}. 
On the event $\fB_{n,\lambda}$, 
\begin{equation}\label{eqn:35A}
\cW(\cC) \leq n^{\rho}\lambda^{-1/4} \ \ \text{ for all components }
\cC \text{ of }  H_n(\lambda, \delta_1)\, . 
\end{equation}
Let $i_1<i_2<\ldots$ and $\cC_{\res}(i_k;\lambda)$ be as defined before \eqref{eqn:45}.
Then 
\begin{align}\label{eqn:222}
&\pr_1\bigg(
\sum_{j\in \cC} \theta_{j,\lambda}^2 
\geq 
\theta_{1,\lambda}^2+2(1+\delta_1)\sum_{j=2}^n \theta_{j,\lambda}^3\ \ \ 
\text{for some component } \cC \text{ of } H_n(\lambda, \delta_1)
\bigg)\\
&\hskip50pt\leq
\pr_1\big(\fB_{n,\lambda}^c\big)
+
\sum_{k}\pr_1\bigg(
\bigg\{\sum_{j\in \cC_{\res}(i_k;\lambda)} \theta_{j,\lambda}^2 
\geq 
\theta_{1,\lambda}^2+2(1+\delta_1)\sum_{j=2}^n \theta_{j,\lambda}^3\bigg\}
\bigcap \fB_{n,\lambda}
\bigg)\, .\notag
\end{align}
Using \eqref{eqn:35A}, we can choose $\lambda_{\ref{eqn:96}}\geq\lambda_{\ref{it:lem-ss-4}}$ such that $\falln$ and for $\lambda\geq\lambda_{\ref{eqn:96}}$, on the event $\fB_{n,\lambda}$,
\begin{align}\label{eqn:81}
\cW_{\mvx}\big(\cC_{\res}(i_k;\lambda)\big)\leq 1\ \ \text{ for }\ \ 1\leq k\leq |\cV'_{\lambda}|\, .
\end{align}
Analogous to the exploration mentioned before Lemma \ref{lem:domin-ald-limic-1}, conditional on $\cC^n_1(\lambda)$, $\cC_{\res}(i_k;\lambda)$ can be explored in a breadth-first manner starting from $i_k$, and $\cW_{\mvx}\big(\cC_{\res}(i_k;\lambda)\big)$ will have the same distribution as the hitting time of zero by the corresponding breadth-first walk.
Let $\xi^n_{j, (1+\delta_1)\lambda}$, $2\leq j\leq n$, be as in \eqref{eqn:713} with $\lambda$ replaced by $(1+\delta_1)\lambda$.
Using \eqref{eqn:81}, $\falln$ and for $\lambda\in[\lambda_{\ref{eqn:96}},\, \eps_{\ref{eqn:94}} n^\eta]$,
for every $1\leq k\leq |\cV'_{\lambda}|$,
\begin{align}\label{eqn:85}
&\pr_1\bigg(
\bigg\{\sum_{j\in \cC_{\res}(i_k;\lambda)} \theta_{j,\lambda}^2 
\geq 
\theta_{1,\lambda}^2+2(1+\delta_1)\sum_{j=2}^n \theta_{j,\lambda}^3\bigg\}
\bigcap \fB_{n,\lambda}
\bigg)\\
&\hskip30pt
\leq
\pr_1\bigg(
\theta_{i_k,\lambda}^2
+
\sum_{j=i_{k}}^n\theta_{j,\lambda}^2\ind\big\{\xi^n_{j, (1+\delta_1)\lambda}\leq 1\big\}
\geq 
\theta_{1,\lambda}^2 + 2(1+\delta_1)\sum_{j=2}^n \theta_{j,\lambda}^3
\bigg)\notag\\
&\hskip60pt
\leq
\exp\bigg(
-\frac{C}{\sum_{j=i_k}^n\theta_{j, \lambda}^4}
\bigg)
\leq
\exp\bigg(-\frac{C'}{\sum_{j\geq i_k}j^{-4\alpha}}\bigg)
\leq
\exp\bigg(-C''i_k^{4\alpha-1} \bigg)
\notag\, ,
\end{align}
where the second step follows from the bounded difference inequality.
Combining \eqref{eqn:85} with \eqref{eqn:222} and \eqref{eqn:33}, we see that 
$\falln$ and for $\lambda\in[\lambda_{\ref{eqn:96}},\, \eps_{\ref{eqn:94}} n^\eta]$,
on the event $\fE_{n,\lambda}'$,
\begin{align}\label{lem:36}
&\pr_1\bigg(
\sum_{j\in \cC} \theta_{j,\lambda}^2 
\geq 
\theta_{1,\lambda}^2+2(1+\delta_1)\sum_{j=2}^n \theta_{j,\lambda}^3\ \ \ 
\text{for some component } \cC \text{ of } H_n(\lambda, \delta_1)
\bigg)\notag\\
&\hskip70pt
\leq
\frac{C}{\sqrt{\lambda}}
+
\sum_{k\geq\lambda^{1/\eta}/\log^3\lambda}\exp\bigg(-C'k^{4\alpha-1} \bigg)
\leq\frac{C''}{\sqrt{\lambda}}\, .
\end{align}

Now, using Proposition \ref{prop:generate-nr-given-partition}, conditional on $\cC_1^n(\lambda)$, we can generate the components of $H_n(\lambda, \delta_1)$ in two steps: 
first generate the random partition $\big(\cV^{\sss(i)};\, i\geq 1\big)$ of $\cV'_{\lambda}$ into component vertex sets, and then conditionally generate the internal structure of each component.
Here, the relevant parameters are
\begin{equation}\label{eqn:vpi-chpn}
\mvq^{\sss(i)} 
=
\big(q^{\sss(i)}_v;\, v\in \cV^{\sss(i)}\big)
:= 
\bigg(\frac{w_v}{\cW(\cV^{\sss(i)})};\, v\in \cV^{\sss(i)}\bigg)\, , 
\ \ \text{ and }\ \ 
a^{\sss(i)} = p^n_{(1+\delta_1)\lambda} \big(\cW(\cV^{\sss(i)})\big)^2\big/\ell_n\, .
\end{equation}
We will write
\[
\pr_{\partition}(\cdot):=\pr\big(\,\cdot\, \big|\, \cC^n_1(\lambda)\, ,\ (\cV^{\sss(i)};\, i\geq 1)\big)\, ,
\ \ \text{ and }\ \ 
\E_{\partition}(\cdot):=\E\big(\,\cdot\, \big|\, \cC^n_1(\lambda)\, ,\ (\cV^{\sss(i)};\, i\geq 1)\big)\, .
\]
For $\lambda>0$, let
$\fB_{n,\lambda}' 
:= 
\fB_{n,\lambda} \cap \big\{
\sum_{j\in \cC} w_j^2 
\leq C_{\ref{eqn:67}}n^{2\alpha}
\text{ for each component }\cC\text{ of } H_n(\lambda, \delta_1)
\big\}
$,
where
\begin{align}\label{eqn:67}
C_{\ref{eqn:67}}
:=
\sup\,\bigg\{
\big(\theta_{1,\eps_{\ref{eqn:94}} n^\eta}^2+2(1+\delta_1)\sum_{j=2}^n \theta_{j,\eps_{\ref{eqn:94}} n^\eta}^3\big)\cdot \sigma_2^{\sss (n)}
\, :\, 
n\geq 1\bigg\}\, .
\end{align}
Writing
$
q_{\max}^{\sss (i)}
:=
\max \big\{q^{\sss(i)}_v\, :\, v\in \cV^{\sss(i)}\big\}
$,
we get, $\falln$ and for $\lambda\in(0,\, \eps_{\ref{eqn:94}} n^\eta]$, 
on the event $\fB_{n,\lambda}'$,  for $i\geq 1$,
\begin{align}
a^{\sss(i)}q_{\max}^{\sss (i)}
&\leq a^{\sss(i)}\|\mvq^{\sss(i)}\|_2 
= p^n_{(1+\delta_1)\lambda}\cdot
\frac{\big(\cW(\cV^{\sss(i)})\big)^2}{\ell_n}\cdot
\frac{\big(\sum_{v\in \cV^{\sss(i)}}w_v^2\big)^{1/2}}{\cW(\cV^{\sss(i)})} \notag \\ 
&
\leq 
C\cdot \frac{\cW(\cV^{\sss(i)})}{n}\cdot n^{\alpha} 
=
C\cdot\frac{\cW(\cV^{\sss(i)})}{n^{\rho}}
\label{eqn:40}\\
&\hskip40pt
\leq 
\frac{C}{\lambda^{1/4}}\, , \label{eqn:40A}
\end{align}
\ch{where the last step uses \eqref{eqn:35A}.}
Hence, $\falln$ and for $\lambda\in[1,\, \eps_{\ref{eqn:94}} n^\eta]$, 
on the event $\fB_{n,\lambda}'$,  for $i\geq 1$,
\begin{align}
&\pr_{\partition}\big(\surplus(\cC_i) \geq 2\big) 
\leq
\big(a^{\sss(i)}\big)^2\cdot\exp\big(a^{\sss(i)}q_{\max}^{\sss (i)}\big)\cdot
\bE\big[
\|f_{T_{\mvq^{\sss(i)}}}\|_{\infty}^2\cdot\exp\big(a^{\sss(i)} \|f_{T_{{\mvq}^{\sss(i)}}}\|_{\infty} \big)
\big]
\notag\\
&\hskip30pt
=
\big(a^{\sss(i)}\|\mvq^{\sss(i)}\|_2\big)^2\cdot
\exp\big(a^{\sss(i)}q_{\max}^{\sss (i)}\big)\cdot
\bE\bigg[
\frac{\|f_{T_{\mvq^{\sss(i)}}}\|_{\infty}^2}{\|\mvq^{\sss(i)}\|_2^2}\cdot
\exp\bigg(
a^{\sss(i)}\|\mvq^{\sss(i)}\|_2\cdot
\frac{\|f_{T_{{\mvq}^{\sss(i)}}}\|_{\infty} }{\|\mvq^{\sss(i)}\|_2}
\bigg)
\bigg]
\notag\\
&\hskip60pt
\leq
C\big(\cW(\cV^{\sss(i)})\big)^2n^{-2\rho}\, ,
\label{eqn:38}
\end{align}
where the first step uses \eqref{prop:surplus}, 
and the final step uses Lemma \ref{lem:tilt-bound} together with the bounds in both \eqref{eqn:40} and \eqref{eqn:40A}.
Consequently, $\falln$ and for $\lambda\in[1,\, \eps_{\ref{eqn:94}} n^\eta]$, 
on the event $\fB_{n,\lambda}'$,
\[
\pr_{\partition}\bigg(
\max\surplus\big(H_n(\lambda, \delta_1)\big)\geq 2
\bigg) 
\leq 
C\sum_i \frac{\big(\cW(\cV^{\sss(i)})\big)^2}{n^{2\rho}}
\leq
\frac{C}{\sqrt{\lambda}}.
\]
Thus, to complete the proof it is enough to show that 
$\falln$ and $\lambda \in [\lambda_{\ref{eqn:96}},\, \eps_{\ref{eqn:94}} n^{\eta}]$,
\begin{equation}\label{eqn:424}
\pr\big((\fB_{n,\lambda}')^c\big) \leq C\lambda^{-1/2}.
\end{equation}
To this end, observe that 
$\falln$ and for $\lambda \in [\lambda_{\ref{eqn:96}},\, \eps_{\ref{eqn:94}} n^{\eta}]$,
\begin{align*}
\pr\big((\fB_{n,\lambda}')^c\big) 
&\leq 
\pr\big((\fE_{n,\lambda}')^c\big) 
+ 
\E\big[\ind_{\fE_{n,\lambda}'}\cdot \pr_1\big((\fB_{n,\lambda}')^c\big)\big]\\
&\leq 
\pr\big((\fE_{n,\lambda}')^c\big) + C\lambda^{-1/2}
\leq 
\exp\big(-C'\log^{3\alpha}\lambda\big) + \exp(-C'\lambda)+ C\lambda^{-1/2}\, ,
\end{align*}
where the second inequality uses \eqref{eqn:33} and \eqref{lem:36},
and the final inequality uses Proposition \ref{prop:lem-ss-4}.
This completes the proof of \eqref{eqn:424}, and hence of Proposition \ref{prop:lem-ss-9}.
\qed

\subsection{Proof of Proposition \ref{prop:ss-10}}\label{sec:scanning}
Let $\delta_1$ be as chosen in \eqref{eqn:944}.
Let $\delta_2$ and $\eps_{\ref{eqn:333}}$ be given by
\[
\delta_2:=\delta_1/2\, ,\ \ \text{ and }\ \ 
\eps_{\ref{eqn:333}}:=\eps_{\ref{eqn:94}}/2\, .
\]
For $\lambda\geq 0$, let $\tbar H_n(\lambda, \delta_2)$ be the random graph on the vertex set
$\tbar \cV'_{\lambda}:=[n]\setminus V\big(\tbar\cC_1^n(\lambda)\big)$ 
obtained by placing edges independently between $i, j\in\tbar\cV'_{\lambda}$ with probability
$1\wedge\big(p^n_{(1+\delta_2)\lambda}w_i w_j/\ell_n\big)$.

Since
$1-e^{-u}\leq u\wedge 1$ for all $u\geq 0$, 
$\cC^n_1(\lambda)$ can be coupled with $\tbar\cC^n_1(\lambda)$ so that the former is a subgraph of the latter.
Further, we can choose 
$\lambda_{\ref{eqn:80}}\geq \lambda_{\ref{eqn:94}}\vee\lambda_{\ref{eqn:96}}$ 
such that $\falln$ and 
$\lambda \in [\lambda_{\ref{eqn:80}},\, \eps_{\ref{eqn:333}} n^{\eta}]$, 
\begin{align}\label{eqn:80}
p^n_{(1+\delta_2)\lambda} w_i w_j/\ell_n 
\leq 
1-\exp\big(-p^n_{(1+\delta_1)\lambda} w_i w_j/\ell_n\big)
\end{align}
for all $i,j\in [n]$.
Then it follows that $\falln$ and 
$\lambda \in [\lambda_{\ref{eqn:80}},\, \eps_{\ref{eqn:333}} n^{\eta}]$, 
\begin{equation}\label{eqn:41}
d_{\rH}\bigg(
\tbar M^n_{\lambda}\, ,\, \tbar M^n_{(1+\delta_2)\lambda}
\bigg) 
\stod\,
1+\mathrm{LP}\big(\tbar H_n(\lambda, \delta_2)\big)
\stod\, 
1+\mathrm{LP}\big(H_n(\lambda, \delta_1)\big)\, ,
\end{equation}
where $H_n(\lambda, \delta_1)$ is as in the setting of Proposition \ref{prop:lem-ss-7}, 
\ch{and $\mathrm{LP}(\cdot)$ is as defined below \eqref{eqn:234}}.
It is easy to check that for any finite graph $H$, 
\begin{equation}\label{eqn:50}
1+\mathrm{LP}(H) 
\leq 
1+8\cdot\diam(H) \cdot \big(\max\surplus(H)+1\big)
\leq
8\big(\diam(H)+1\big)\big(\max\surplus(H)+1\big)\, .
\end{equation}
Combining \eqref{eqn:41} and \eqref{eqn:50}, we see that for any $\Delta\in(0, 1/2]$, $\falln$ and 
$\lambda \in [\lambda_{\ref{eqn:80}},\, \eps_{\ref{eqn:333}} n^{\eta}]$, 
\begin{align}\label{eqn:90}
&\pr\bigg(
d_{\rH}\bigg(
\tbar M^n_{\lambda}\, ,\, \tbar M^n_{(1+\delta_2)\lambda}
\bigg) 
\geq 
\frac{24 n^{\eta}}{\lambda^{1-\Delta/2}}
\bigg)\\
&\hskip30pt\leq 
\pr\bigg(1+\diam\big(H_n(\lambda, \delta_1)\big) \geq \frac{n^\eta}{\lambda^{1-\Delta/2}}\bigg) 
+ 
\pr\bigg(\max\surplus\big(H_n(\lambda, \delta_1)\big)\geq 2\bigg)\notag\\
&\hskip30pt\leq 
\pr\bigg(\diam\big(H_n(\lambda, \delta_1)\big) \geq \frac{n^\eta}{(2\lambda)^{1-\Delta/2}}\bigg) 
+ 
\pr\bigg(\max\surplus\big(H_n(\lambda, \delta_1)\big)\geq 2\bigg)\notag\\
&
\hskip30pt\leq
C \bigg(
\lambda^{\frac{\tau+1}{\tau-3}}\exp\big(-C' \lambda^{\Delta/2}\big)
+\lambda^{-1/2}
\bigg)
\leq
C\lambda^{-1/2} 
\sup_{u\geq 1}\bigg[ u^{\frac{1}{2}+\frac{\tau+1}{\tau-3}}\exp\big(-C' u^{\Delta/2}\big)+1\bigg]
\notag
\, ,
\end{align}
where the third inequality follows from Proposition \ref{prop:lem-ss-7} and Proposition \ref{prop:lem-ss-9}.
Now for any $\lambda_{\ref{eqn:80}}\leq \lambda\leq \eps_{\ref{eqn:333}} n^{\eta}$, 
\begin{equation}\label{eqn:826}
d_{\rH}\bigg(
\tbar M^n_{\lambda}\, ,\, \tbar M^n_{\eps_{\ref{eqn:333}} n^{\eta}}
\bigg) 
\leq 
\sum_{j=0}^{k_0} 
d_{\rH}\bigg(
\tbar M^n_{(1+\delta_2)^j\lambda}\, ,\, \tbar M^n_{(1+\delta_2)^{j+1}\lambda}
\bigg)\, ,
\end{equation}
where 
$k_0 = k_0(n,\lambda) = 
\min\big\{j\geq 0\, :\, (1+\delta_2)^{j+1}\lambda \geq \eps_{\ref{eqn:333}} n^{\eta}\big\}$.
Using \eqref{eqn:826} and \eqref{eqn:90}, the proof of Proposition \ref{prop:ss-10} can be completed via a simple union bound.
Note that
\begin{align}\label{eqn:101}
\sum_{j=0}^{k_0}\frac{24n^{\eta}}{\big((1+\delta_2)^j \lambda\big)^{1-\Delta/2}}
\leq 
\frac{24n^{\eta}}{\lambda^{1-\Delta/2}} \sum_{j=0}^\infty \frac{1}{(1+\delta_2)^{3j/4}} 
=
\frac{C_{\ref{eqn:101}}n^{\eta}}{\lambda^{1-\Delta/2}}
\leq 
\frac{n^{\eta}}{\lambda^{1-\Delta}}\, ,
\end{align}
where the last step is true for
$\lambda\geq \lambda_{\ref{eqn:555}}=\lambda_{\ref{eqn:555}}(\Delta)$, 
and $\lambda_{\ref{eqn:555}}\geq\lambda_{\ref{eqn:80}}$ is chosen in a way so that 
$C_{\ref{eqn:101}}\leq \lambda_{\ref{eqn:555}}^{\Delta/2}$.
We then choose $n_{\ref{eqn:333}}(\Delta)\geq 1$ so that the interval 
$[\lambda_{\ref{eqn:555}}(\Delta),\, \eps_{\ref{eqn:333}} n^{\eta}]$ is nonempty for 
$n\geq n_{\ref{eqn:333}}(\Delta)$.
The rest of the argument is routine.
\qed

\subsection{Proof of Proposition \ref{prop:ss-11}}\label{sec:45}
To simplify notation, we will write 
\[
\gamma_n = \eps_{\ref{eqn:333}} n^{\eta}\, .
\] 
Before starting the proof of Proposition \ref{prop:ss-11}, we need two elementary lemmas.
	
\begin{lemma}\label{lem:ss-lem-12}
For $C>0$, write $\fC_n(C)$ for the event that 
$\cW(\cC) \leq C\log{n}$ for all components $\cC$ of $\barGn(\gamma_n)$ other than $\tbar\cC^n_1(\gamma_n)$.
Then for every $\kappa >0$, there exists $C_{\ref{eqn:30}}=C_{\ref{eqn:30}}(\kappa) <\infty$ 
and $n_{\ref{eqn:30}}=n_{\ref{eqn:30}}(\kappa)$ such that for $n\geq n_{\ref{eqn:30}}$,
\begin{align}\label{eqn:30}
\pr\big(\fC_n(C_{\ref{eqn:30}})\big) \geq 1-n^{-\kappa}\, .
\end{align}
\end{lemma}
The proof of Lemma \ref{lem:ss-lem-12} will be given in Section \ref{sec:proof-ss-lem-12}.

\begin{lemma}\label{lem:ss-lem17}
{\bf\upshape (a)}
Suppose $k, m_0\in\bZ_{>0}$, and 
$x_0, x_1, \ldots, x_k \in (0,1)$ with $x_0\leq \min\big\{x_j\, :\,1\leq j\leq k\big\}$. 
Further, assume that 
$Z_i^{\sss(0)}$, $1\leq i\leq m_0$, and $Z^{\sss(1)}, \ldots, Z^{\sss(k)}$ are independent random variables such that 
$Z_i^{\sss(0)} \sim \Unif\,[x_0,1]$, $i=1,\ldots, m_0$, and
$Z^{\sss(j)} \sim \Unif\, [x_j,1]$, $j=1,\ldots, k$.
Then 
\[
\pr\bigg(\min_{1\leq i\leq m_0} Z_i^{\sss(0)} < \min_{1\leq j\leq k} Z^{\sss(j)}\bigg)
\geq 
\frac{m_0}{m_0+k}. 
\]

\vskip4pt

\noindent{\bf\upshape (b)}
Suppose in addition to the parameters and random variables in (a), 
we have, for some $\tbar k, \tbar m_0\in\bZ_{>0}$, another collection of numbers 
$\tbar x_0, \tbar x_1, \ldots, \tbar x_{\tbar k} \in (0,1)$ with $\tbar x_0\leq \min\big\{\tbar x_l\, :\,1\leq l\leq \tbar k\big\}$. 
Further, assume that 
$\tbar Z_i^{(0)}$, $1\leq i\leq \tbar m_0$, and $\tbar Z^{(1)}, \ldots, \tbar Z^{(\tbar k)}$ are independent random variables that are also independent of the collection of random variables in {\upshape (a)} such that 
$\tbar Z_i^{\sss(0)} \sim \Unif\, [\tbar x_0,1]$, $i=1,\ldots, \tbar m_0$, and
$\tbar Z^{\sss(l)} \sim \Unif\, [\tbar x_l,1]$, $l=1,\ldots, \tbar k$.
Then 
\[
\pr\bigg(
\min_{1\leq i\leq m_0} Z_i^{\sss(0)} \wedge \min_{1\leq i\leq \tbar m_0} \tbar{Z}_i^{\sss(0)} 
< 
\min_{1\leq j\leq k} Z^{\sss(j)} \wedge \min_{1\leq l\leq\tbar k} \tbar{Z}^{(l)}  
\bigg)
\geq 
\frac{m_0}{m_0+k} \wedge \frac{\tbar{m}_0}{\tbar{m}_0+\tbar k}\, . 
\]
\end{lemma}

\noindent{\bf Proof:}
We only prove (b). 
We start with a specific construction of the random variables in (b). 
First generate $m_0+k$ i.i.d. $\Unif\, [x_0,1]$ random variables 
$Y_1^{\sss(0)}, \ldots, Y_{m_0}^{\sss(0)}, Y^{\sss(1)}, \ldots, Y^{\sss(k)}$,
and independent of this collection, generate 
$\tbar m_0+\tbar k$ i.i.d. $\Unif\, [\tbar x_0,1]$ random variables 
$\tbar Y_1^{(0)}, \ldots, \tbar Y_{\tbar m_0}^{(0)}, \tbar Y^{(1)}, \ldots, \tbar Y^{(\tbar k)}$.
Since 
$x_0\leq \min\big\{x_j\, :\,1\leq j\leq k\big\}$ and
$\tbar x_0\leq \min\big\{\tbar x_l\, :\,1\leq l\leq \tbar k\big\}$,
we can construct the random variables in (b) on the same probability space as the above collection of random variables such that
\begin{gather*}
Z_i^{\sss(0)} = Y_i^{\sss(0)}\ \ \text{ for }\ \ i=1,\ldots, m_0\, ,\ \ \ 
\tbar{Z}_i^{(0)} = \tbar{Y}_i^{(0)}\ \ \text{ for }\ \ i=1, \ldots, \tbar{m}_0\, ,\\
Z^{\sss(j)} \geq Y^{\sss(j)}\ \ \text{ for }\ \ j=1, \ldots, k\, , 
\ \ \text{ and }\ \ 
\tbar{Z}^{(l)} \geq \tbar{Y}^{(l)}\ \ \text{ for }\ \ l=1, \ldots, \tbar k\, .
\end{gather*}
Then note that the probability in Lemma \ref{lem:ss-lem17} (b) is lower bounded by 
\[
\pr\bigg(
\min_{1\leq i\leq m_0} Y_i^{\sss(0)} \wedge \min_{1\leq i\leq \tbar m_0} \tbar{Y}_i^{\sss(0)} 
< 
\min_{1\leq j\leq k} Y^{\sss(j)} \wedge \min_{1\leq l\leq\tbar k} \tbar{Y}^{(l)}  
\bigg)
=:
\pr(\fM)\, ,
\]
say. 
Let
$\mathpzc{Min}
:=
\min\big\{Y_1^{\sss(0)}, \ldots, Y_{m_0}^{\sss(0)}, Y^{\sss(1)}, \ldots, Y^{\sss(k)}\big\}$,
and similarly define $\overline{\mathpzc{Min}}$.
By symmetry, 
\[
\pr\big(\fM\, \big|\, {\mathpzc{Min}} < \overline{\mathpzc{Min}}\big) =\frac{m_0}{m_0+k}\, ,
\ \ \text{ and }\ \  
\pr\big(\fM\, \big|\, \overline{\mathpzc{Min}} < {\mathpzc{Min}}\big) = \frac{\tbar{m}_0}{\tbar{m}_0 +\tbar k}\, .
\]
Hence the claim follows. 
\qed

\medskip

Given $\barGn(\gamma_n)$ and the edge weights $\mvU\big|_{\barGn(\gamma_n)}$, we can construct the restriction of $\tbar M^n$ to $\barGn(\gamma_n)$.
For any vertex $v\notin V\big(\tbar\cC^n_1(\gamma_n)\big)$, in order to find the path in $\tbar M^n$ that connects $\tbar M^n_{\gamma_n}$ to the restriction of $\tbar M^n$ to $\tbar\cC^n_v(\gamma_n)$, we can proceed via the following algorithm. 

\vskip5pt

\noindent {\bf Algorithm 1:} In this algorithm, we will join connected components of $\barGn(\gamma_n)$  sequentially using edges from $E(\barGn)\setminus E(\barGn(\gamma_n))$. 
We will refer to a collection of connected components joined by such edges as a ``cluster." 
	
\begin{enumeratea}
\item 
Look at the edges in $E(\barGn)\setminus E(\barGn(\gamma_n))$ going out of $\tbar\cC^n_v(\gamma_n)$. 
Add the edge with the minimum weight to $\tbar\cC^n_v(\gamma_n)$, thereby connecting it to another component of $\barGn(\gamma_n)$. 
\item 
Repeat sequentially with the current cluster of $v$. 
At the $k+1$-th step, we look at the edges (if any) that 
are in $E(\barGn)\setminus E(\barGn(\gamma_n))$, and
have one endpoint in the current cluster of $v$ (i.e., the cluster after the addition of the $k$-th edge) and one endpoint outside of this cluster.
We then add the edge with the minimum weight among all the outgoing edges from the cluster, i.e., including the ones considered in the first $k$ steps. 
\item 
Stop if $v$ gets connected to $\tbar\cC^n_1(\gamma_n)$, or if there are no outgoing edges from the current cluster of $v$. 
In the latter case, $v\notin \cC(1, \barGn)$. 
\end{enumeratea}

Using Lemma \ref{lem:mst-minimax-criterion}, it is easy to see that if $v\in \cC(1, \barGn)$, then every edge added in Algorithm 1 is an edge in $\tbar M^n$.
Now, the two-layered randomness present in the problem (presence/absence of edges and edge weights) makes Algorithm 1 hard to analyze directly. 
Loosely speaking, the problem arises from the fact that if $\cW(\cC)$ is small for some component $\cC$ of $\barGn(\gamma_n)$ other than $\tbar\cC^n_1(\gamma_n)$, 
then we cannot say that there are enough edges in $E(\barGn)\setminus E(\barGn(\gamma_n))$ that connect $\cC$ to $\tbar\cC^n_1(\gamma_n)$ with high probability-- something we need for the argument to work.
So instead we will work with the following modified algorithm. 
Let
\begin{align}\label{eqn:tn-def}
\pzt_n:=\frac{n^{\eta/4}}{(\log n)^{1/3}}\, .
\end{align}

\noindent{\bf Algorithm 2:} 
The word ``cluster" will have the same meaning as in Algorithm 1.
\begin{enumeratea}
\item 
Look at the edges in $E(\barGn)\setminus E(\barGn(\gamma_n))$ that connect $\tbar\cC^n_v(\gamma_n)$ to $[n]\setminus \big(V\big(\tbar\cC^n_v(\ch{\gamma_n})\big)\cup V\big(\tbar\cC^n_1(\gamma_n)\big)\big)$,
i.e., we do not check for edges that connect $\tbar\cC^n_v(\gamma_n)$ to $\tbar\cC^n_1(\gamma_n)$.
Add the edge with the minimum weight to $\tbar\cC^n_v(\gamma_n)$, thereby connecting it to another component of $\barGn(\gamma_n)$. 
\item \label{it:check} 
Repeat for another $(\pzt_n-1)$ many steps in a manner similar to Algorithm 1, but without checking for edges that connect the current cluster to 
$\tbar\cC^n_1(\gamma_n)$, or until there are no outgoing edges from the current cluster to vertices that are not in $\tbar\cC^n_1(\gamma_n)$.
\item 
After the $\pzt_n$-th edge has been added, if the algorithm has not terminated already, look at the edges (if any) that 
are in $E(\barGn)\setminus E(\barGn(\gamma_n))$, and
have one endpoint in the current cluster and one endpoint outside of this cluster.
Let us emphasize that at this step, we are checking for edges that connect the current cluster to $\tbar\cC^n_1(\gamma_n)$.
We then add the edge with the minimum weight among all the outgoing edges from the current cluster, thereby connecting it to another component of $\barGn(\gamma_n)$. 
Stop if this component is $\tbar\cC^n_1(\gamma_n)$. 
\item 
Else, ignore the edges (if any) found at the $(\pzt_n+1)$-th step between $\tbar\cC^n_1(\gamma_n)$ and the cluster of $v$, and continue as in \eqref{it:check} for another $\pzt_n$ many steps.
Thus, for $j=\pzt_n+2,\ldots, 2\pzt_n+1$, the $j$-th edge added will be between the current cluster of $v$ and a vertex that is not in $\tbar\cC^n_1(\gamma_n)$.
At step $2(\pzt_n+1)$, again check for new edges that connect the current cluster of $v$ to its complement (including edges that connect the cluster of $v$ to $\tbar\cC^n_1(\gamma_n)$), and at this step the edges found in step $(\pzt_n+1)$ that connect to $\tbar\cC^n_1(\gamma_n)$ will again be considered.
\item 
Continue while checking for possible new edges that connect the current cluster of $v$ to $\tbar\cC^n_1(\gamma_n)$ every $(\pzt_n+1)$ steps. 
Stop if we connect to $\tbar\cC^n_1(\gamma_n)$ (this is only possible at step $j(\pzt_n+1)$ for some $j\geq 1$), or if no new edge can be added to the current cluster.
\end{enumeratea}
The next lemma states a simple property relating the two algorithms.

\begin{lemma}\label{lem:simple}
For $v\notin\tbar\cC^n_1(\gamma_n)$ and any $i_0\geq 1$, if $v$ gets connected to $\tbar\cC^n_1(\gamma_n)$ after the addition of $i_0+1$ many edges in Algorithm 1, then Algorithm 1 and Algorithm 2 coincide up to the addition of $i_0$ edges. 
\end{lemma}

Define the event 
\begin{align}
\fB_n
:= \bigg\{
&
\cW\big(\tbar\cC^n_1(\gamma_n)\big) \geq (1-2\kappa_0)s^{\sss(n)}(\gamma_n)\cdot n^{\rho}\cdot\big(\sigma_2^{\sss(n)}\big)^{1/2}\, ,\notag\\
& 
w_n|\cC|\leq\cW(\cC) \leq C_{\ref{eqn:30}}(2)\cdot \log{n} \text{ for all other components }
\cC \text{ of } \barGn(\gamma_n)\, ,\text{ and} \notag\\
&
\text{degree of all } v\in [n]\setminus V\big(\tbar\cC^n_1(\gamma_n)\big) \text{ in }\barGn\text{ is } \leq (\log n)^{3/2} \bigg\}\, . \label{eqn:fbn-def}
\end{align}
Here, the first line corresponds to the event of interest in Proposition \ref{prop:lem-ss-3}, 
and the second line corresponds to the event in Lemma \ref{lem:ss-lem-12} with $\kappa=2$. 
By \eqref{it:lem:ss-8}, $\falln$,
\begin{equation}\label{eqn:43}
\pr\bigg(j\in \tbar\cC^n_1(\gamma_n)\ \text{ for }\ 1\leq j\leq \frac{\gamma_n^{1/\eta}}{(\log \gamma_n)^3}\bigg) 
\geq 
1-\exp\big(-C\log^{3\alpha}n\big)\, .  
\end{equation}
By Assumption \ref{ass:wts}\,(iii) and \eqref{eqn:up-bnd}, 
for any $j> \gamma_n^{1/\eta}(\log \gamma_n)^{-3}$, 
$w_j \leq 2^{\alpha}A_2(n/j)^{\alpha} \leq  C\log^{3\alpha}n$ for all large $n$. 
Now the degree of a vertex $v$ in $\barGn$ is distributed as 
$\sum_{j\neq v}\Bern\big(w_v w_j/\ell_n\big)$, 
where the summands are independent random variables.
Thus, using \eqref{eqn:43}, a union bound and Bennett's inequality \cite{boucheron2013concentration} shows that $\falln$,
\begin{equation}\label{eqn:312}
\pr\bigg(
\text{degree of some } v\in [n]\setminus V\big(\tbar\cC^n_1(\gamma_n)\big) \text{ in }\barGn\text{ is } > (\log n)^{3/2}
\bigg) 
\leq n\exp(-C\log^{3\alpha}n)\, .
\end{equation}
Combining \eqref{eqn:312}, Proposition \ref{prop:lem-ss-3}, and Lemma \ref{lem:ss-lem-12} yields, $\falln$,
\begin{equation}\label{eqn:44}
\pr\big(\fB_n^c\big) \leq 2/n^2\, . 
\end{equation}
Write
\[
\pr_2(\cdot) = \pr\bigg(\,\cdot\, \big|\, 
\big(\barGn(\gamma_n)\, ,\ \barGn \text{ restricted to } [n]\setminus V\big(\tbar\cC^n_1(\gamma_n)\big)\big)
\bigg)\, ,
\] 
and let $\E_2[\cdot]$ denote the corresponding expectation. 
Under $\pr_2$, for any realization of 
$\big(\barGn(\gamma_n)\, ,\ \barGn \text{ restricted to } [n]\setminus V\big(\tbar\cC^n_1(\gamma_n)\big)$,
the rest of the edges in $\barGn$ and the weights associated to the edges in
$E(\barGn)\setminus E(\barGn(\gamma_n))$
can be generated as follows:
Let 
\[
E_{\star}=\big\{\{i, j\}\, :\, 
1\leq i\neq j\leq n\, , 
\text{ at least one of } i\text{ or } j\text{ is in } V\big(\tbar\cC^n_1(\gamma_n)\big)\, ,
\text{ and }\{i, j\}\notin E\big(\tbar\cC^n_1(\gamma_n)\big)
\big\}\, .
\]
Let $\bar q_{ij}$ be as in \eqref{eqn:300}.
Place the edges $\{i, j\}\in E_{\star}$ independently with respective probabilities
\begin{align}\label{eqn:200}
\frac{\big(1-p^n_{\gamma_n}\big)\bar q_{ij}}{\big(1-p^n_{\gamma_n}\big)\bar q_{ij} + (1- \bar q_{ij})}
=
\big(1+O(n^{-\eta})\big)\cdot\big(1-p^n_{\gamma_n}\big)\frac{w_i w_j}{\ell_n}\, ,
\end{align}
where the $O(n^{-\eta})$ term is uniform over $\{i, j\}\in E_{\star}$, and $p^n_{\bullet}$ is as in \eqref{eqn:pnlam-def}.
Adding these edges would generate the complete set of edges $E(\barGn)$.
Now assign i.i.d. $\Unif\,[ p^n_{\gamma_n},\, 1]$ weights to the edges in 
$E(\barGn)\setminus E(\barGn(\gamma_n))$.

Fix a vertex $v\notin\tbar\cC^n_1(\gamma_n)$. 
Suppose Algorithm 2 started from $\tbar\cC^n_v(\gamma_n)$ terminates after the addition of $K$ edges;
let $j_0\geq 0$ be such that
$j_0\big(\pzt_n+1\big)\leq K< (j_0+1)\big(\pzt_n+1\big)$.
Let $T_j := (j\wedge j_0)\big(\pzt_n+1\big)$ for $j\geq 0$. 
If the path in $\tbar M^n$ connecting $v$ to $\tbar M^n_{\gamma_n}$ has length 
$\geq \big(\log n\big)^{-1/6}n^{\eta}$, 
then $\falln$, on the event $\fB_n$, Algorithm 1 started from $\tbar\cC^n_v(\gamma_n)$ will run for at least 
\[
\frac{n^{\eta}}{\big(\log n\big)^{1/6}}
\times
\frac{w_n}{w_n+C_{\ref{eqn:30}}(2)\log n}
\geq
\frac{C w_n n^{\eta}}{\big(\log n\big)^{7/6}}+1
\]
many steps.
Using Lemma \ref{lem:simple}, we see that Algorithm 2 started from $\tbar\cC^n_v(\gamma_n)$ will run for at least $Cw_n n^{\eta}\big(\log n\big)^{-7/6}$ many steps. 
Using \eqref{eqn:tn-def} and Assumption \ref{ass:wts}\,(iv), we see that in this case, for all large $n$, the events 
$\fA_j$, $j=1, 2,\ldots, n^{\eta/2}(\log n)^{7/12}$, will take place, where 
\[
\fA_j 
= \big\{
T_j =  j\big(\pzt_n+1\big) \text{ and the } T_j \text{-th edge in Algorithm 2 does not connect to } \tbar\cC^n_1(\gamma_n)
\big\}\, . 
\]
Both $T_j$ and $\fA_j$ depend on $v$, but we will suppress this dependence to simplify notation.
\begin{lemma}\label{lem:45}
For all large $n$, on the event $\fB_n$,
$\pr_2\big(\fA_{g(n)}\big) \leq 1/n^3$ 
for any $v\notin\tbar\cC^n_1(\gamma_n)$, where $g(n)=n^{\eta/2}(\log n)^{7/12}$.
\end{lemma}

\noindent{\bf Proof:}
The proof recursively analyzes $\fA_j$, $j\geq 1$. 
We will discuss how to analyze $\fA_1$ and $\fA_2$ in detail; 
the argument for a general $j$ is similar. 
Let $S_1$ be the sequence of the first $(T_1-1)\vee 0$ edges added to the cluster of $v$ in Algorithm 2 (arranged in the order they were added). 
To bound $\pr_2\big(\fA_1\big)$, we will actually prove an upper bound on 
$\pr_2\big(\fA_1\,\big|\, S_1 = {\cs_1}\big)$ that is uniform over the choice of $\cs_1$ with $\pr_2\big(S_1 = \cs_1\big)>0$ and $\text{length}(\cs_1) = \pzt_n$, where $\text{length}(\cs_1)$ denotes the number of edges in the sequence $\cs_1$.
Call such a choice of $\cs_1$ \emph{tenable}.

Since the status of the edges (presence or absence) connecting to $\tbar\cC^n_1(\gamma_n)$ are not checked in the first $\pzt_n$ steps of Algorithm 2, the probability of these edges being present under 
$\pr_2(\cdot| S_1 = \cs_1)$ is the same as in \eqref{eqn:200} for any tenable choice of edges $\cs_1$. 
Let $\cE_1\subseteq E(\barGn)$ be the set of edges found between the cluster of $v$ and $\tbar\cC^n_1(\gamma_n)$ in the $T_1$-th step in Algorithm 2. 
For any tenable $\cs_1$, 
\begin{equation*}
\E_2\big[ |\cE_1|\, \big|\, S_1 = \cs_1 \big]
= 
\big(1+O(n^{-\eta})\big)\cdot\big(1-p^n_{\gamma_n}\big)
\sum_{i\in \tbar\cC^n_1(\gamma_n)}\ 
\sum_{\small \substack{j \text{ in cluster of}\\ v \text{ at time } \pzt_n}} \frac{w_i w_j}{\ell_n}\, .
\end{equation*}
Now, the weight of the cluster of $v$ after the $\pzt_n$-th edge has been added in Algorithm 2 is at least 
$\pzt_n w_n$.
Further, for all large $n$, on the event $\fB_n$, $\cW\big(\tbar\cC^n_1(\gamma_n)\big) \geq Cn$ by virtue of \eqref{eqn:ss-8}, and consequently,
for any tenable $\cs_1$, 
\begin{align}\label{eqn:500}
\E_2\big[|\cE_1| \big| S_1 = \cs_1\big] \geq 2C_{\ref{eqn:500}}\pzt_n w_n\, . 
\end{align}
Letting $\fE_1 = \big\{|\cE_1|\geq C_{\ref{eqn:500}}\pzt_n w_n\big\}$, 
using Bennet's inequality \cite{boucheron2013concentration} and Assumption \ref{ass:wts}\,(iv),
we see that there exist $n_{\ref{eqn:46}}\geq 1$ and $C_{\ref{eqn:46}}>0$ depending only on $\pmtr$ such that for all $n\geq n_{\ref{eqn:46}}$, on the event $\fB_n$, for any tenable $\cs_1$, 
\begin{equation}\label{eqn:46}
\pr_2\big(\fE_1^c\, \big|\, S_1 = \cs_1\big) 
\leq 
\exp\big(-C_{\ref{eqn:46}}\big(\log n\big)^{7/6}\big)\, .
\end{equation}
Hence, for all $n\geq n_{\ref{eqn:46}}$, on the event $\fB_n$, for any tenable $\cs_1$,
\begin{align}
\pr_2\big(\fA_1\, \big|\, S_1 = \cs_1\big) 
&
\leq 
\pr_2\big(\fA_1\cap \fE_1\, \big|\, S_1 = \cs_1\big) 
+ 
\pr_2\big(\fE_1^c\, \big|\, S_1 = \cs_1\big)\notag\\
&
\leq 
\pr_2\big(\fA_1\, \big|\, \fE_1\cap\big\{S_1 = \cs_1\big\}\big) 
+
\exp\big(-C_{\ref{eqn:46}}\big(\log n\big)^{7/6}\big)\, . \label{eqn:508}
\end{align}
We will bound the first term on the right side of \eqref{eqn:508} with the help of Lemma \ref{lem:ss-lem17}\,(a).
To this end, note that on $\fE_1\cap\big\{S_1 = \cs_1\big\}$, 
in the $T_1$-th step, we have found $m_0 \geq C_{\ref{eqn:500}}\pzt_n w_n$ many edges connecting the cluster of $v$ to $\tbar\cC^n_1(\gamma_n)$.
The weights associated with these edges are i.i.d. $\Unif[x_0, 1]$ random variables, where $x_0=p^n_{\gamma_n}$.
Next, on $\fB_n$, the number of vertices in the cluster of $v$ after the addition of the $\pzt_n$-th edge is at most 
$C_{\ref{eqn:30}}(2)\cdot (1+\pzt_n)\log{n}/w_n$, and the degree of each of these vertices in $\barGn$ is at most $(\log n)^{3/2}$.
Hence, the number of outgoing edges from the cluster of $v$ at this point that do not connect to $\tbar\cC^n_1(\gamma_n)$ is
$k\leq C_{\ref{eqn:30}}(2)\cdot (1+\pzt_n)(\log{n})^{5/2}/w_n$.
Conditional on the values of the weights $U_{ij}$, $\{i, j\}\in\cs_1$, the weights of these $k$ edges are independent $\Unif\, [x_j, 1]$ random variables for some $x_j\geq x_0$, $j=1,\ldots, k$.
Thus, using Lemma \ref{lem:ss-lem17} (a), we get,
on the event $\fB_n$, for any tenable $\cs_1$,
\begin{equation}\label{eqn:509}
\pr_2\big(\fA_1\, \big|\, \fE_1\cap\big\{S_1 = \cs_1\big\}\big) 
\leq 
1
-\frac{C_{\ref{eqn:500}}\pzt_n w_n}{C_{\ref{eqn:500}}\pzt_n w_n 
+ 
C_{\ref{eqn:30}}(2)\cdot (1+\pzt_n)(\log{n})^{5/2}/w_n}\, .
\end{equation} 
Combining \eqref{eqn:508} and \eqref{eqn:509}, we see that there exists $n_{\ref{eqn:47}}\geq n_{\ref{eqn:46}}$ such that on the event $\fB_n$, 
\begin{align}\label{eqn:47}
&\pr_2\big(\fA_1\big)
\leq
\max_{\cs_1\text{ tenable}}\, 
\pr_2\big(\fA_1\, \big|\, S_1 = \cs_1\big) \\
&
\leq
1
-\bigg(
\frac{C_{\ref{eqn:500}}\pzt_n w_n}{C_{\ref{eqn:500}}\pzt_n w_n 
+ 
C_{\ref{eqn:30}}(2)\cdot (1+\pzt_n)(\log{n})^{5/2}/w_n}
\bigg)
+
\exp\big(-C_{\ref{eqn:46}}\big(\log n\big)^{7/6}\big)
\leq 
1-\frac{(\log{n})^{1/2}}{C_{\ref{eqn:47}} n^{\eta/2}}\, ,\notag
\end{align}
where the last step uses Assumption \ref{ass:wts}\,(iv),
and holds for $n\geq n_{\ref{eqn:47}}$.

Now let us move to $\fA_2$. 
We first note that $\pr_2\big(\fA_2\big) = \pr_2\big(\fA_2\cap \fA_1\big)$. 
Thus, 
\begin{equation}\label{eqn:48}
\pr_2\big(\fA_2\big) 
\leq 
\pr_2\big(\fA_2\, \big|\, \fA_1\cap \fE_1\big) \cdot \pr_2\big(\fA_1\big) 
+ 
\pr_2\big(\fE_1^c\big)\, .
\end{equation}
Let $S_2^{\sss (1)}$ be the sequence of the first $(T_2-1)\vee 0$ edges added in Algorithm 2 (arranged in the order they were added), 
and let $S_2^{\sss (2)}$ be the collection of edges found between $\tbar\cC^n_1(\gamma_n)$ and the cluster of $v$ in the $T_1$-th step in Algorithm 2.
Write $S_2=\big(S_2^{\sss (1)}, S_2^{\sss (2)} \big)$.
Call $\cs_2=\big(\cs_2^{\sss (1)}, \cs_2^{\sss (2)} \big)$ tenable if $\text{length}(\cs_2^{(1)})=2\pzt_n+1$,
$|\cs_2^{\sss (2)}|\geq C_{\ref{eqn:500}}\pzt_n w_n$, and 
$\pr\big(S_2 = \cs_2\big) >0$.
To bound $\pr_2\big(\fA_2\, \big|\, \fA_1\cap \fE_1\big)$ appearing on the right side of \eqref{eqn:48}, it is enough to prove a bound on $ \pr_2\big(\fA_2\, \big|\, S_2 = \cs_2\big)$ that is uniform over all tenable $\cs_2$.

Let $\cE_2$ be the collection of {\bf new} edges found between $\tbar\cC^n_1(\gamma_n)$ and the cluster of $v$ in the $T_2$-th step;
we emphasize that these edges were not present in $\cE_1$. 
Let $\fE_2 = \big\{|\cE_2|\geq C_{\ref{eqn:500}}\pzt_n w_n\big\}$. 
Then, for all $n\geq n_{\ref{eqn:46}}$, on the event $\fB_n$, for any tenable $\cs_2$
\begin{align}\label{eqn:724}
\pr_2\big(\fA_2\, \big|\, S_2 = \cs_2\big) 
&\leq 
\pr_2\big(\fA_2\cap \fE_2\, \big|\, S_2 = \cs_2\big) 
+ 
\pr_2\big(\fE_2^c\, \big|\, S_2 = \cs_2\big)\notag\\
&\leq
\pr_2\big(\fA_2\, \big|\, \fE_2\cap\{S_2 = \cs_2\}\big) 
+
\exp\big(-C_{\ref{eqn:46}}\big(\log n\big)^{7/6}\big)\, ,
\end{align}
where the last step follows from an argument similar to the one leading to \eqref{eqn:46}.
We will bound the first term on the right side of \eqref{eqn:724} with the help of Lemma \ref{lem:ss-lem17} (b).
To this end, note that on $\fE_2\cap\big\{S_2 = \cs_2\big\}$, 
in the $T_2$-th step, we have found $\tbar m_0 \geq C_{\ref{eqn:500}}\pzt_n w_n$ many new edges connecting the cluster of $v$ to $\tbar\cC^n_1(\gamma_n)$.
The weights associated with these edges are i.i.d. $\Unif\, [\tbar x_0, 1]$ random variables, where 
$\tbar x_0=p^n_{\gamma_n}$.
By an argument similar to the one used while analyzing $\fA_1$, on the event $\fB_n$, 
the number of outgoing edges from the cluster of $v$ that do not connect to $\tbar\cC^n_1(\gamma_n)$ and were found after the $T_1$-th step is
$\tbar k\leq C_{\ref{eqn:30}}(2)\cdot (1+\pzt_n)(\log{n})^{5/2}/w_n$.
Conditional on the values of the weights $U_{ij}$, $\{i, j\}\in\cs_2^{\sss (1)}$, the weights of these $\tbar k$ edges are independent $\Unif\, [\tbar x_l, 1]$ random variables for some $\tbar x_l\geq \tbar x_0$,  $l=1,\ldots, \tbar k$.
Note also that conditional on the values of the weights $U_{ij}$, $\{i, j\}\in\cs_2^{\sss (1)}$, the weights associated with the $m_0\geq C_{\ref{eqn:500}}\pzt_n w_n$ edges found between $\tbar\cC^n_1(\gamma_n)$ and the cluster of $v$ in the $T_1$-th step are i.i.d. $\Unif\, [x'_0, 1]$ random variables, where $x'_0$ is the weight of the $(\pzt_n+1)$-th edge in $\cs_2^{\sss (1)}$, and 
the weights associated with the $k'$ (say) outgoing edges from the cluster of $v$ that do not connect to $\tbar\cC^n_1(\gamma_n)$ and were found in the first $T_1$ steps are independent $\Unif\, [x_j', 1]$ random variables for some $x_j'\geq x_0'$, $j=1,\ldots, k'$.
Thus, using Lemma \ref{lem:ss-lem17} (b), we get,
on the event $\fB_n$, for any tenable $\cs_2$,
\[
\pr_2\big(\fA_2\, \big|\, \fE_2\cap\big\{S_2 = \cs_2\big\}\big) 
\leq 
1
-\frac{C_{\ref{eqn:500}}\pzt_n w_n}{C_{\ref{eqn:500}}\pzt_n w_n 
+ 
C_{\ref{eqn:30}}(2)\cdot (1+\pzt_n)(\log{n})^{5/2}/w_n}\, ,
\]
which combined with \eqref{eqn:724} shows that for $n\geq n_{\ref{eqn:47}}$, on the event $\fB_n$,
\begin{align}\label{eqn:5099}
\pr_2\big(\fA_2\, \big|\, \fA_1\cap \fE_1\big)
\leq
\max_{\cs_2\text{ tenable}}\,
\pr_2\big(\fA_2\, \big|\, S_2 = \cs_2\big) 
\leq 
1-\frac{(\log{n})^{1/2}}{C_{\ref{eqn:47}} n^{\eta/2}}\, .
\end{align}
Now going back to \eqref{eqn:48} and using \eqref{eqn:46}, \eqref{eqn:47}, and \eqref{eqn:5099}, we see that for $n\geq n_{\ref{eqn:47}}$, on the event $\fB_n$,
\[
\pr_2\big(\fA_2\big)
\leq 
\bigg(1-\frac{(\log{n})^{1/2}}{C_{\ref{eqn:47}} n^{\eta/2}}\bigg)^2
+ 
\exp\big(-C_{\ref{eqn:46}}\big(\log n\big)^{7/6}\big)\, .
\]

Turning to the events $\fA_j$ for $3\leq j\leq g(n)$,
we define $\fE_j$ in a manner analogous to $\fE_1$ and $\fE_2$.
Proceeding similarly, we can show that for $n\geq n_{\ref{eqn:47}}$ and $2\leq j\leq g(n)$, on the event $\fB_n$,
\begin{align*}
\pr_2\big(\fA_{j+1}\big) 
&
=
\pr_2\bigg(\bigcap_{i=1}^{j+1}\fA_i\bigg)
\leq 
\pr_2\bigg(\fA_{j+1}\, \bigg|\, \bigcap_{i=1}^{j}\big(\fA_i\cap \fE_i\big)\bigg) \cdot \pr_2\big(\fA_j\big) 
+ 
\sum_{i=1}^j\pr_2\big(\fE_i^c\big)\, \\
&
\leq
\bigg(1-\frac{(\log{n})^{1/2}}{C_{\ref{eqn:47}} n^{\eta/2}}\bigg)\cdot \pr_2\big(\fA_j\big) 
+
j\exp\big(-C_{\ref{eqn:46}}\big(\log n\big)^{7/6}\big)\\
&
\leq
\bigg(1-\frac{(\log{n})^{1/2}}{C_{\ref{eqn:47}} n^{\eta/2}}\bigg)^{j+1}
+
\frac{j(j+1)}{2} \exp\big(-C_{\ref{eqn:46}}\big(\log n\big)^{7/6}\big)\, .
\end{align*}
For the third step, we need to use the analogue of Lemma \ref{lem:ss-lem17} for $(j+1)$ collections of independent uniform random variables.
However, this generalization is straightforward.
This completes the proof of Lemma \ref{lem:45}.
\qed

\medskip

\noindent{\bf Completing the Proof of Proposition \ref{prop:ss-11}:} 
Combining \eqref{eqn:44} with Lemma \ref{lem:45} and using a union bound over 
$v\notin\tbar\cC^n_1(\gamma_n)$, we get, for all large $n$,
\[
\pr\bigg(
d_{\rH}\big(\tbar M^n_{\gamma_n}\, ,\, \tbar M^n\big) 
\geq 
\big(\log n\big)^{-1/6}n^{\eta}  
\bigg) 
\leq 3/n^2\, .
\]
This completes the proof. 
\qed

\section{Proof of Theorem \ref{thm:mst-convg}}\label{sec:proof-main-thm}
We will complete the proof in five steps. 
As observed in Section \ref{sec:different-model}, it is enough prove the result for $M^n$ and $\tbar M^n$.

\subsection{Existence of the scaling limit}\label{sec:existence-scaling-limit}
Our aim in this section is to show that there exists a random compact $\bR$-tree $\mst^{\mvtheta^*}$ whose law depends only on $\mvtheta^*$ such that \eqref{eqn:666} is satisfied if we replace $\mvM^n$ by $M^n$ or $\tbar M^n$.

Note that \eqref{eqn:ss-2}, the first relation in \eqref{eqn:ss-1}, Assumption \ref{ass:wts}\,(iii), and \eqref{eqn:up-bnd} imply that $\mvx^{\sss (n)}$, $n\geq 1$, satisfies \cite[Assumption 1.6]{SB-vdH-SS-PTRF} with the limiting sequence $\mvtheta^*$.
Now consider the random graph $\cG_n(\lambda)$ as in \eqref{eqn:712}.
Then \cite[Theorem 1.8]{SB-vdH-SS-PTRF} shows that the largest components (arranged according to their masses) of $\cG_n(\lambda)$, endowed with the graph distance scaled by $n^{-\eta}$ and the probability measure that assigns mass proportional to $x_v$ to each vertex $v$ in a component, converges in distribution to a sequence of random metric measure spaces with respect to the product topology induced by the Gromov-weak topology on each coordinate.
We claim that in a similar manner, we can show that for every $\lambda\in\bR$ there exists a random metric measure space $\pzS_{\lambda}^{\mvtheta^*}$ whose law depends only on $\lambda$ and $\mvtheta^*$ such that as $n\to\infty$,
\begin{align}\label{eqn:63}
n^{-\eta}\cdot\cC\big(1, \cG_n(\lambda)\big)
\weakc
\pzS_{\lambda}^{\mvtheta^*}\ \ \text{ w.r.t. the Gromov-weak topology}\, ,
\end{align}
where the measure on 
$n^{-\eta}\cdot\cC\big(1, \cG_n(\lambda)\big)$ 
assigns mass $x_v\big[\cW_{\mvx}\big(\cC\big(1, \cG_n(\lambda)\big)\big)\big]^{-1}$ to each vertex $v$ in $\cC\big(1, \cG_n(\lambda)\big)$.
We briefly explain how this can be done.
The proof of \cite[Theorem~1.8]{SB-vdH-SS-PTRF}, applied to the special case $\cG_n(\lambda)$, can be divided into two steps:
\begin{inparaenumi}
\item
In \cite[Theorem 4.5]{SB-vdH-SS-PTRF}, it is proved that under some assumptions, a related connected random graph has a scaling limit in the Gromov-weak topology.
\item\label{it:222} 
Then the arguments of \cite[Section 5.1]{SB-vdH-SS-PTRF} show that the maximal connected components of $\cG_n(\lambda)$ satisfy the assumptions in \cite[Theorem~4.5]{SB-vdH-SS-PTRF}.
\end{inparaenumi}
A key result needed to complete step \eqref{it:222} is \cite[Proposition~9]{aldous-limic} (restated in \cite[Proposition~5.2]{SB-vdH-SS-PTRF}), which shows that a breadth-first walk of $\cG_n(\lambda)$ defined in \cite{aldous-limic} converges in distribution w.r.t. the Skorohod $J_1$ topology to a limiting process.
In a similar way (see also the proof of \cite[Theorem~2.4]{SBVHVJL12}), one can prove that for any $\lambda\geq 0$, the breadth-first walk of $\cG_n(\lambda)$ started from the vertex $1$ as given by \eqref{eqn:8888} satisfies
\begin{align}\label{eqn:777}
n^{\eta}\cdot Z_{\lambda}^{\sss n, (1)}(u)
\weakc
Z_{\lambda}^{\sss (1)}(u)
:=
\theta_1^*+\lambda u+\sum_{j=2}^{\infty}\theta_j^*\big(\ind_{\{\xi_j\leq u\}} -u\theta_j^*\big)
\, ,\ \ \ u\geq 0\, 
\end{align}
with respect to the Skorohod $J_1$ topology on $D[0, \infty)$, where $\xi_j\sim \EXP(\theta_j^*)$, $j\geq 1$, are independent random variables.
Now using \eqref{eqn:777}, the rest of the arguments in \cite[Section 5.1]{SB-vdH-SS-PTRF} repeated verbatim would show that the vertex weights in $\cC\big(1, \cG_n(\lambda)\big)$ satisfy the assumptions of \cite[Theorem 4.5]{SB-vdH-SS-PTRF}, which would then yield \eqref{eqn:63}.
For concreteness, we describe the construction of $\pzS_{\lambda}^{\mvtheta^*}$.

\vskip3pt

\begin{construction}\label{construction:S-lambda}
Define 
$
\cZ_{\lambda}:=\inf\big\{u\geq 0:  Z_{\lambda}^{\sss (1)}(u)\leq 0\big\}
$,
and
$
\cB_{\lambda}:=\{1\}\cup\big\{j\geq 2\, :\, \xi_j\leq\cZ_{\lambda} \big\}
$.
Let 
\[
\tbar\gamma:=
\cZ_{\lambda}\cdot\bigg(\sum_{j\in\cB_{\lambda}}\big(\theta_j^*\big)^2\bigg)^{1/2}\, ,\ \ 
\text{ and }\ \
\tbar\mvtheta:=
\bigg(
\theta_j^*\cdot\bigg(\sum_{i\in\cB_{\lambda}}\big(\theta_i^*\big)^2\bigg)^{-1/2}\, ;\ j\in\cB_{\lambda}
\bigg)\, .
\]
Set
\[
\pzS_{\lambda}^{\mvtheta^*}
=
\cZ_{\lambda}\cdot\bigg(\sum_{i\in\cB_{\lambda}}\big(\theta_i^*\big)^2\bigg)^{-1/2}
\cdot\cG_{\infty}\big(\tbar\mvtheta,\, \tbar\gamma \big)\, ,
\]
where $\cG_{\infty}(\cdot\, ,\, \cdot)$ is as given in \cite[Definition 2.2]{SB-vdH-SS-PTRF}.
\end{construction}

Next, we claim that under Assumption \ref{ass:wts}, $\pzS_{\lambda}^{\mvtheta^*}$ is compact almost surely, and further, as $n\to\infty$,
\begin{align}\label{eqn:64}
n^{-\eta}\cdot\cC\big(1, \cG_n(\lambda)\big)
\weakc
\pzS_{\lambda}^{\mvtheta^*}\ \ \text{ w.r.t. the GHP topology.}
\end{align}
To deduce \eqref{eqn:64} from \eqref{eqn:63}, we need to prove that 
$n^{-\eta}\cdot\cC\big(1, \cG_n(\lambda)\big)$
satisfies a global lower mass-bound; 
see \cite[Theorem 6.1]{AthLorWin14} and \cite[Lemma 6.1]{SB-vdH-SS-PTRF}.
The global lower mass-bound was established in \cite{SB-vdH-SS-PTRF} under the stronger \cite[Assumption 1.1]{SB-vdH-SS-PTRF}.
This stronger assumption was needed to prove \cite[Lemma 6.7]{SB-vdH-SS-PTRF}, which relied on a tail bound on the heights of branching processes established in \cite[Theorem 2]{Kortchemski}.
The result in \cite[Theorem 2]{Kortchemski} was proved for a (conditioned) branching process with a given offspring distribution, whereas in the random graph setting, one needs to consider branching processes with varying offspring distributions for different $n$; e.g., the offspring distributions of interest in this paper are $\Poi(V_n)$, $n\geq 1$, where $V_n$ is as defined around \eqref{def:V-n}.
To get around this difficulty, \cite[Assumption 1.1]{SB-vdH-SS-PTRF} was used in the proof of \cite[Lemma 6.6]{SB-vdH-SS-PTRF} to stochastically upper bound the offspring distributions for different $n$ by a single offspring distribution to which \cite[Theorem 2]{Kortchemski} is applicable.
However, instead of \cite[Lemma 6.7]{SB-vdH-SS-PTRF}, we can now appeal to Proposition \ref{prop:lem-6-bp-height} (which relies on the techniques developed in the more recent work \cite{addario2019most}), and avoid the use of the stronger \cite[Assumption 1.1]{SB-vdH-SS-PTRF}.
Then the rest of the argument from \cite{SB-vdH-SS-PTRF} carries over in an identical way under Assumption \ref{ass:wts} to establish the desired global lower mass-bound for $n^{-\eta}\cdot\cC\big(1, \cG_n(\lambda)\big)$, which then yields \eqref{eqn:64}.
Let us also note that the construction of $\pzS_{\lambda}^{\mvtheta^*}$ given above together with the almost sure compactness of $\pzS_{\lambda}^{\mvtheta^*}$ shows that $\pzS_{\lambda}^{\mvtheta^*}$ is an $\bR$-graph almost surely.
Further, \eqref{eqn:712} and \eqref{eqn:64} imply
\begin{align}\label{eqn:64-a}
n^{-\eta}\cdot\cC_1^n(\lambda)
\weakc
\pzS_{\lambda}^{\mvtheta^*}\ \ \text{ w.r.t. the GH topology.}
\end{align}

The next result shows that Theorem \ref{thm:mst-from-ghp-convergence} can be applied to the sequence $n^{-\eta}\cdot\cC^n_1(\lambda)$, $n\geq 1$.
Let $\cA_r$ be as defined around \eqref{eqn:65}.

\begin{lem}\label{lem:22}
	Fix $\lambda\geq 0$.
	Then we can construct $\pzS_{\lambda}^{\mvtheta^*}$, 
	$n^{-\eta}\cdot\cC^n_1(\lambda)$, $n\geq 1$, and
	a positive random variable $R$ on the same probability space such that
	$n^{-\eta}\cdot\cC^n_1(\lambda)
	\convas
	\pzS_{\lambda}^{\mvtheta^*}$ w.r.t. the GH topology as $n\to\infty$, and 
	$~\pr\big(n^{-\eta}\cdot\cC^n_1(\lambda)\in\cA_R\ \text{ eventually}\,\big)=1$.
\end{lem}

\begin{rem}
We will omit the proof of Lemma~\ref{lem:22}, as this result is essentially contained in \cite{SB-vdH-SS-PTRF}.
We only make a brief comment on how this result follows.
Recall the notation $\big(k(X), e(X)\big)$ for an $\bR$-graph $X$ from Section~\ref{sec:r-tree} 
\ch{and the probability distribution $\bP_{\con}$ given by \eqref{eqn:pr-con-vp-a-cV-def}.
For $m\geq 1$, $a>0$, and a probability mass function $\vp=(p_1,\ldots, p_m)$ on $[m]$ with 
$p_1\geq p_2\geq\cdots\geq p_m>0$, 
write $\tilde\cG_m(\vp, a)$ for a random graph with distribution $\bP_{\con}(\cdot, \vp, a, [m])$.
Let $\cG_{\infty}(\cdot\, ,\, \cdot)$ be as given in \cite[Definition~2.2]{SB-vdH-SS-PTRF}.}
Now, as mentioned at the beginning of this section, the key ingredient in the proof of \eqref{eqn:63} is \cite[Theorem 4.5]{SB-vdH-SS-PTRF}, which shows that under \cite[Assumption 4.4]{SB-vdH-SS-PTRF}, 
\begin{align}\label{eqn:709}
\|\vp\|_2\cdot\tilde\cG_m(\vp, a)\weakc\cG_{\infty}\big(\mvtheta, \gamma\big)
\end{align}
with respect to Gromov-weak topology as $m\to\infty$.
It also follows from the proof of \cite[Lemma 4.12]{SB-vdH-SS-PTRF} and \cite[(4.23)]{SB-vdH-SS-PTRF} that under \cite[Assumption 4.4]{SB-vdH-SS-PTRF},
\begin{align}\label{eqn:710}
\spls\big(\tilde\cG_m(\vp, a)\big)\weakc\spls\big(\cG_{\infty}(\mvtheta, \gamma)\big)
\end{align}
jointly with the convergence in \eqref{eqn:709}.
Now, from the proof of \cite[Theorem 4.5]{SB-vdH-SS-PTRF}, it can be easily deduced that 
\[
\|\vp\|_2\cdot\bigg(\len\big(\core\big(\tilde\cG_m(\vp, a)\big)\big), \, 
\min_{e\in e(\tilde\cG_m(\vp, a))}\len(e)
\bigg)
\weakc
\bigg(\len\big(\core\big(\cG_{\infty}(\mvtheta, \gamma)\big), \, 
\min_{e\in e(\cG_{\infty}(\mvtheta, \gamma))}\len(e)
\bigg)
\]
jointly with the convergences in \eqref{eqn:709} and \eqref{eqn:710}.
These results directly translate into the corresponding convergences for $\cC_1^n(\lambda)$, which then yield the claim in Lemma \ref{lem:22}.
\end{rem}

\ch{We will now deduce the existence of the MST scaling limit from Lemma~\ref{lem:22} and Theorem~\ref{thm:mst-crit-gh}.
As mentioned in Section~\ref{sec:proof-strategy}, the argument for combining these two results to get the GH scaling limit of $\tbar M^n$ is similar to the one used in \cite{AddBroGolMie13}.}
Using Lemma~\ref{lem:22} and Theorem~\ref{thm:mst-from-ghp-convergence}, we get
\begin{align}\label{eqn:75-a}
n^{-\eta}\cdot\cbd^{\infty}\big(\cC^n_1(\lambda)\big)
\weakc
\cb^{\infty}\big(\pzS_{\lambda}^{\mvtheta^*}\big)
=:\mst^{\mvtheta^*}_{\lambda}
\ \ \text{ w.r.t. the GH topology}
\end{align}
as $n\to\infty$.
Similar to Lemma \ref{lem:minus-coup}, for any $\lambda\geq 0$, there exist couplings of $G_n(\lambda)$ and $\barGn(\lambda)$ such that 
$\pr\big(G_n(\lambda)\neq\barGn(\lambda)\big)\to 0$ as $n\to\infty$.
This fact, together with \eqref{eqn:75-a} and Lemma \ref{lem:cycle-breaking-gives-mst}, gives, for any $\lambda\geq 0$,
\begin{align}\label{eqn:75}
n^{-\eta}\cdot\tbar M^n_{\lambda}
\weakc
\mst^{\mvtheta^*}_{\lambda}
\ \ \text{ w.r.t. the GH topology}
\end{align}
as $n\to\infty$.
As discussed below Definition \ref{def:graph-vgnp}, 
$\tbar M_{\lambda}^n$ is a subtree of $\tbar M_{\lambda'}^n$ whenever $0\leq\lambda\leq\lambda'$.
For $0\leq\lambda\leq\lambda'$, consider the sequence of marked spaces 
$\big(n^{-\eta}\cdot\tbar M^n_{\lambda'}\, ,\, n^{-\eta}\cdot\tbar M^n_{\lambda}\big)$, $n\geq 1$.
Since $\mst^{\mvtheta^*}_{\lambda'}$ is compact almost surely,
\cite[Proposition 9]{miermont2009tessellations} implies that this \ch{sequence} is tight in the marked GH topology.
Thus, by passing to a subsequence, we can assume that
\begin{align}\label{eqn:51}
\big(n^{-\eta}\cdot\tbar M^n_{\lambda'}~ ,\, n^{-\eta}\cdot\tbar M^n_{\lambda}\big)
\weakc
\big(\mst^{\mvtheta^*}_{\lambda'},\, \mst^{\mvtheta^*}_{\lambda}\big)
\end{align}
with respect to the marked GH topology.
On the right side of \eqref{eqn:51}, we have a coupling of $\mst^{\mvtheta^*}_{\lambda'}$ and $\mst^{\mvtheta^*}_{\lambda}$ in which the latter is a subspace of the former.
Since 
$
d_{\rH}\big(\tbar M^n_{\lambda}, \tbar M^n_{\lambda'}\big)
\leq
d_{\rH}\big(\tbar M^n_{\lambda}, \tbar M^n\big)
$,
an application of Theorem \ref{thm:mst-crit-gh} shows that for all $\Delta\in(0, 1/2]$,
in the coupling of \eqref{eqn:51},
\begin{align}\label{eqn:92}
\pr\big(
d_{\rH}\big(\mst^{\mvtheta^*}_{\lambda},\, \mst^{\mvtheta^*}_{\lambda'}\big)
>
\lambda^{-1+\Delta}
\big)
\leq
C_{\ref{eqn:555}}(\Delta)\cdot\lambda^{-1/2}
\end{align}
for $\lambda'\geq\lambda\geq\lambda_{\ref{eqn:555}}(\Delta)$. 
This implies that for $\lambda'\geq\lambda\geq\lambda_{\ref{eqn:555}}(\Delta)$,
\begin{align}\label{eqn:52}
d_{\rP\rr}\bigg(
\mathpzc{Law}\big(\mst^{\mvtheta^*}_{\lambda}\big),\, 
\mathpzc{Law}\big(\mst^{\mvtheta^*}_{\lambda'}\big)
\bigg)
\leq
\lambda^{-1+\Delta}+C_{\ref{eqn:555}}(\Delta)\cdot\lambda^{-1/2}\, ,
\end{align}
where $\mathpzc{Law}(\cdot)$ denotes the law of a random object, and 
$d_{\rP\rr}(\cdot\, ,\,\cdot)$ denotes the Prokhorov distance.
Thus, 
$\big(\mathpzc{Law}(\mst^{\mvtheta^*}_{\lambda})\, ,\ \lambda\geq 0\big)$ is a Cauchy sequence in the space of probability measures on $\fS_{\GH}$.
Since $\fS_{\GH}$ is Polish, the space of probability measures on $\fS_{\GH}$ is complete under  $d_{\rP\rr}(\cdot\, ,\,\cdot)$.
Hence, there exists a random compact metric space $\mst^{\mvtheta^*}$ whose law depends only on $\mvtheta^*$ such that
$
\mst^{\mvtheta^*}_{\lambda}\weakc\mst^{\mvtheta^*}
$ 
w.r.t. the GH topology as $\lambda\to\infty$.
Since $\mst^{\mvtheta^*}_{\lambda}$ is an $\bR$-tree for any $\lambda$, \cite[Lemma 2.1]{evans-pitman-winter} shows that $\mst^{\mvtheta^*}$ is an $\bR$-tree almost surely. 
From \eqref{eqn:52}, we further get
\begin{align}\label{eqn:52-a}
d_{\rP\rr}\bigg(
\mathpzc{Law}\big(\mst^{\mvtheta^*}_{\lambda}\big),\, 
\mathpzc{Law}\big(\mst^{\mvtheta^*}\big)
\bigg)
\leq
\lambda^{-1+\Delta}+C_{\ref{eqn:555}}(\Delta)\cdot\lambda^{-1/2}
\end{align}
for $\lambda\geq\lambda_{\ref{eqn:555}}(\Delta)$.
Now, Theorem \ref{thm:mst-crit-gh} and \eqref{eqn:52-a} coupled with an application of the triangle inequality yield
\begin{align*}
&\limsup_{n\to\infty}\, d_{\rP\rr}\bigg(
\mathpzc{Law}\big(n^{-\eta}\cdot\tbar M^n\big),\, 
\mathpzc{Law}\big(\mst^{\mvtheta^*}\big)
\bigg)\\
&\hskip50pt
\leq
2\big(\lambda^{-1+\Delta}+C_{\ref{eqn:555}}(\Delta)\cdot\lambda^{-1/2}\big)
+
\limsup_{n\to\infty}\, d_{\rP\rr}\bigg(
\mathpzc{Law}\big(n^{-\eta}\cdot\tbar M^n_{\lambda}\big)\, ,\, 
\mathpzc{Law}\big(\mst^{\mvtheta^*}_{\lambda}\big)
\bigg)
\end{align*}
for any $\lambda\geq\lambda_{\ref{eqn:555}}(\Delta)$.
Now using \eqref{eqn:75} and letting $\lambda\to\infty$ gives 
\begin{align}\label{eqn:9999}
n^{-\eta}\cdot\tbar M^n\weakc\mst^{\mvtheta^*}\ \ \text{ w.r.t. the GH topology}
\end{align}
as $n\to\infty$.
Finally, Lemma \ref{lem:minus-coup} shows that \eqref{eqn:9999} continues to hold if we replace 
$\tbar M^n$ by $M^n$.

\subsection{The degrees of points in $\mst^{\mvtheta^*}$}
In this section we will prove Theorem \ref{thm:mst-convg}\,\eqref{it:a}.
Recall the definition of degree from Section \ref{sec:r-tree} and the notation $\cL(\cdot)$ and $\cH(\cdot)$ from \eqref{eqn:leaf-hub-def}.
From Construction \ref{construction:S-lambda} and the construction of the space 
$\cG_{\infty}(\cdot\, ,\, \cdot)$ using an inhomogeneous continuum random tree (ICRT) given in \cite[Section 2.3.1]{SB-vdH-SS-PTRF}, it follows that 
\ch{the set of points in $\pzS_{\lambda}^{\mvtheta^*}$ with infinite degree is countably infinite}, 
and all other points in $\pzS_{\lambda}^{\mvtheta^*}$ either have degree $1$ or $2$.
Since 
$
\mst^{\mvtheta^*}_{\lambda}=\cb^{\infty}\big(\pzS_{\lambda}^{\mvtheta^*}\big)
$,
the same is true of $\mst^{\mvtheta^*}_{\lambda}$, i.e.,
\ch{the set $\cH(\mst^{\mvtheta^*}_{\lambda})$ is countably infinite} and 
$\deg\big(x\, ;\, \mst^{\mvtheta^*}_{\lambda}\big)\in\big\{1, 2 \big\}$ for all 
$x\in\mst^{\mvtheta^*}_{\lambda}\setminus\cH(\mst^{\mvtheta^*}_{\lambda})$.

By an argument similar to the one leading to \eqref{eqn:51}, we can assume that for any $\lambda\geq 1$,
\begin{align}\label{eqn:53}
\big(n^{-\eta}\cdot\tbar M^n~ ,\, n^{-\eta}\cdot\tbar M^n_{\lambda}\big)
\weakc
\big(\mst^{\mvtheta^*},\, \mst^{\mvtheta^*}_{\lambda}\big)
\end{align}
with respect to the marked GH topology along a suitable subsequence (which may depend on $\lambda$).
In the coupling between $\mst^{\mvtheta^*}_{\lambda}$ and $\mst^{\mvtheta^*}$ on the right side of \eqref{eqn:53}, the former is a closed subset of the latter.
Hence, $\mst^{\mvtheta^*}$ has infinitely many points of infinite degree almost surely.

Now, Theorem \ref{thm:mst-crit-gh} and \eqref{eqn:53} imply that for all $\lambda\geq 1$, 
\begin{align}\label{eqn:ss-53a}
\pr\big(
d_{\rH}\big(\mst^{\mvtheta^*}_{\lambda},\, \mst^{\mvtheta^*}\big)
>
\lambda^{-1/2}
\big)
\leq
C_{\ref{eqn:ss-53a}}\lambda^{-1/2}\, .
\end{align}
in the coupling of \eqref{eqn:53}. 
In this coupling, on the event
\begin{align}\label{eqn:ss-53b}
\fD_{\lambda}
:=
\big\{ 
d_{\rH}\big(\mst^{\mvtheta^*}_{\lambda},\, \mst^{\mvtheta^*}\big)
\leq
\lambda^{-1/2}
\big\}\, ,
\end{align}
$\mst^{\mvtheta^*}$ can be obtained by attaching countably many $\bR$-trees each having diameter at most $2\lambda^{-1/2}$ to $\mst^{\mvtheta^*}_{\lambda}$.
Hence, on $\fD_{\lambda}$, any $x\in\mst^{\mvtheta^*}$ that satisfies
\begin{enumerateA}
\item\label{it:1}
$3\leq \deg\big(x\, ;\, \mst^{\mvtheta^*}\big)<\infty$, and
\item[]\label{it:2}\hskip-20pt(B$_{\lambda}$)
all $\deg\big(x\, ;\, \mst^{\mvtheta^*}\big)$ of the components in $\mst^{\mvtheta^*}\setminus\{x\}$ have diameter $>2\lambda^{-1/2}$
\end{enumerateA}
must also satisfy $x\in\mst^{\mvtheta^*}_{\lambda}$, and each of the 
$\deg\big(x\, ;\, \mst^{\mvtheta^*}\big)$ components in $\mst^{\mvtheta^*}\setminus\{x\}$ must have a non-empty intersection with $\mst^{\mvtheta^*}_{\lambda}$.
But this implies that 
$\deg\big(x\, ;\, \mst^{\mvtheta^*}_{\lambda}\big)\geq 3$, which in turn implies that 
$\deg\big(x\, ;\, \mst^{\mvtheta^*}_{\lambda}\big)=\infty$, and consequently,
$\deg\big(x\, ;\, \mst^{\mvtheta^*}\big)=\infty$ -- a contradiction.
Hence, on the event $\fD_{\lambda}$, there does not exist $x\in\mst^{\mvtheta^*}$ satisfying \eqref{it:1} and (B$_{\lambda}$) as above.
This observation combined with \eqref{eqn:ss-53a} shows that for all $\lambda\geq 1$,
\begin{align}
\pr\bigg(
\pzZ^{\mvtheta^*}>2\lambda^{-1/2}
\bigg)
\leq
C_{\ref{eqn:ss-53a}}\lambda^{-1/2}\, ,\label{eqn:800}
\end{align}
where
\[
\pzZ^{\mvtheta^*}
:=
\sup\bigg\{
\min\big\{\diam(\cT)\, :\, \cT\text{ component of }\mst^{\mvtheta^*}\setminus\{x\} \big\}
 \, :\, x\in\mst^{\mvtheta^*},\ 3\leq\deg\big(x; \mst^{\mvtheta^*}\big)<\infty \bigg\}\, .
\]
Here supremum of an empty set is defined as zero.
Then $\pzZ^{\mvtheta^*}=0$ if and only if there does not exist 
$x\in\mst^{\mvtheta^*}$ with $3\leq\deg\big(x; \mst^{\mvtheta^*}\big)<\infty$.
Now, letting $\lambda\to\infty$ in \eqref{eqn:800} shows that $\pzZ^{\mvtheta^*}=0$ almost surely, which concludes the proof of Theorem \ref{thm:mst-convg}\,\eqref{it:a}.

\subsection{Leaves and hubs of $\mst^{\mvtheta^*}$}
In this section we will prove Theorem \ref{thm:mst-convg}\,\eqref{it:b}.
Consider $\tbar\mvtheta$ as in Construction \ref{construction:S-lambda}.
Then in the ICRT corresponding to $\tbar\mvtheta$, both the set of leaves and the set of hubs are dense almost surely. 
(We refer the reader to \cite{pitman-camarri, aldous-pitman-entrance} for background on ICRTs.)
From the construction of the space $\cG_{\infty}(\cdot,\, \cdot)$ given in \cite[Section 2.3.1]{SB-vdH-SS-PTRF}, the same is true for $\cG_{\infty}\big(\tbar\mvtheta, \tbar\gamma\big)$, and consequently for $\pzS_{\lambda}^{\mvtheta^*}$.
Since
$\mst^{\mvtheta^*}_{\lambda}=\cb^{\infty}\big(\pzS_{\lambda}^{\mvtheta^*}\big)$,
both $\cL\big(\mst^{\mvtheta^*}_{\lambda}\big)$ and $\cH\big(\mst^{\mvtheta^*}_{\lambda}\big)$ are dense in $\mst^{\mvtheta^*}_{\lambda}$.

Fix $k\geq 1$, and consider $\lambda\geq 1$ large so that $3\lambda^{-1/2}<1/k$.
Consider the coupling between $\mst^{\mvtheta^*}_{\lambda}$ and $\mst^{\mvtheta^*}$ as in \eqref{eqn:53}.
On the event $\fD_{\lambda}$ defined in \eqref{eqn:ss-53b}, any point in $\mst^{\mvtheta^*}$ is within distance $\lambda^{-1/2}$ from a point in $\mst^{\mvtheta^*}_{\lambda}$, and consequently within distance $2\lambda^{-1/2}$ from a point in $\cH\big(\mst^{\mvtheta^*}_{\lambda}\big)$.
Since $\cH\big(\mst^{\mvtheta^*}_{\lambda}\big)\subseteq\cH\big(\mst^{\mvtheta^*}\big)$,
using \eqref{eqn:ss-53a}, we see that 
\[
\pr\bigg(
\mst^{\mvtheta^*}
=
\bigcup_{x\in\cH(\mst^{\mvtheta^*})} B\big(x, k^{-1}\, ;\, \mst^{\mvtheta^*}\big)
\bigg)
\geq
1-C_{\ref{eqn:ss-53a}}\lambda^{-1/2}\, .
\]
Letting $\lambda\to\infty$, we get
\[
\pr\bigg(\forall k\geq 1\, ,\ \ 
\mst^{\mvtheta^*}
=
\bigcup_{x\in\cH(\mst^{\mvtheta^*})} B\big(x, k^{-1}\, ;\, \mst^{\mvtheta^*}\big)
\bigg)
=1\, ,
\]
which shows that $\cH(\mst^{\mvtheta^*})$ is dense in $\mst^{\mvtheta^*}$ almost surely.

We now turn to the leaves of $\mst^{\mvtheta^*}$.
Once again, fix $k\geq 1$, and consider $\lambda\geq 1$ large so that $3\lambda^{-1/2}<1/k$.
Arguing as before, on the event $\fD_{\lambda}$, any $x\in\mst^{\mvtheta^*}$ is within distance $2\lambda^{-1/2}$ from a point $x'\in\cL\big(\mst^{\mvtheta^*}_{\lambda}\big)$.
If $x'\in\cL\big(\mst^{\mvtheta^*}\big)$, then we are done.
Otherwise, $\mst^{\mvtheta^*}\setminus\{x'\}$ has a component $\cT$ that is a subset of $\mst^{\mvtheta^*}\setminus \mst^{\mvtheta^*}_{\lambda}$.
Further, $\cT$ contains a point $x''\in\cL\big(\mst^{\mvtheta^*}\big)$ such that $x''$ is within distance $3\lambda^{-1/2}$ of $x$.
The rest is routine.

\subsection{The upper Minkowski dimension}\label{sec:upper}
The aim of this section is to prove
\begin{align}\label{eqn:223}
\odim\big(\mst^{\mvtheta^*}\big)\leq 1/\eta\ \ \text{ almost surely.}
\end{align}
We will make use of the following result in our proof.
\begin{prop}\label{lem:ss-14}
There exist $\lambda_{\ref{eqn:6}}\geq\lambda_{\ref{eqn:ss-7a}}$ and 
$\eps_{\ref{eqn:6}}\in(0, \eps_{\ref{eqn:ss-7a}})$ such that for all large $n$ and $\lambda\in[\lambda_{\ref{eqn:6}},\, \eps_{\ref{eqn:6}}n^{\eta}]$,
\begin{align}\label{eqn:6}
\pr\bigg(
\surplus\big(\cC^n_1(\lambda)\big)\geq C_{\ref{eqn:6}}\lambda^{1/\eta}
\bigg)
\leq
\exp\big(-C\lambda\big)\, .
\end{align}
\end{prop}

The proof of Proposition \ref{lem:ss-14} relies on the following result which will also be useful in Section \ref{sec:lower-bound-on-dimenension}.

\begin{lem}\label{lem:1}
There exist $C_{\ref{eqn:1}}, \lambda_{\ref{eqn:1}}\geq 1$ and $\eps_{\ref{eqn:1}}>0$ such that for all large $n$, 
$\lambda\in[\lambda_{\ref{eqn:1}},\, \eps_{\ref{eqn:1}}n^{\eta}]$ and $u\geq C_{\ref{eqn:1}}\lambda^{1/(\tau-3)}$,
\begin{align}\label{eqn:1}
\pr\bigg(\cW_{\mvx}\big(\cC\big(1, \cG_n(\lambda)\big)\big)\geq u\bigg)
\leq 
\exp\big(-Cu\big)\, .
\end{align}
\end{lem}

The proof of Lemma \ref{lem:1} will be given at the end of this section.

\medskip

\noindent{\bf Proof of Proposition \ref{lem:ss-14}:}
In view of \eqref{eqn:712}, it is enough to prove the claimed bound for 
$\spls\big(\cC\big(1, \cG_n(\lambda)\big)\big)$.
Using \eqref{eqn:ss-5} and \eqref{eqn:ss-7}, $\falln$ and for $(\lambda, u)\in\pzI^{\sss(n), 2}$,
\[
\Phi_{\lambda}^{\sss (n)}(u)
\leq 
\lambda u-Cu^{\tau-2}
\leq 
C'\lambda^{\frac{\tau-2}{\tau-3}}\, ,
\]
which together with \eqref{eqn:ss-7a} and the strict concavity of $\Phi_{\lambda}^{\sss (n)}(\cdot)$ implies that $\falln$, and 
$\lambda\in[\lambda_{\ref{eqn:ss-7a}},\, \eps_{\ref{eqn:ss-7a}}n^{\eta}]$,
\begin{align}\label{eqn:ss-55}
\sup_{u\geq 0}\,
\Phi_{\lambda}^{\sss (n)}(u)
=
\sup_{0\leq u\leq s^{\sss (n)}(\lambda)}
\Phi_{\lambda}^{\sss (n)}(u)
\leq 
C_{\ref{eqn:ss-55}}\lambda^{\frac{\tau-2}{\tau-3}}\, .
\end{align}

Recall from \eqref{eqn:8888} the breadth-first walk $Z_{\lambda}^{\sss n, (1)}(\cdot)$ of 
$\cC(1, \cG_n(\lambda))$ started from the vertex $1$.
Let $\text{Gen}_i^{\sss (1)}$ be the set of vertices in the $i$-th generation of the breadth-first tree of $\cC(1, \cG_n(\lambda))$ rooted at the vertex $1$.
Then on the event 
$
\big\{\cW_{\mvx}\big(\cC\big(1, \cG_n(\lambda)\big)\big)
\leq
C_{\ref{eqn:1}}\lambda^{1/(\tau-3)}
\big\}
$,
in the coupling of Lemma \ref{lem:domin-ald-limic-1},
\begin{align}\label{eqn:58-a}
\mathrm{MaxGen}^{(1)}:=
\max\, \bigg\{\sum_{v\in\text{Gen}_j^{\sss(1)} } n^{\eta}x_v\, :\, j\geq 0\bigg\}
\leq
\sup_{0\leq u\leq C_{\ref{eqn:1}}\lambda^{1/(\tau-3)}} 
n^{\eta}\cdot Z_{\lambda}^{\sss n, (1)}(u)\, .
\end{align}
Finally, note that \eqref{eqn:ss-3} and \eqref{eqn:9D} imply, for any $u\geq 0$,
\begin{align}\label{eqn:58-b}
n^{\eta}Z_{\lambda}^{\sss n, (1)}(u)
=
n^{\eta}x_1
+
\sum_{j=2}^n \frac{\theta_{j,\lambda}}{(1+\lambda n^{-\eta})}
\bigg(\ind\big\{\xi^n_j\leq u\big\} - \pr\big(\xi^n_j \leq u\big)\bigg)
+
\Phi_{\lambda}^{\sss (n)}(u)\, .
\end{align}
Combining \eqref{eqn:ss-55}, \eqref{eqn:58-a}, \eqref{eqn:58-b} with Lemma \ref{lem:1}, and applying Lemma \ref{lem:ss-lem-3A} with 
$s=C_{\ref{eqn:1}}\lambda^{1/(\tau-3)}$ and $y=C_{\ref{eqn:ss-55}}\big(C_{\ref{eqn:1}}\big)^{-\frac{\tau-3}{2}}\lambda$,
we see that there exists $\lambda_{\ref{eqn:ss-56}}\geq\lambda_{\ref{eqn:ss-7a}}\vee\lambda_{\ref{eqn:1}}$ 
such that $\falln$ and $\lambda\in[\lambda_{\ref{eqn:ss-56}},\, \eps_{\ref{eqn:6}}n^{\eta}]$,
\begin{align}\label{eqn:ss-56}
\pr\bigg(
\bigg\{
\cW_{\mvx}\bigg(
\big(\cC(1, \cG_n(\lambda))\big)\bigg)
\geq 
C_{\ref{eqn:1}}\lambda^{\frac{1}{\tau-3}}
\bigg\}
\bigcup
\bigg\{
\mathrm{MaxGen}^{(1)}
\geq 2C_{\ref{eqn:ss-55}}\lambda^{\frac{\tau-2}{\tau-3}}
\bigg\}
\bigg)
\leq
\exp\big(-C\lambda\big)\, ,
\end{align}
where
$
\eps_{\ref{eqn:6}}:=\eps_{\ref{eqn:ss-7a}}\wedge\eps_{\ref{eqn:1}}
$.
Conditional on the breadth-first tree $T$ (note that $T$ is a plane tree with vertex labels) of $\cC(1, \cG_n(\lambda))$, the rest of the edges in $\cC(1, \cG_n(\lambda))$ can be generated by placing edges independently with probability $1-\exp\big(-(\lambda+n^{\eta})x_u x_v\big)$ for every $u, v\in V(T)$, where $u\neq v$, and either
(i) $u, v\in\text{Gen}_j^{\sss(1)}$ for some $j\geq 1$, or 
(ii) $u\in\text{Gen}_{j+1}^{\sss(1)}$ and $v\in\text{Gen}_j^{\sss(1)}$ for some $j\geq 1$, and $v$ lies on the right side of the ancestral line of $u$.
Using \eqref{eqn:88}, we see that conditional on $T$, on the complement of the event in \eqref{eqn:ss-56}, 
for any $\lambda\in[\lambda_{\ref{eqn:ss-56}},\, \eps_{\ref{eqn:6}}n^{\eta}]$,
\[
\spls\big(\cC(1, \cG_n(\lambda))\big)
\stod
\Poi\big((\eps_{\ref{eqn:6}}+1)\sum_{u\in V(T)}x_u\cdot 2\cdot\text{MaxGen}^{(1)}\big)
\stod
\Poi\big(C\lambda^{1/\eta}\big)\, .
\]
Now the proof of Proposition \ref{lem:ss-14} can be completed by combining \eqref{eqn:ss-56} with standard tail bounds for a Poisson random variable.
\qed

\medskip

We are now ready to prove the claimed upper bound on $\odim\big(\mst^{\mvtheta^*}\big)$.
Recall the  process $\mtbp$ from Definition \ref{defn:mtbp}.
Fix $\Delta\in(0, 1/2]$.
Then for any $\lambda\geq 2^{\eta}$,
\begin{align}\label{eqn:4}
&\pr\bigg(\diam\big(G_n(\lambda)\setminus\big[\lambda^{\frac{1}{\eta}(1+\Delta)}\big]\big)
\geq n^{\eta}/\lambda\bigg)
\leq
\sum_{\lambda^{(1+\Delta)/\eta}\leq i\leq n}
\pr\bigg(\diam\bigg(\cC\big(i, G_n(\lambda)\setminus [i-1]\big)\bigg)\geq n^{\eta}/\lambda\bigg)
\notag\\
&\hskip100pt\leq
\sqrt{\e}\sum_{\lambda^{(1+\Delta)/\eta}\leq i\leq n}
\pr\bigg(
\height\bigg(\mtbp_n^i\big([n]\setminus [i-1]\big)\bigg)\geq n^{\eta}/(2\lambda)
\bigg)\, ,
\end{align}
where the last step uses arguments similar to the ones used in the proof of Lemma \ref{lem:26}.
To bound the summands in the last step in \eqref{eqn:4}, we can use an argument similar to the one used in Section \ref{sec:diam}.

\begin{defn}\label{def:two-layer-cake}
For $\lambda\geq 1$ and $2\leq i\leq n$, consider the following (potentially) two layer process $\twolbp_n^{i}(\lambda)$: 
\begin{enumeratea}
\item {\bf Layer 1:} 
Start $\mtbp_n^{i}([n]\setminus [i-1])$, and run this process up to generation $n^{\eta}/(4\lambda)$. 
Call this the first layer. 
If there is at least one vertex in generation $n^{\eta}/(4\lambda)$, then we say that the first layer has been fully activated. 
\item {\bf Layer 2:} 
If the first layer is fully activated, then starting from every vertex $v$ in generation  $n^{\eta}/(4\lambda)$, run independent $\mtbp_n^{\mathrm{type}(v)}\big([n]\setminus \{1\}\big)$ processes.
Call this the second layer. 
If any of these branching processes survives up to generation $n^{\eta}/(4\lambda)$, then we say that the second layer has been fully activated.
\end{enumeratea}
\end{defn}

Let $\gamma_*$ be as in Proposition \ref{prop:lem-6-bp-height}.
Then $\falln$, for any $2\leq i\leq n$ and $\lambda\in[1/\gamma_*\, ,\, n^{\eta}/10]$,
\begin{align}\label{eqn:3}
&\pr\bigg(
\height\bigg(\mtbp_n^i\big([n]\setminus [i-1]\big)\bigg)\geq n^{\eta}/(2\lambda)
\bigg)
\leq
\pr\bigg(
\height\big(\twolbp_n^i(\lambda)\big)\geq n^{\eta}/(2\lambda)
\bigg)\notag\\
&\hskip20pt\leq
Cw_i\bigg(\sum_{j\geq i}w_j^2/(\nu_n\ell_n)\bigg)^{\frac{n^{\eta}}{4\lambda}}
\big(4\lambda\big)^{\frac{1}{\tau-3}}\cdot n^{-\alpha}
\leq
C' i^{-\alpha}\cdot\lambda^{\frac{1}{\tau-3}}\cdot\exp\big(-C''i^{\eta}/\lambda\big)\, ,
\end{align}
where we have used Proposition \ref{prop:lem-6-bp-height} and arguments similar to the ones used in the proof of Lemma \ref{lem:fnl-bound}.
Combining \eqref{eqn:3} with \eqref{eqn:4}, and calculations similar to the ones in the proof of Lemma \ref{lem:fnl-bound} show that $\falln$ and for $\lambda\in[2^{\eta}\vee(\gamma_*)^{-1},\, n^{\eta}/10]$,
\begin{align}\label{eqn:7}
\pr\bigg(\diam\big(G_n(\lambda)\setminus\big[\lambda^{\frac{1}{\eta}(1+\Delta)}\big]\big)
\geq n^{\eta}/\lambda\bigg)
\leq
C_{\ref{eqn:7}}(\Delta)\cdot\exp\big(-C\lambda^{\Delta}\big)\, .
\end{align}

Recall the notation $\cN(\cdot\, ,\, \cdot)$ introduced at the beginning of Section \ref{sec:res}.
Now note that on the complement of the event in \eqref{eqn:7},
\begin{align}\label{eqn:8}
\cN\big(\cC_1^n(\lambda),\, 2n^{\eta}/\lambda\big)
\leq 
\lambda^{(1+\Delta)/\eta}\, .
\end{align}
By \eqref{eqn:75-a}, for any $k\geq 1$ and $\lambda\geq 1$,
\begin{align}\label{eqn:9}
\pr\big(
\cN\big(\mst_{\lambda}^{\mvtheta^*},\, 3/\lambda\big)\geq k
\big)
\leq
\limsup_{n\to\infty}\,
\pr\bigg(
\cN\bigg(\cbd^{\infty}\big(\cC_1^n(\lambda)\big),\, 2n^{\eta}/\lambda\bigg)\geq k
\bigg)\, .
\end{align}
Since $\cbd^{\infty}\big(\cC_1^n(\lambda)\big)$ can be covered by 
$\cN\big(\cC_1^n(\lambda), 2n^{\eta}/\lambda\big)
+
\spls\big(\cC_1^n(\lambda)\big)$
many $2n^{\eta}/\lambda$ balls, combining \eqref{eqn:7}, \eqref{eqn:8}, and \eqref{eqn:9} with Proposition \ref{lem:ss-14}, we see that there exists $\lambda_{\ref{eqn:700}}\geq 1$ such that for all $\lambda\geq\lambda_{\ref{eqn:700}}$,
\begin{align}\label{eqn:700}
\pr\big(
\cN\big(\mst_{\lambda}^{\mvtheta^*},\, 3/\lambda\big)
>
2\lambda^{(1+\Delta)/\eta}
\big)
\leq
\big(1+C_{\ref{eqn:7}}(\Delta)\big)\exp\big(-C\lambda^{\Delta}\big)\, .
\end{align}
In the coupling on the right side of \eqref{eqn:53},
\[
\cN\bigg(\mst^{\mvtheta^*},\, 3\lambda^{-1}
+d_{\rH}\big(\mst^{\mvtheta^*},\, \mst_{\lambda}^{\mvtheta^*}\big)\bigg)
\leq
\cN\big(\mst_{\lambda}^{\mvtheta^*},\, 3\lambda^{-1}\big)\, ,
\]
which combined with \eqref{eqn:700} and Theorem \ref{thm:mst-crit-gh} shows that for all $\lambda\geq\lambda_{\ref{eqn:555}}(\Delta)\vee\lambda_{\ref{eqn:700}}$,
\begin{align}\label{eqn:701}
\pr\big(
\cN\big(\mst^{\mvtheta^*},\, 4\lambda^{-1+\Delta}\big)
>
2\lambda^{(1+\Delta)/\eta}
\big)
\leq
C_{\ref{eqn:701}}(\Delta)\cdot\lambda^{-1/2}\, .
\end{align}
Replacing $\lambda$ by $k^4$ in \eqref{eqn:701} for large integer values of $k$, an application of the Borel-Cantelli lemma gives
\begin{align}\label{eqn:702}
\limsup_{k\to\infty}\, 
\frac{\log \cN\big(\mst^{\mvtheta^*},\, 4k^{-4(1-\Delta)}\big)}{\log\big(k^{4(1-\Delta)}/4\big)}
\leq
\frac{1+\Delta}{\eta(1-\Delta)}\ \ \ \text{ almost surely.}
\end{align}
Sandwiching $\delta$ between $4k^{-4(1-\Delta)}$ and $4(k+1)^{-4(1-\Delta)}$ and letting $\delta\downarrow 0$,  we can show using \eqref{eqn:702} that $\odim\big(\mst^{\mvtheta^*}\big)\leq (1+\Delta)/(\eta(1-\Delta))$ almost surely.
We complete the proof of \eqref{eqn:223} by letting $\Delta\downarrow 0$.

\medskip

\noindent{\bf Proof of Lemma \ref{lem:1}:}
By \eqref{eqn:ss-7}, $\falln$ and for $(\lambda, u)\in\pzI^{\sss (n), 2}$,
\begin{align}\label{eqn:2}
\varphi_{\lambda}^{\sss (n)}(u)\geq 2C_{\ref{eqn:2}}u^{\tau-3}\, .
\end{align}
Thus, using \eqref{eqn:ss-5} and \eqref{eqn:2}, we see that 
$\falln$, 
\begin{align*}
\Phi_{\lambda}^{\sss(n)}(u)
\leq 
\lambda u-10\cdot u \varphi_{\lambda}^{\sss (n)}(u)/11
\leq
u\big(\lambda-20 C_{\ref{eqn:2}}u^{\tau-3}/11\big)
\leq
-C_{\ref{eqn:2}}u^{\tau-2}\, ,
\end{align*}
whenever $(\lambda, u)\in\pzI^{\sss (n), 2}$ and $9C_{\ref{eqn:2}}u^{\tau-3}\geq 11\lambda$,
and consequently,
\begin{align*}
&\pr\bigg(\cW_{\mvx}\big(\cC\big(1, \cG_n(\lambda)\big)\big)\geq u\bigg)
\leq
\pr\big(Z_{\lambda}^{\sss n, (1)}(u)\geq 0\big)\\
&\hskip40pt\leq
\pr\bigg(n^{\eta} x_1
+\sum_{j=2}^n\theta_{j, \lambda}\bigg(
\ind_{\{\xi_j^n\leq u\}}-\pr\big(\xi_j^n\leq u\big)
\bigg)
\geq 
C_{\ref{eqn:2}}u^{\tau-2}
\bigg)\, ,
\end{align*}
where the first step uses Lemma \ref{lem:domin-ald-limic-1}, and the last step uses \eqref{eqn:58-b}.
Now an application of Lemma \ref{lem:ss-lem-3A} will yield the bound in \eqref{eqn:1}, and this bound will be valid for all large $n$, for
$\lambda\in[\lambda_{\ref{eqn:1}},\, \eps_{\ref{eqn:1}}n^{\eta}]$ and 
$
u\in\big[C_{\ref{eqn:1}}\lambda^{1/(\tau-3)}\, ,\, n^{\alpha}\big(\sigma_2^{\sss (n)}\big)^{1/2}/(A_2 2^{\alpha+1})\big]
$,
where $C_{\ref{eqn:1}}=\big(11/(9 C_{\ref{eqn:2}})\big)^{1/(\tau-3)}$, and $\lambda_{\ref{eqn:1}}$ and $\eps_{\ref{eqn:1}}$ are appropriately chosen constants.
Finally, this bound can be extended to all $u\geq C_{\ref{eqn:1}}\lambda^{1/(\tau-3)}$ by simply noting that the probability on the left side of \eqref{eqn:1} is zero if $u>\sum_{i}x_i\asymp n^{\alpha}$.
\qed

\subsection{The lower Minkowski dimension}\label{sec:lower-bound-on-dimenension}
In this section we will prove
\begin{align}\label{eqn:223-a}
\udim\big(\mst^{\mvtheta^*}\big)\geq 1/\eta\ \ \text{ almost surely,}
\end{align}
which combined with \eqref{eqn:223} will complete the proof of Theorem \ref{thm:mst-convg}\,\eqref{it:c}.
To that end, let us first introduce some notation.
For disjoint $A, B\subseteq [n]$ and $r\geq 1$, we write $A\connects{r} B$ to mean there exist $1\leq t\leq r$, $v_0\in A$, $v_t\in B$, and $v_1,\ldots,v_{t-1}\in [n]\setminus\{1\}$ such that the edges $\{v_i, v_{i+1}\}$, $i=0,\ldots, t-1$, are present in $G_n(\lambda)$.
If $A=\{i\}$, then we simply write $i\connects{r} B$ instead of $\{i\}\connects{r} B$.

Next, note that if $\cN(X, u)\leq k$ for some metric space $(X, d)$, $u\geq 0$, and $k\geq 1$, then for any $x_1,\ldots, x_{k+1}\in X$, there exist $1\leq i<j\leq k+1$ such that $d(x_i, x_j)\leq 2u$.
Hence, for any $x_1,\ldots, x_{2k}\in X$,
\begin{align}\label{eqn:703}
\#\big\{2\leq i\leq 2k\, :\, \min_{j\in [i-1]}d(x_i, x_j)\leq 2u \big\}\geq k\, ,
\end{align}
since otherwise, we would be able to find $1= i_1<i_2<\ldots<i_{k+1}\leq 2k$ such that the minimum of the pairwise distances between these $k+1$ points is more than $2u$.

For $\Delta\in(0, 1/2]$, choose $\lambda_{\ref{eqn:ss-59}}(\Delta)\geq\lambda_{\ref{it:lem-ss-4}}\vee\lambda_{\ref{eqn:1}}$ such that $\lambda^{1/\eta}\big(\log\lambda\big)^{-3}>\lambda^{(1-\Delta)/\eta}$ for all $\lambda\geq\lambda_{\ref{eqn:ss-59}}(\Delta)$.
For $\Delta\in(0, 1/2]$ and $\lambda\geq\lambda_{\ref{eqn:ss-59}}(\Delta)$, choose $n_{\ref{eqn:ss-59}}(\lambda)\geq 1$ such that 
$\lambda<\big(\eps_{\ref{eqn:ss-7a}}\wedge\eps_{\ref{eqn:1}}\big)n^{\eta}$ for all 
$n\geq n_{\ref{eqn:ss-59}}(\lambda)$.
(Thus, these thresholds are chose in a way so that the bounds in \eqref{it:lem:ss-8} and \eqref{eqn:1} are applicable whenever $\lambda\geq\lambda_{\ref{eqn:ss-59}}$ and $n\geq n_{\ref{eqn:ss-59}}$.)
Letting
\[
\cD_i:=\big\{[i]\subseteq\cC_1^n(\lambda)\, ,\ \text{ and }\  i\connectsn [i-1]\big\}\, ,
\]
we see that for $\lambda\geq\lambda_{\ref{eqn:ss-59}}$, $n\geq n_{\ref{eqn:ss-59}}$, and 
\begin{align}\label{eqn:def-pza}
\pza:=2\Delta/(\tau-3)\, ,
\end{align}
\begin{align}\label{eqn:ss-59}
&\pr\bigg(
\cN\big(n^{-\eta}\cdot\cC^n_1(\lambda)\, ,\, \big(2\lambda^{1+\pza}\big)^{-1}\big)
\leq 2^{-1}\lambda^{\frac{1}{\eta}(1-\Delta)}
\bigg)\notag\\
&\hskip30pt\leq
\pr\bigg(
\bigg\{
\#\big\{
2\leq i\leq\lambda^{\frac{1}{\eta}(1-\Delta)}\, :\, i\connectsn [i-1]
\big\}
\geq 2^{-1}\lambda^{\frac{1}{\eta}(1-\Delta)}
\bigg\}\bigcap
\bigg\{
i\in\cC_1^n(\lambda)\ \forall\ i\leq\lambda^{\frac{1}{\eta}(1-\Delta)}
\bigg\}
\bigg)\notag\\
&\hskip60pt
+
\pr\bigg(
\exists i\leq\lambda^{\frac{1}{\eta}}\big(\log\lambda\big)^{-3} \text{ such that }i\notin\cC_1^n(\lambda)
\bigg)\notag\\
&\hskip90pt
\leq
\pr\bigg(
\bigg\{
\#\big\{
2\leq i\leq\lambda^{\frac{1}{\eta}(1-\Delta)}\, :\, \cD_i\ \text{ holds}\,
\big\}
\geq 2^{-1}\lambda^{\frac{1}{\eta}(1-\Delta)}
\bigg\}
\bigg)
+\exp\big(-C\log^{3\alpha}\lambda\big)
\notag\\
&\hskip120pt\leq
2\lambda^{-\frac{1}{\eta}(1-\Delta)}\cdot
\sum_{i=2}^{\lambda^{\frac{1}{\eta}(1-\Delta)}} \pr\big(\cD_i\big)
+\exp\big(-C\log^{3\alpha}\lambda\big)\, ,
\end{align}
where the first step uses the observation around \eqref{eqn:703}, the second step uses \eqref{it:lem:ss-8}, and the last step uses Markov's inequality.

\begin{figure}
	\centering
	\includegraphics[trim=2cm 7.5cm 1cm 13.5cm, clip=true, angle=0, scale=.82]{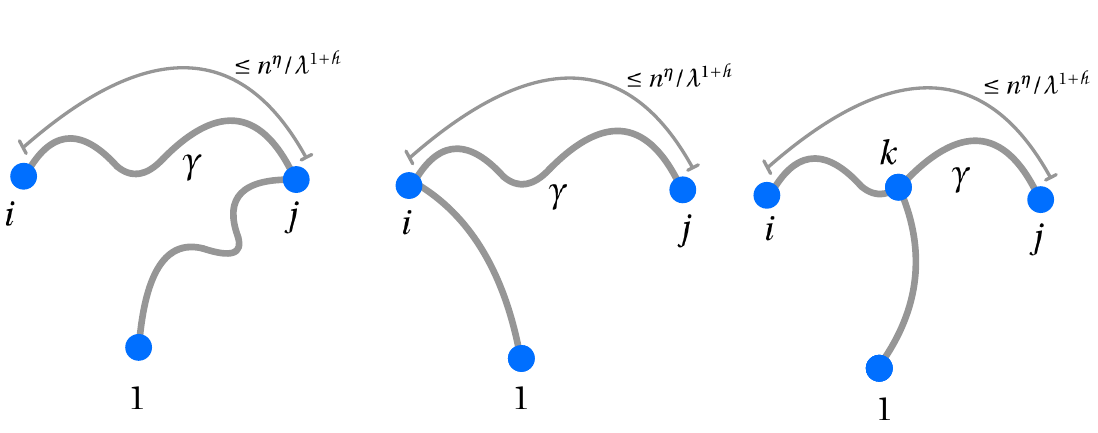}
	\captionof{figure}{The vertex $j\in[i-1]$ satisfies $\dist(i, j)=\min_{j'\in [i-1]}\dist(i, j')$, where $\dist$ denotes the graph distance in $G_n(\lambda)$. 
		In the leftmost figure, $j=1$ is possible, and $2\leq j\leq i-1$ in the other two figures.
		In the rightmost figure, $i+1\leq k\leq n$.}
	\label{fig:1}
\end{figure}

On the event $\cD_i$, consider a shortest path $\gamma$ from $i$ to $[i-1]$.
Then depending on where the shortest path from $1$ to $\gamma$ meets $\gamma$, we may consider three possibilities as in Figure \ref{fig:1}.
Hence,
\begin{align}\label{eqn:ss-60}
\cD_i
&\subseteq
\bigg(\bigcup_{j=1}^{i-1}\bigg(\big\{j\in\cC_1^n(\lambda)\big\}\circ
\big\{j\connectsn i\big\}\bigg)\bigg)\,
\bigcup\,
\bigg(\bigcup_{j=2}^{i-1}\bigg(\big\{i\in\cC_1^n(\lambda)\big\}\circ
\big\{j\connectsn i\big\}\bigg)\bigg)\notag\\
&\hskip50pt
\bigcup
\bigg(\bigg(
\bigcup_{j=2}^{i-1}\ \bigcup_{k=i+1}^{n}
\big\{k\in\cC_1^n(\lambda)\big\}\circ\big\{k\connectsn i\big\}\circ\big\{k\connectsn j\big\}
\bigg)\bigg)\, ,
\end{align}
where $\circ$ denotes disjoint occurrence of events.
Writing $\sum_\square$ for sum over $t=0, 1,\ldots, -1+n^{\eta}/\lambda^{1+\pza}$, we have, for any $i\neq j$,
\begin{align}\label{eqn:ss-61}
\pr\big(i\connectsn j\big)
&\leq
\sum\displaystyle_{\square}\ \sum_{2\leq v_1,\ldots, v_t\leq n}
\bigg(p_{\lambda}^n\frac{w_i w_{v_1}}{\ell_n}\bigg)
\bigg(p_{\lambda}^n\frac{w_{v_1} w_{v_2}}{\ell_n}\bigg)\ldots
\bigg(p_{\lambda}^n\frac{w_{v_t}w_j}{\ell_n}\bigg)\notag\\
&\leq
\sum\displaystyle_{\square}\ \bigg(1+\frac{\lambda}{n^{\eta}}\bigg)^{t+1}
\bigg(\frac{w_i w_j}{\nu_n\ell_n}\bigg)
\sum_{2\leq v_1,\ldots, v_t\leq n}\bigg(\frac{w_{v_1}^2}{\nu_n\ell_n}\bigg)\ldots \bigg(\frac{w_{v_t}^2}{\nu_n\ell_n}\bigg)\notag\\
&\leq
\sum\displaystyle_{\square}\ \bigg(1+\frac{\lambda}{n^{\eta}}\bigg)^{\frac{n^{\eta}}{\lambda^{1+\pza}}}
\bigg(\frac{w_i w_j}{\nu_n\ell_n}\bigg)
\leq\frac{n^{\eta}}{\lambda^{1+\pza}}\cdot\e\cdot\bigg(\frac{w_i w_j}{\nu_n\ell_n}\bigg)
\leq\frac{C}{\lambda^{1+\pza}i^{\alpha}j^{\alpha}}\, ,
\end{align}
where the third step uses \eqref{eqn:741-a}, and the last step uses Assumption \ref{ass:wts}\,(iii) and \eqref{eqn:up-bnd}.
Note also that for $2\leq j\leq n$,
\begin{align}\label{eqn:ss-62a}
\pr\big(j\in\cC_1^n(\lambda)\big)
&\leq
\sum_{j'\in[n]}\pr\big(j'\in\cC_1^n(\lambda)\big)\pr\big(\{j,j'\}\in E(G_n(\lambda))\big)
\leq
\sum_{j'\in[n]}p_{\lambda}^n\frac{w_j w_{j'}}{\ell_n}\pr\big(j'\in\cC_1^n(\lambda)\big)\notag\\
&=
p_{\lambda}^n\frac{w_j}{\ell_n}\cdot\bE\big[\cW\big(\cC_1^n(\lambda)\big)\big]
\leq 
Cp_{\lambda}^n\frac{w_j}{\ell_n}\cdot n^{\rho}\lambda^{1/(\tau-3)}
\leq
C'\frac{\lambda^{1/(\tau-3)}}{j^{\alpha}}\, ,
\end{align}
where the penultimate uses Lemma \ref{lem:1} and is valid for $\lambda\geq\lambda_{\ref{eqn:ss-59}}$ and $n\geq n_{\ref{eqn:ss-59}}$.
Thus, combining \eqref{eqn:ss-60}, \eqref{eqn:ss-61}, and \eqref{eqn:ss-62a}, an application of the BK inequality shows, for $\lambda\geq\lambda_{\ref{eqn:ss-59}}$, $n\geq n_{\ref{eqn:ss-59}}$, and $2\leq i\leq n$,
\begin{align}\label{eqn:704}
\pr\big(\cD_i\big)
&\leq
C\sum_{j=1}^{i-1}
\bigg(\ind_{\{j=1\}}+\frac{\lambda^{1/(\tau-3)}}{j^{\alpha}}\bigg)
\frac{1}{\lambda^{1+\pza}i^{\alpha}j^{\alpha}}
+
C\sum_{j=2}^{i-1}\sum_{k=i+1}^n
\bigg(\frac{\lambda^{1/(\tau-3)}}{k^{\alpha}}\cdot
\frac{1}{\lambda^{1+\pza}i^{\alpha}k^{\alpha}}\cdot
\frac{1}{\lambda^{1+\pza}j^{\alpha}k^{\alpha}}
\bigg)\notag\\
&\hskip80pt\leq
C'\bigg(
\lambda^{\frac{4-\tau}{\tau-3}-\pza}\cdot i^{\frac{\tau-4}{\tau-1}}
+
\lambda^{\frac{7-2\tau}{\tau-3}-2\pza}\cdot i^{\frac{2\tau-7}{\tau-1}}
\bigg)\, ,
\end{align}
where the last step follows from some routine calculations.
Combining \eqref{eqn:704} with \eqref{eqn:ss-59} shows that for $\lambda\geq\lambda_{\ref{eqn:ss-59}}$ and $n\geq n_{\ref{eqn:ss-59}}$,
\begin{align}\label{eqn:705}
&\pr\bigg(
\cN\big(n^{-\eta}\cdot\cC^n_1(\lambda)\, ,\, \big(2\lambda^{1+\pza}\big)^{-1}\big)
\leq 2^{-1}\lambda^{\frac{1}{\eta}(1-\Delta)}
\bigg)\notag\\
&\hskip20pt\leq
\exp\big(-C\log^{3\alpha}\lambda\big)+
C'\bigg[
\lambda^{\frac{4-\tau}{\tau-3}-\pza}\cdot
\bigg(\lambda^{\frac{1}{\eta}(1-\Delta)}\bigg)^{\frac{\tau-4}{\tau-1}}
+
\lambda^{\frac{7-2\tau}{\tau-3}-2\pza}\cdot \bigg(\lambda^{\frac{1}{\eta}(1-\Delta)}\bigg)^{\frac{2\tau-7}{\tau-1}}
\bigg]\notag\\
&\hskip40pt=
\exp\big(-C\log^{3\alpha}\lambda\big)+
C'\bigg[
\lambda^{\Delta\big(\frac{4-\tau}{\tau-3}\big)-\pza}
+
\lambda^{\Delta\big(\frac{7-2\tau}{\tau-3}\big)-2\pza}
\bigg]\notag\\
&\hskip60pt\leq
\exp\big(-C\log^{3\alpha}\lambda\big)
+
C''\lambda^{-\Delta/(\tau-3)}\, ,
\end{align}
where the last step uses \eqref{eqn:def-pza}.

By \eqref{eqn:64-a}, for any $k\geq 1$ and $\lambda\geq 1$,
\begin{align*}
\pr\bigg(
\cN\big(\pzS_{\lambda}^{\mvtheta^*},\, \big(4\lambda^{1+\pza}\big)^{-1}\big)\leq k
\bigg)
\leq
\limsup_{n\to\infty}\,
\pr\bigg(
\cN\bigg(n^{-\eta}\cdot\cC_1^n(\lambda),\, \big(2\lambda^{1+\pza}\big)^{-1}\bigg)\leq k
\bigg)\, ,
\end{align*}
which together with \eqref{eqn:705} yields, for $\lambda\geq\lambda_{\ref{eqn:ss-59}}$,
\begin{align}\label{eqn:706}
\pr\bigg(
\cN\big(\pzS_{\lambda}^{\mvtheta^*},\, \big(4\lambda^{1+\pza}\big)^{-1}\big)
\leq
2^{-1}\lambda^{\frac{1}{\eta}(1-\Delta)}
\bigg)
\leq 
\exp\big(-C\log^{3\alpha}\lambda\big)
+
C'\lambda^{-\Delta/(\tau-3)}\, .
\end{align}
Note that 
$
\cN\big(\pzS_{\lambda}^{\mvtheta^*}, r\big)
\stod
\cN\big(\mst_{\lambda}^{\mvtheta^*}, r\big)
\stod
\cN\big(\mst^{\mvtheta^*}, r\big)
$
for any $r>0$, where the first relation follows since 
$\mst^{\mvtheta^*}_{\lambda}=\cb^{\infty}\big(\pzS_{\lambda}^{\mvtheta^*}\big)$,
and the second relation follows from the existence of the coupling between $\mst^{\mvtheta^*}_{\lambda}$ and $\mst^{\mvtheta^*}$ as appears on the right side of \eqref{eqn:53}.
Hence, \eqref{eqn:706} continues to hold if we replace $\pzS_{\lambda}^{\mvtheta^*}$ by $\mst^{\mvtheta^*}$.
Let $g=g(\Delta):=2(\tau-3)/\Delta$.
Using \eqref{eqn:706} with $\mst^{\mvtheta^*}$ replacing $\pzS_{\lambda}^{\mvtheta^*}$, and letting $\lambda\to\infty$ along the sequence $\big(k^g,\, k\geq 1\big)$, we get, via an application of the Borel-Cantelli lemma,
\[
\liminf_{k\to\infty}\
\frac{\log\cN\bigg(\mst^{\mvtheta^*},\, \big(4k^{g(1+\pzh)}\big)^{-1} \bigg)}{\log\big(4k^{g(1+\pzh)}\big)}
\geq
\frac{1-\Delta}{\eta(1+\pza)}\ \ \text{ almost surely}
\]
for any $\Delta\in(0, 1/2]$.
Sandwiching $\delta$ between $\big(4k^{g(1+\pzh)}\big)^{-1} $ and $\big(4(k+1)^{g(1+\pzh)}\big)^{-1}$ and letting $\delta\downarrow 0$, we conclude that for any $\Delta\in(0, 1/2]$,
\[
\udim\big(\mst^{\mvtheta^*}\big)
\geq
\frac{1-\Delta}{\eta(1+\pza)}\ \ \text{ almost surely.}
\]
We complete the proof of \eqref{eqn:223-a} by letting $\Delta\downarrow 0$, and using \eqref{eqn:def-pza}.

\section{Discussion}
\label{sec:disc}
The random graph models considered in this paper are closely related to the multiplicative coalescent--a fact that was crucial in our proof.
It is easy to argue heuristically that around the point of phase transition, many standard models of dynamic random graphs evolve roughly like the multiplicative coalescent. 
This intuition was formalized in \cite{SBSSXW-universal,baslingker-bhamidi-broutin-sen-wang}, where it was shown that the metric scaling limits for maximal components in the critical regime for a number of random graph models have the same limit (up to constant scaling factors) as that of the \erdos random graph established in \cite{BBG-12}. 
These models include 
critical percolation on edge weighted graphs converging to an $L^3$ graphon,
Aldous and Pittel's RGIV model \cite{aldous-pittel},
the stochastic block model,
and the configuration model (under appropriate moment assumptions).
In the heavy-tailed regime as considered in this paper, metric scaling limits for maximal components of critical inhomogeneous random graphs were first established in \cite{SB-vdH-SS-PTRF}. 
These results were then leveraged in \cite{SB-SD-vdH-SS} to establish the scaling limit of the critical configuration model.

We expect a similar program to be carried out building on the results of this paper for establishing universality of the MST for a host of random graph models. 
Two important models for which this problem remains open are 
(i) the configuration model and (simple) random graphs with given degree sequence with tail exponent $\tau\in(3, 4)$, and
(ii) a sequence of edge weighted graphs converging to a graphon whose leading eigenfunction is an element of $L^p[0, 1]\setminus L^3[0, 1]$ for some $p\in(2, 3)$.
We expect that for both these models, the scaling limit of the MST on the giant components of these models, under suitable assumptions, will be the same as the ones obtained in this paper.

We now briefly remark on the assumptions in this paper.
As mentioned before, Assumption \ref{ass:wts}~(i) implies supercriticality of the random graph model, whereas Assumption \ref{ass:wts}~(ii) corresponds to the condition in \cite[Display (19)]{aldous-limic}.
It should be possible to relax Assumption \ref{ass:wts}~(iii) and Assumption \ref{ass:wts}~(iv).
For example, our proof of Proposition \ref{prop:lem-6-bp-height}--a key ingredient used to establish tail bounds on the diameter outside the component of the vertex $1$--does not require the full force of Assumption \ref{ass:wts}~(iii).
This proof only uses the relation
\[
\pr\big(V_n\geq u\big)\asymp 1/u^{\tau-2}\  \ \text{ for all }\ \ n\geq 2,\ \ \text{ and }\ \ u\in[1 ,\, v_2]\, .
\]
Similarly, it should be possible to relax Assumption \ref{ass:wts}\,(iii) and Assumption \ref{ass:wts}\,(iv) and still carry out the other steps in our proof at the cost of a more intricate analysis.
In the context of critical inhomogeneous random graphs, \cite{broutin2020limits} establishes convergence of the maximal components as well as compactness of the scaling limit under a weaker integrability condition \cite[Display (9)]{broutin2020limits}.
For the critical heavy-tailed configuration model, \cite{SB-SD-vdH-SS-2020} proves GHP convergence of the maximal components under \cite[Assumptions 1 and 2]{SB-SD-vdH-SS-2020}, and \cite{conchon-goldschmidt} works with i.i.d. degrees having a power-law.
However, all these papers deal only with the critical regime.
It is not clear what the necessary and sufficient conditions are for results such as Theorem \ref{thm:mst-convg} to hold, as proving such a result requires careful control over the complement of the maximal component in the critical window as well as the barely supercritical regime.
Establishing such a necessary and sufficient condition remains a challenging open problem.

We close this section with a discussion on the possible extension of the convergence in \eqref{eqn:666} with respect to the stronger GHP topology. 
Let $\mvmu^n$ (resp. $\tbar{\mvmu}^n$) denote the uniform probability measure on the vertices of $\mvM^n$ (resp. $\tbar{\mvM}^n$).
View $\big(n^{-\eta}\cdot \mvM^n,\, \mvmu^n\big)$ and $\big(n^{-\eta}\cdot \tbar{\mvM}^n,\, \tbar{\mvmu}^n\big)$ as random metric measure spaces.
\begin{conj}\label{conj:mst-ghp-convg}
	Under Assumption \ref{ass:wts} on the weight sequence, there exists a random compact measured $\bR$-tree $\big(\mst^{\mvtheta^*},\, \mu\big)$ whose law depends only on $\mvtheta^*$ such that
	\begin{align}\label{eqn:666-a}
		\big(n^{-\eta}\cdot \mvM^n,\, \mvmu^n\big) 
		\weakc
		\big(\mst^{\mvtheta^*},\, \mu\big) , \ \ \text{ as }\ \ \ n\to\infty\, ,
	\end{align}
	with respect to the GHP topology.
	Further, almost surely,
	the measure $\mu$ is nonatomic, i.e.,
	\begin{align}\label{eqn:667}
		\pr\big(\mu(\{x\})=0\ \text{ for all }\ x\in\mst^{\mvtheta^*}\big)=1\, ,
	\end{align}
	and $\mu$ is concentrated on the set of leaves of $\mst^{\mvtheta^*}$, i.e.,
	\begin{align}\label{eqn:668}
		\pr\bigg( \mu\big(\cL(\mst^{\mvtheta^*}) \big) = 1 \bigg)=1\, .
	\end{align}
	Moreover, \eqref{eqn:666-a} continues to hold under Assumption \ref{ass:wts} if we replace the left side by $\big(n^{-\eta}\cdot \tbar{\mvM}^n,\, \tbar{\mvmu}^n\big)$.
\end{conj}
In the context of the complete graph, the analogues of \eqref{eqn:666-a} and \eqref{eqn:667} were proved in \cite{AddBroGolMie13} and \cite{addarioberry-sen} respectively.
In order to prove \eqref{eqn:666-a}, it is enough to establish the following:
Let $\tbar T_{\lambda; i}^n$~, $i=1,\ldots, \tbar k_{\lambda}^n$, be the trees in the forest obtained by removing from $\tbar M^n$ the vertices in $\tbar M_{\lambda}^n$ and all edges incident to the vertices in $\tbar M_{\lambda}^n$.
Let $\tbar Y^n_{\lambda}:= \max_{1\leq i\leq \tbar k_{\lambda}^n} \big|\tbar T_{\lambda; i}^n\big|$.
Then for all $\eps>0$,
\begin{align}\label{eqn:486}
	\lim_{\lambda\to\infty}\limsup_{n\to\infty}\, 
	\pr\big( \tbar Y^n_{\lambda}>\eps n\big)=0\, .
\end{align}
As mentioned in Section~\ref{sec:proof-strategy}, in the setting of the complete graph, the analogoue of \eqref{eqn:486} was established in \cite[Lemma~4.11]{AddBroGolMie13}. 
We briefly describe why proving \eqref{eqn:486} suffices.

Consider a degree sequence $\vd=\big(d_1, \ldots, d_n\big)$, and let $\cmnd$ (resp. $\cgnd$) be a configuration model (resp. uniform simple random graph) with degree sequence $\vd$.
(The reader can consult \cite[Section 1.2 and Display (6.48)]{addarioberry-sen} for a quick overview of the definitions and the relevant properties of these models. 
We refer the reader to \cite{Hofs15, Hofs17} for a more detailed treatment.)
Let $U_e$, $e\in E(\cmnd)$, be i.i.d. $\mathrm{Uniform}[0, 1]$ random variables conditional on $\cmnd$.
Let $M^{\vd}$ be the MST of the component of the vertex $1$ constructed using these edge weights.
Fix $p\in (0, 1)$, and let $\cmnd_p$  the graph with vertex set $[n]$ and edge set 
$\big\{e\in E(\cmnd)\, :\, U_e\leq p \big\}$.
Let $\cC^{\vd}_p(1)$ denote the component of the vertex $1$ in $\cmnd_p$, and $M^{\vd}_p$ be the restriction of $M^{\vd}$ to $\cC^{\vd}_p(1)$.
For each vertex $v$ in $\cC^{\vd}_p(1)$, let $d_{v,p}$ denote the degree of $v$ in $\cC^{\vd}_p(1)$, and define
$d_{v,p}^{\avail}:=d_v-d_{v,p}$.
Note that $M^{\vd}$ can be viewed as $M^{\vd}_p$ together with a collection of trees each of which is attached to a vertex of $M^{\vd}_p$ via an edge; 
for every $v\in V(M^{\vd}_p)$, let $\fT_{v, p}^{(i)}$, $1\leq i\leq r_{v,p}$, denote the trees that are attached to $v$ (arranged following some deterministic rule).
For every $v\in V(M^{\vd}_p)$, append $(d_{v,p}^{\avail}-r_{v,p})$ many zeros to the sequence $\big(\big|\fT_{v, p}^{(i)}\big|,~1\leq i\leq r_{v,p}\big)$ and let $\big(\Delta_{v, p}^{(i)},~1\leq i\leq d_{v,p}^{\avail}\big)$ be a uniform permutation of the resulting sequence; use independent permutations for different $v\in V(M^{\vd}_p)$ that are also independent of all the other random variables being considered.
Then conditional on $\cmnd_p$ and $M^{\vd}_p$, the family
\begin{align}\label{eqn:487}
	\big(\Delta_{v, p}^{(i)}\, ;\, 1\leq i\leq d_{v, p}^{\avail}, v\in \cC^{\vd}_p(1)\big)\
	\text{ of random variables is exchangeable}\,.
\end{align}
The proof of \eqref{eqn:487} is similar to the proof of \cite[Display (7.4)]{addarioberry-sen}.

Now, using \cite[Corollary 2.12]{janson-equiv}, it is enough to prove \eqref{eqn:666-a} for an IRG with an equivalent kernel.
Let $\widehat{G}_n$ be the random graph on $[n]$ obtained by placing an edge between $i$ and $j$ independently for each $i< j\in [n]$ with probability 
\[
\widehat q_{ij}:=
\frac{w_i w_j}{w_i w_j+\sum_{k=1}^n w_k}\, .
\] 
This is referred to as the Britton-Deijfen-Martin-L\"{o}f model \cite{britton-deijfen-martinlof}.
This model has the following nice property:
For $v\in [n]$, let $D_v$ denote the degree of $v$ of $\widehat G_n$, and define 
$\mvD_n:=\big(D_1,\ldots, D_n\big)$. 
Then
\begin{align}\label{eqn:6668}
\big(\widehat G_n\, \big|\, \mvD_n=\vd\big)
\ \equald\
\cgnd
\ \equald\  
\big(\cmnd\ \big|\ \cmnd\text{ is simple}\big)
\, .
\end{align}
See \cite[Theorem 7.18]{Hofs17} for a proof of \eqref{eqn:6668}.
For any degree sequence $\vd$, let $g(\vd):=\pr\big(\cmnd\text{ is simple}\big)$.
Then using \cite[Theorem 1.4]{janson-probability-of-being-simple}, it can be shown that there exists $c>0$ such that 
\begin{align}\label{eqn:486-c}
\lim_{n\to\infty}\ \pr\big(g(\mvD_n)\geq c\big)=1\, .
\end{align}
Define $\widehat{M}^n$, $\widehat G_n(\lambda)$, $\widehat{M}^n_{\lambda}$, and $\widehat Y^n_{\lambda}$ in a manner analogous to $\tbar M^n$, $\barGn(\lambda)$, $\tbar M^n_{\lambda}$, and $\tbar Y^n_{\lambda}$ respectively.
Once again, we can use \cite[Corollary 2.12]{janson-equiv} to transfer the results proved for the models considered earlier in this paper to the random graph $\widehat G_n$.
Then \eqref{eqn:486} and Theorem \ref{thm:mst-crit-gh} imply that for any $\eps>0$,
\begin{align}\label{eqn:486-a}
\lim_{\lambda\to\infty}\limsup_{n\to\infty}\, 
\bigg( 
\pr\big( \widehat Y^n_{\lambda}>\eps n\big)
+
\pr\big(d_{\rH}\big(\widehat M^n,\, \widehat M^n_{\lambda}\big)>\eps n^{\eta}\big)
\bigg)
=0\, .
\end{align}
Further, by \eqref{eqn:64}, Lemma \ref{lem:cycle-breaking-gives-mst}, and \cite[Theorem 3.3]{AddBroGolMie13}, for every $\lambda\geq 0$,
\begin{align}\label{eqn:486-b}
\big(n^{-\eta}\cdot\widehat M^n_{\lambda}\, ,\ \widehat\mu^{n, \mvw}_{\lambda}\big)
\weakc
\mst^{\mvtheta^*}_{\lambda}
\ \ \text{ w.r.t. the GHP topology}
\end{align}
as $n\to\infty$, where $\widehat\mu^{n, \mvw}_{\lambda}$ is the probability measure that assigns mass proportional to $w_v$ to each vertex $v$ in $\widehat M^n_{\lambda}$, and
$\mst^{\mvtheta^*}_{\lambda}
=
\cb^{\infty}\big(\pzS_{\lambda}^{\mvtheta^*}\big)$
is endowed with the measure inherited from $\pzS_{\lambda}^{\mvtheta^*}$.

For $v\in [n]$, let $D_{v,\lambda}$ denote the degree of $v$ in $\widehat G_n(\lambda)$, and let $D_{v, \lambda}^{\avail}:= D_v-D_{v,\lambda}$.
Using \eqref{eqn:486-b}, it is not too difficult to show that 
\begin{align*}
	\big(n^{-\eta}\cdot\widehat M^n_{\lambda}\, ,\ \widehat\mu^{n, \avail}_{\lambda}\big)
	\weakc
	\mst^{\mvtheta^*}_{\lambda}
	\ \ \text{ w.r.t. the GHP topology}
\end{align*}
as $n\to\infty$, where $\widehat\mu^{n, \avail}_{\lambda}$ is the probability measure that assigns mass proportional to $D_{v, \lambda}^{\avail}$ to each vertex $v$ in $\widehat M^n_{\lambda}$.
Let $\widehat\mu^n$ denote the uniform probability measure on the vertices of $\widehat M^n$.
Then using \eqref{eqn:486-a}, \eqref{eqn:486-b}, \eqref{eqn:6668}, \eqref{eqn:487}, \eqref{eqn:486-c}, and \cite[Lemma 7.5]{bhamidi-sen}, it can be shown that for all $\eps>0$,
\[
\lim_{\lambda\to\infty}\limsup_{n\to\infty}\, 
\pr\bigg( d_{\GHP}\bigg( 
\big(n^{-\eta}\cdot\widehat M^n\, ,\ \widehat\mu^n\big)\, ,\
\big(n^{-\eta}\cdot\widehat M^n_{\lambda}\, ,\ \widehat\mu^{n, \avail}_{\lambda}\big)
\bigg)>\eps \bigg)
=0\, .
\]
The rest is routine.

\appendix

\section{}\label{sec:appendix}

\subsection{Proof of Lemma \ref{lem:5}} \label{sec:proof-lem-5}
By Assumption \ref{ass:wts}, for all $n\geq 4$ and $2\leq j\leq n/2$,
\begin{align*}
\frac{\sum_{i=2}^j v_i}{\sum_{i=2}^n v_i}
=
\frac{1}{\ell_n}\sum_{i=2}^j w_i
\asymp
\frac{1}{n}\sum_{i=2}^j\left(\frac{n}{i}\right)^{\alpha}
\asymp
\left(\frac{j}{n}\right)^{1-\alpha}
\asymp
\frac{1}{v_j^{\tau-2}}\, ,
\end{align*}
and hence, $\pr\big(V_n\geq u\big)\asymp 1/u^{\tau-2}$ for all $n\geq 4$ and $u\in(v_{n/2}, v_2]$.
Further, if $v_{n/2}\geq 1$, then for $u\in [1, v_{n/2}]$,
\begin{align*}
\frac{v_{n/2}^{\tau-2}}{u^{\tau-2}}
\geq
\pr\big(V_n\geq u\big)
\geq
\pr\big(V_n\geq v_{n/2}\big)
\geq
\frac{C}{v_{n/2}^{\tau-2}}
\geq
\frac{C}{\big(A_2\cdot 2^{\alpha}\big)^{\tau-2}}
\geq
\frac{C}{\big(2^{\alpha} A_2 u\big)^{\tau-2}}
\, .
\end{align*}
Thus, 
$\pr\big(V_n\geq u\big)\asymp 1/u^{\tau-2}$ for all $n\geq 4$ and $u\in [1, v_2]$.

Let us first prove the claimed lower bound in \eqref{eqn:2323}.
For $n\geq 4$ and $k\geq 2$, we have
\begin{align}
\pr\big(\Poi(V_n)\geq k\big) 
&= \E\bigg[\int_0^{V_n} e^{-u} \frac{u^{k-1}}{(k-1)!} du\bigg] \notag \\
&= \int_0^{v_2} e^{-u} \frac{u^{k-1}}{(k-1)!} \pr(V_n\geq u) du 
\geq 
\int_{1}^{v_2} e^{-u} \frac{u^{k-1}}{(k-1)!} \pr(V_n\geq u) du \notag \\
& 
\geq C \int_{1}^{v_2} e^{-u} \frac{u^{k-1}}{(k-1)!} \frac{1}{u^{\tau-2}} du 
=:
C\big(\fT_1 - \fT_2-\fT_3\big) \, , \label{eqn:15} 
\end{align}
where $\fT_1, \fT_2$, and $\fT_3$ are respectively the integrals $\int_{0}^{\infty}$, $\int_{0}^{1}$, and $\int_{v_2}^{\infty}$ of the integrand in the penultimate step.
Then for $k\geq 2$,
\begin{equation}
\label{eqn:18}
\fT_1
=
\frac{\Gamma(k-\tau+2)}{\Gamma(k)}\asymp \frac{k^{k-\tau+3/2}}{k^{k-1/2}} = \frac{1}{k^{\tau -2}}\, ,
\end{equation}
and
\begin{equation}
\label{eqn:17}
\fT_2
= 
\int_0^{1} e^{-u} \frac{u^{k-1}}{(k-1)!} \frac{1}{u^{\tau-2}} du
\leq \frac{1}{(k-1)!\times (k-\tau+2)}\, .
\end{equation}
Next, note that for any $\delta>0$, $u\mapsto u^{\delta}e^{-u/2}$ is decreasing on $[2\delta, \infty)$, and consequently, for $3\leq k\leq v_2/2$, 
\begin{align*}
\fT_3 
\leq 
\frac{e^{-v_2/2} v_2^{k-\tau+1}}{(k-1)!} \int_{v_2}^\infty e^{-u/2} du = 2\frac{e^{-v_2} v_2^{k-\tau+1}}{(k-1)!} \asymp \frac{e^{-v_2}(ev_2)^k}{k^k}\cdot \frac{\sqrt{k}}{v_2^{\tau -1}}\, ,
\end{align*}
where the last step follows from Stirling's approximation for $(k-1)!$. 
Using the fact that
$\sup_{x\geq 1} (ev_2)^x/x^x=e^{v_2}$,
we get, for $3\leq k\leq v_2/2$,
\begin{equation}
\label{eqn:20}
\fT_3\leq C\sqrt{k}/v_2^{\tau -1}. 
\end{equation}
From \eqref{eqn:15}, \eqref{eqn:18}, \eqref{eqn:17}, and \eqref{eqn:20}, we see that 
there exists $k_0\geq 3$ such that for all large $n$ and for $k_0\leq k\leq v_2/2$,
\[
\pr\big(\Poi(V_n)\geq k\big) 
\geq 
\frac{C}{k^{\tau-2}} - \frac{C'\sqrt{k}}{v_2^{\tau-1}} \geq \frac{C}{2k^{\tau -2}}\, ,
\]
where we have used the scaling asymptotics of $v_2$ from Assumption \ref{ass:wts}. 
Choosing a smaller constant $C$ if necessary, the above bound can be extended to $k\in \big\{1,\ldots, k_0\big\}$ and to all $n\geq 1$.

Now we prove the claimed upper bound in \eqref{eqn:2323}. 
Proceeding as in \eqref{eqn:15}, we see that for any $n\geq 4$ and $k\geq 2$,
\begin{align*}
\pr\big(\Poi(V_n)\geq k\big) 
&= 
\int_0^{v_2} e^{-u} \frac{u^{k-1}}{(k-1)!} \pr(V_n\geq u) du \\
&\leq
\int_0^1 \frac{u^{k-1}}{(k-1)!} du
+
\int_{1}^{v_2} e^{-u} \frac{u^{k-1}}{(k-1)!} \pr(V_n\geq u) du \\
&\leq
\frac{1}{k!}
+
C\int_{1}^{v_2} e^{-u} \frac{u^{k-1}}{(k-1)!} \frac{1}{u^{\tau-2}} du 
\leq
\frac{1}{k!}
+
C\fT_1
\leq
\frac{C'}{k^{\tau-2}}\, ,
\end{align*}
where the last step uses \eqref{eqn:18}.
This completes the proof.
\qed

\subsection{Proof of Lemma \ref{lem:ss-lem-12}}\label{sec:proof-ss-lem-12}
We can construct $\barGn(\gamma_n)$ in the following steps:
\begin{enumeratea}
\item\label{it:aa}
Generate $\tbar\cC_1^n(\gamma_n/2)$.
\item\label{it:bb}
Place edges independently between $i, j\in [n]\setminus V\big(\tbar\cC_1^n(\gamma_n/2)\big)$, $i\neq j$, with respective probabilities $p_{\gamma_n}^n \tbar q_{ij}$, where $\tbar q_{ij}$ is as in \eqref{eqn:300}.
\item\label{it:cc}
Place edges independently between $i\in V\big(\tbar\cC_1^n(\gamma_n/2)\big)$ and $j\in [n]\setminus V\big(\tbar\cC_1^n(\gamma_n/2)\big)$ with respective probabilities
\begin{align*}
\frac{\big(p_{\gamma_n}^n-p_{\gamma_n/2}^n\big)\tbar q_{ij}}{1-p_{\gamma_n/2}^n\tbar q_{ij}}
\geq
\frac{\eps_{\ref{eqn:333}} w_i w_j}{2\nu_n \ell_n}\, ,
\end{align*}
where the last step holds for all large $n$.
\end{enumeratea}
Since $1-e^{-u}\leq u\wedge 1$ for all $u\geq 0$, $\cW\big(\cC_1^n(\gamma_n/2)\big)\stod\cW\big(\tbar\cC_1^n(\gamma_n/2)\big)$.
Further, $\eps_{\ref{eqn:333}}<\eps_{\ref{eqn:ss-7a}}$.
Hence, an application of Proposition \ref{prop:lem-ss-7} and \eqref{eqn:ss-8} shows that 
\begin{align}\label{eqn:708}
\pr\big(\cW\big(\tbar\cC_1^n(\gamma_n/2)\big)\geq C_{\ref{eqn:708}}n\big)
\geq 
1-\exp\big(-Cn^{\alpha}\big)
\end{align}
for all large $n$.
Suppose for some $K_{\star}>0$, 
$\cW(\cC_{\star})>K_{\star}\log n$ for some component $\cC_{\star}$ of $\barGn(\gamma_n)$, and $\cC_{\star}\neq\tbar\cC_1^n(\gamma_n)$.
Then the component $\cC_{\star}$ will be constructed in step \eqref{it:bb} above, and there will be no edges between $\cC_{\star}$ and $\tbar\cC_1^n(\gamma_n/2)$ in step \eqref{it:cc}.
Now, conditional on steps \eqref{it:aa} and \eqref{it:bb}, on the event $\big\{\cW\big(\tbar\cC_1^n(\gamma_n/2)\big)\geq C_{\ref{eqn:708}}n\big\}$,
the expected number of edges between $\cC_{\star}$ and $\tbar\cC_1^n(\gamma_n/2)$ in step \eqref{it:cc} is
\begin{align}\label{eqn:707}
\geq 
\frac{\eps_{\ref{eqn:333}}}{2\nu_n\ell_n}\times C_{\ref{eqn:708}}n\times K_{\star}\log n
\geq
C_{\ref{eqn:707}}K_{\star}\log n
\end{align}
for all large $n$.
Now the proof can be completed by combining \eqref{eqn:708} with an application of Bennett's inequality \cite{boucheron2013concentration}.

\section*{Acknowledgments}
We thank two anonymous reviewers for their helpful comments and suggestions on a previous version of the paper.
SB was partially supported by NSF grants DMS-1613072, DMS-1606839 and ARO grant W911NF-17-1-0010. 
SS was partially supported by the Infosys Foundation, Bangalore, MATRICS grant MTR/2019/000745 from SERB, and the DST FIST program - 2021 [TPN–700661].

\bibliographystyle{plain}
\bibliography{mst_heavy}


\end{document}